\documentclass[12pt,letterpaper]{article}
\PassOptionsToPackage{unicode}{hyperref}
\PassOptionsToPackage{hyphens}{url}
\PassOptionsToPackage{dvipsnames,svgnames,x11names}{xcolor}
\usepackage{amsmath,amssymb}
\usepackage{iftex}
\ifPDFTeX
  \usepackage[T1]{fontenc}
  \usepackage[utf8]{inputenc}
  \usepackage{textcomp} % provide euro and other symbols
\else % if luatex or xetex
  \usepackage{unicode-math}
  \defaultfontfeatures{Scale=MatchLowercase}
  \defaultfontfeatures[\rmfamily]{Ligatures=TeX,Scale=1}
\fi
\usepackage{lmodern}
\ifPDFTeX\else  
\fi
\IfFileExists{upquote.sty}{\usepackage{upquote}}{}
\IfFileExists{microtype.sty}{% use microtype if available
  \usepackage[]{microtype}
  \UseMicrotypeSet[protrusion]{basicmath} % disable protrusion for tt fonts
}{}
\usepackage{xcolor}
\makeatletter
\ifx\paragraph\undefined\else
  \let\oldparagraph\paragraph
  \renewcommand{\paragraph}{
    \@ifstar
      \xxxParagraphStar
      \xxxParagraphNoStar
  }
  \newcommand{\xxxParagraphStar}[1]{\oldparagraph*{#1}\mbox{}}
  \newcommand{\xxxParagraphNoStar}[1]{\oldparagraph{#1}\mbox{}}
\fi
\ifx\subparagraph\undefined\else
  \let\oldsubparagraph\subparagraph
  \renewcommand{\subparagraph}{
    \@ifstar
      \xxxSubParagraphStar
      \xxxSubParagraphNoStar
  }
  \newcommand{\xxxSubParagraphStar}[1]{\oldsubparagraph*{#1}\mbox{}}
  \newcommand{\xxxSubParagraphNoStar}[1]{\oldsubparagraph{#1}\mbox{}}
\fi
\makeatother

\usepackage{longtable,booktabs,array}
\usepackage{calc} % for calculating minipage widths
\usepackage{etoolbox}
\makeatletter
\patchcmd\longtable{\par}{\if@noskipsec\mbox{}\fi\par}{}{}
\makeatother
\IfFileExists{footnotehyper.sty}{\usepackage{footnotehyper}}{\usepackage{footnote}}
\makesavenoteenv{longtable}
\usepackage{graphicx}
\makeatletter
\def\maxwidth{\ifdim\Gin@nat@width>\linewidth\linewidth\else\Gin@nat@width\fi}
\def\maxheight{\ifdim\Gin@nat@height>\textheight\textheight\else\Gin@nat@height\fi}
\makeatother
\setkeys{Gin}{width=\maxwidth,height=\maxheight,keepaspectratio}
\makeatletter
\def\fps@figure{htbp}
\makeatother

\makeatletter
\@ifpackageloaded{caption}{}{\usepackage{caption}}
\AtBeginDocument{%
\ifdefined\contentsname
  \renewcommand*\contentsname{Table of contents}
\else
  \newcommand\contentsname{Table of contents}
\fi
\ifdefined\listfigurename
  \renewcommand*\listfigurename{List of Figures}
\else
  \newcommand\listfigurename{List of Figures}
\fi
\ifdefined\listtablename
  \renewcommand*\listtablename{List of Tables}
\else
  \newcommand\listtablename{List of Tables}
\fi
\ifdefined\figurename
  \renewcommand*\figurename{Figure}
\else
  \newcommand\figurename{Figure}
\fi
\ifdefined\tablename
  \renewcommand*\tablename{Table}
\else
  \newcommand\tablename{Table}
\fi
}
\@ifpackageloaded{float}{}{\usepackage{float}}
\floatstyle{ruled}
\@ifundefined{c@chapter}{\newfloat{codelisting}{h}{lop}}{\newfloat{codelisting}{h}{lop}[chapter]}
\floatname{codelisting}{Listing}

\makeatother
\makeatletter
\@ifpackageloaded{caption}{}{\usepackage{caption}}
\@ifpackageloaded{subcaption}{}{\usepackage{subcaption}}
\makeatother

\ifLuaTeX
  \usepackage{selnolig}  % disable illegal ligatures
\fi
\usepackage[]{natbib}
\usepackage{bookmark}

\IfFileExists{xurl.sty}{\usepackage{xurl}}{} % add URL line breaks if available
\hypersetup{
  colorlinks=true, linkcolor=blue, filecolor=Maroon,
  citecolor=Blue, urlcolor=Blue, pdfcreator={LaTeX}}
\usepackage{amsthm,multirow,verbatim,tikz}
\usepackage[noend]{algpseudocode}
\usepackage[section]{placeins}
\newfloat{algorithm}{tbp}{loa}
\floatname{algorithm}{Algorithm}
\algnewcommand\algorithmicendprocedure{\textbf{end procedure}}
\algdef{SE}[PROC]{Procedure}{EndProcedure}[2]{\algorithmicprocedure\ #1(#2)}{\algorithmicendprocedure}
\newtheorem{theorem}{Theorem}
\newtheorem{assumption}{Assumption}
\newtheorem{lemma}{Lemma}
\newtheorem{proposition}{Proposition}
\newtheorem{corollary}{Corollary}
\theoremstyle{definition}
\newtheorem{definition}{Definition}
\theoremstyle{remark}

\newcolumntype{L}[1]{>{\raggedright\arraybackslash}p{#1}}

\newcommand{\anon}{1}
\AtBeginDocument{\typeout{JCGS-LAYOUT: width=\the\textwidth; height=\the\textheight; oddsidemargin=\the\oddsidemargin; topmargin=\the\topmargin}}

\hypersetup{pdftitle={Scaffold-Constrained Subset Dynamic Programming for Exact SSE Clustering},
pdfauthor={Yordan P. Raykov; Max A. Little},
pdfkeywords={Bayesian mixture models; Connected partitions; Nearest-neighbour graphs},
hidelinks}
\if0\anon\hypersetup{pdfauthor={}}\fi
\begin{document}
\def\spacingset#1{\renewcommand{\baselinestretch}{#1}\small\normalsize}
\spacingset{1}
\if1\anon
{
 \title{\bf Scaffold-Constrained Subset Dynamic Programming for Exact SSE Clustering}
 \author{Yordan P. Raykov\hspace{.2cm}\\
  School of Mathematical Sciences, University of Nottingham\\
  and \\
  Max A. Little \\
  School of Computer Science, University of Birmingham}
 \date{}
 \maketitle
} \fi

\if0\anon
{
 \bigskip
 \bigskip
 \bigskip
 \begin{center}
  {\LARGE\bf Scaffold-Constrained Subset Dynamic Programming for Exact SSE Clustering}
\end{center}
 \medskip
} \fi

\bigskip
\begin{abstract}
Exact Euclidean \(K\)-means partitions \(n\) observations into \(K\) unlabelled clusters, but the unrestricted search is generally exponential. We use data-derived geometric graphs to precondition an exact subset dynamic program: as a result only connected vertex subsets are admitted as clusters, while sum-of-squared-errors (SSE) loss is unchanged. A remaining-set recurrence minimises fixed-\(K\) or penalised SSE, with exact factorisation over the connected components of each remaining set. The central question we study is how much computational support can be removed while preserving an unrestricted optimum. Graph inclusion gives monotone coverage and support relations, and a bottleneck threshold identifies the first covering graph in a nested hierarchy. For fixed \(K\) and dimension, under compact ball support and density bounds, retaining \(q=O(\log n)\) nearest neighbours per observation preserves an empirical SSE optimum with probability tending to one, using an \(O(\log n/n)\) fraction of complete-graph edges. Truncated Gaussian mixtures with unequal weights and covariances satisfy these conditions. The rate we provide is a sufficient upper bound rather than a result implying polynomial optimisation complexity. Objective-matched synthetic and full-data comparisons assess coverage, compression, and reference-label agreement. As a secondary application, we illustrate how the proposed scaffold preconditioning can
be utilized to improve the efficiency of split-merge proposals that preserve unrestricted
mixture posteriors.
\end{abstract}

\noindent{\it Keywords:} Bayesian mixture models; Connected partitions;
Nearest-neighbour graphs

\vfill
\newpage

\spacingset{1.}

\section{Introduction}
\label{sec:intro}

Model-based clustering and finite-mixture estimation are fundamentally
problems over unlabeled set partitions. Relabelling clusters changes neither
their statistical content nor objectives such as within-cluster sum of
squared errors (SSE). Even in low dimensions finding the global optima of the SSE criteria is of combinatorial difficulty and
exact Euclidean \(K\)-means optimization is NP-hard in general \citep{aloise2009nphard,mahajan2012planar}.

In practice, this difficulty is handled by local optimization or approximate
inference. Lloyd-type algorithms are simple and scalable but depend on their
initialization \citep{macqueen1967multivariate,lloyd1982least};
\(k\)-means++ improves initialization but still targets a local optimum
\citep{arthur2007k}. Mixture formulations instead use expectation-maximization,
variational approximations, or Markov chain Monte Carlo (MCMC)
\citep{dempster1977maximum,blei2006variational,neal2000markov,jain2004split}.
Bayesian cluster analysis makes the partition-valued perspective explicit:
point estimates and credible sets can be defined on unlabeled partitions
\citep{wade2018bayesian}. In either setting, optimizing a criterion and
recovering the scientific groups of interest are different questions.
Unequal component weights and scales, anisotropy, outliers, and overlap
illustrate why an SSE optimum need not coincide with generating mixture
labels \citep{raykov2016simple}.

One-dimensional SSE clustering provides an instructive special case where one can utilize the geometry to significantly simplify the optimization problem. After sorting the observations, an optimal fixed-\(K\) partition can be chosen with contiguous blocks.
Each candidate cluster is then an interval, and a Bellman recurrence
combines its cost with an already solved prefix problem
\citep{bellman1961approximation,bellman1969curve,wang2011ckmeans,
nielsen2014optimal}. Weighted and Bregman-type extensions use related
structure \citep{gronlund2017fast,song2020weightedkmeans}.
The gain in the 1-$d$ case comes from a geometrically justified restriction of candidate
clusters, not a replacement of the loss.

In higher dimensions, we will be restricting clusters to connected vertex sets of a
data-derived graph, which we will call a geometric scaffold. The criterion remains
ordinary SSE, at fixed \(K\) or with a penalty on the number of clusters.
A path gives interval blocks, a tree admits partitions obtained by cutting
edges, and the complete graph admits every partition. Related connectivity
restrictions occur in regionalization and spatial clustering
\citep{assuncao2006efficient,guo2008regionalization,duque2012maxp}, spatial
product-partition models \citep{page2016spatial}, and graph-partition
optimization and sampling \citep{voice2012coalition,meila2018connected}.
Here the graph can be varied with the observations, allowing its effects
on statistical estimation and computational difficulty to be studied together.

Additive clustering belongs to the established complete set-partitioning
framework \citep{yeh1986setpartition,michalak2016complete}; exact SSE
optimization has also been approached through mathematical programming
and column generation \citep{dumerle2000interior}. We formulate a
subset dynamic program (DP) whose subproblems are indexed by the observations
remaining to be partitioned. Its remaining-set recurrence selects the block
containing the least-indexed remaining observation and reuses the solution
for the set left after that block is removed. Connected-block restrictions
retain this decomposition, and disconnected remaining sets admit exact
factorization into component problems. Sections~\ref{sec:estimator-class} and \ref{sec:subset-lattice} develop the
partition formulation and its computational consequences.

The choice of scaffold determines which partitions can be considered in the optimization. Removing edges can reduce the number of candidate clusters and subproblems, but may also exclude every unrestricted SSE optimum. We therefore study how much this search can be reduced while retaining at least one such optimum. This coverage-compression trade-off depends on the arrangement of the edges, not simply on their number.
Section~\ref{sec:global-exactness} first describes when a nested sequence of scaffolds contains an unrestricted optimum, then studies this question under a sampling model. For fixed dimension and \(K\), a distribution supported on a compact ball with density bounded above and away from zero admits a sufficient neighbour count \(q=O(\log n)\): the resulting scaffold preserves an empirical SSE optimum with probability tending to one as \(n\) increases. These conditions include truncated Gaussian mixtures with unequal weights and covariances. The result concerns optimisation of SSE, without requiring agreement with the generating mixture labels. We provided a sufficient asymptotic rate rather than a sharp prescription for choosing \(q\) in a particular sample.

The empirical study examines both the computational and statistical consequences of introducing a scaffold. We assess how far the partition search can be reduced while preserving an unrestricted SSE optimum, and distinguish successful optimization of this criterion from recovery of meaningful groups. Larger-data comparisons examine the practical scope and limitations of exact scaffold-constrained fitting. A complementary application considers whether the same geometric information can improve exploration of mixture posteriors through split-merge proposals \citep{jain2004split,lee2015treeguided}. Here the scaffold guides proposed moves rather than restricting the posterior, allowing us to assess its contribution to sampling efficiency without changing the inferential target.

\section{Clustering objectives and scaffold constraints}
\label{sec:estimator-class}
\label{subsec:subset-lattice-preliminaries}

Let \(V=[n]\) index observations
\(\mathbf{X}=(\mathbf{x}_1,\ldots,\mathbf{x}_n)\in(\mathbb R^d)^n\).
Vertex indexing will remain distinct even when observation coordinates coincide. A partition
\(\ell=\{A_1,\ldots,A_m\}\) consists of nonempty, disjoint subsets whose
union is \(V\). The partition defines an equivalence relation on observations through
\(i\sim_\ell j\) if and only if \(i,j\in A\) for some \(A\in\ell\).
Conversely, the equivalence classes of any equivalence relation on \(V\)
form a unique partition. An allocation vector
\(\mathbf z\in\{1,\ldots,m\}^n\) encodes the
same relation through \(i\sim_\ell j\) if and only if \(z_i=z_j\).
Permuting its cluster labels changes neither the equivalence classes nor
which observations are grouped together. The partition, rather than its
labelled encoding, is therefore the object of estimation.

Let \(\Pi_K(V)\), \(1\le K\le n\), be the family of partitions into exactly
\(K\) nonempty blocks. Each such partition has \(K!\) allocation vectors using all
\(K\) labels, corresponding to the permutations of those labels. Assigning each of the \(n\) observations one of \(K\) labels gives \(K^n\) allocation vectors. Some leave labels unused and therefore represent partitions with fewer than \(K\) nonempty clusters.
If a specified set of \(j\) labels is excluded, \((K-j)^n\) allocations
remain. To count partitions into exactly \(K\) nonempty clusters, we first exclude allocations that leave any label unused, using the inclusion-exclusion principle. Dividing by \(K!\) then removes the repeated representations obtained by permuting cluster labels. The resulting count is the Stirling number of the second kind,
\begin{equation*}
|\Pi_K(V)|=S(n,K)
=\frac{1}{K!}\sum_{j=0}^{K}(-1)^j\binom Kj(K-j)^n.
\end{equation*}
For fixed \(K\ge2\), the \(j=0\) term dominates as \(n\to\infty\), so
\(S(n,K)\sim K^n/K!\); already \(S(n,2)=2^{n-1}-1\).
If the cluster count is not fixed, the family
\(\Pi(V)=\bigcup_{K=1}^n\Pi_K(V)\) has size
\(B_n=\sum_{K=1}^nS(n,K)\), the \(n\)-th Bell number.
Therefore, by removing arbitrary labels repeated representations are avoided but this leaves
a large collection of distinct partitions. These counts describe full
partition enumeration, not a lower bound on the cost of exact optimisation:
shared subproblems or geometric structure can reduce that cost.

For a nonempty block \(A\), define its centroid SSE by
\begin{equation*}
\phi_{\mathrm{SSE}}(A)=
\sum_{i\in A}\|\mathbf{x}_i-\boldsymbol{\mu}(A)\|^2,
\qquad
\boldsymbol{\mu}(A)=|A|^{-1}\sum_{i\in A}\mathbf{x}_i.
\end{equation*}
The partition loss is \(\sum_{A\in\ell}\phi_{\mathrm{SSE}}(A)\).

A scaffold is a finite simple undirected graph \(H=(V,E_H)\). The scaffold restricts
which subsets may occur as equivalence classes by requiring each cluster
to induce a connected subgraph, without changing its SSE.
\begin{definition}[\(H\)-realizable partition]
\label{def:Hrealizable-terminal}
A partition is \emph{\(H\)-realizable} if \(H[A]\) is connected for every
block \(A\). Write \(\mathcal F(H)\) for these partitions and
\(\mathcal C(H)=\{A\subseteq V:A\ne\emptyset,\ H[A]\text{ is connected}\}\)
for their admissible blocks.
\end{definition}

Put \(\mathcal F_K(H)=\mathcal F(H)\cap\Pi_K(V)\). For \(1\le K\le n\),
the fixed-count estimator minimizes
\begin{equation}
\mathrm{OPT}(H,K)=
\min_{\ell\in\mathcal F_K(H)}
\sum_{A\in\ell}\phi_{\mathrm{SSE}}(A),
\label{eq:fixedK-estimator}
\end{equation}
where \(\min\emptyset=+\infty\); when feasible, denote a minimizer by
\(\hat\ell_{H,K}\). Alternatively, a DP-means-style penalty
\(\lambda>0\) selects the number of clusters through
\begin{equation}
\hat\ell_{H,\lambda}\in
\arg\min_{\ell\in\mathcal F(H)}
\left\{\sum_{A\in\ell}\phi_{\mathrm{SSE}}(A)+\lambda|\ell|\right\}.
\label{eq:penalized-estimator}
\end{equation}
The complete graph \(K_n\) admits every nonempty subset as a block and
hence every partition. We call optimization \emph{unrestricted} when no
scaffold connectivity constraint is imposed; the chosen cluster count or
penalty is retained.

The scaffold therefore changes the family over which the clustering
criterion is minimised, while leaving the criterion itself unchanged.
Whether this restriction preserves an unrestricted optimum is distinct
from whether the restricted minimum can be computed efficiently.

\section[Subset dynamic programming and scaffold preconditioning]{Subset dynamic programming and scaffold preconditioning}
\label{sec:subset-lattice}
\label{sec:exact-dp-fixed-scaffold}

\subsection{Exact optimization over unlabeled partitions}
\label{subsec:unrestricted-subset-lattice}

At each stage, the \emph{remaining set} \(S\subseteq V\) consists of the
observation indices not yet assigned to selected blocks. These sets belong
to the power set \(2^V\), the collection of all subsets of \(V\). After
choosing a block \(B\subseteq S\), the new remaining set is \(S\setminus B\).
Additivity makes its clustering problem independent of the previously
selected blocks: only \(S\) and, for fixed \(K\), the number of blocks still
required matter. Subset dynamic programming uses this structure to reuse
solutions rather than enumerate complete partitions independently.

Every partition of a nonempty remaining set \(S\) has exactly one block
containing its least-indexed observation, \(a(S)=\min S\). We require this
block to be selected next and call this the \emph{first-block rule}.
The rule selects next the block containing the smallest remaining observation index, giving each unlabeled partition a unique sequence of block selections.
Write \(\Pi_k(S)\) for partitions
of \(S\) into \(k\) nonempty blocks and define
\(F(S,k)=\min_{\ell\in\Pi_k(S)}\sum_{A\in\ell}\phi_{\mathrm{SSE}}(A)\).
The boundary conditions are \(F(\emptyset,0)=0\) and
\(F(S,k)=+\infty\) when \(k=0<|S|\) or \(k>|S|\).
These correspond to \(\Pi_0(\emptyset)=\{\emptyset\}\), with the other
boundary partition families empty.

For \(S\ne\emptyset\) and \(1\le k\le |S|\), the possible first blocks are
\(\mathcal B(S,k)=\{B\subseteq S:a(S)\in B,\ |B|\le |S|-k+1\}\).
The size bound leaves an observation for each remaining block. Removing the
unique class containing \(a(S)\), or adjoining it to a partition of the
remaining set \(S\setminus B\), gives inverse operations and therefore the disjoint decomposition
\begin{equation}
\Pi_k(S)\cong
\bigsqcup_{B\in\mathcal B(S,k)}
\left(\{B\}\times\Pi_{k-1}(S\setminus B)\right).
\label{eq:unrestricted-anchored-decomposition}
\end{equation}
Minimizing the additive loss on this decomposition yields the remaining-set recurrence
\begin{equation}
F(S,k)=\min_{B\in\mathcal B(S,k)}
\{\phi_{\mathrm{SSE}}(B)+F(S\setminus B,k-1)\},
\label{eq:unrestricted-subset-recurrence}
\end{equation}
with unrestricted optimum \(\mathrm{OPT}_{\mathrm{unc}}(K)=F(V,K)\).

The disjoint union in \eqref{eq:unrestricted-anchored-decomposition}
removes redundant orderings of the clusters. For each partition of \(S\),
there is exactly one class containing \(a(S)\). Selecting this class and
repeating the rule on the remaining set gives exactly one sequence of
block choices for each unlabeled partition. Different sequences of earlier
block choices can nevertheless leave the same set \(S\) with the same
number \(k\) of clusters still required. Their minimum remaining loss is
then the same value \(F(S,k)\), which need only be calculated once and
stored for reuse, a procedure called memoization. Algorithm~\ref{alg:unrestricted-subset-dp} evaluates \eqref{eq:unrestricted-subset-recurrence}.

\begin{algorithm}[!tbp]
\caption{Unrestricted subset dynamic programming at fixed \(K\)}
\label{alg:unrestricted-subset-dp}
\begin{algorithmic}[1]
\Require observations \(\mathbf{X}\), block costs \(\phi_{\mathrm{SSE}}\), count \(K\)
\Function{$F$}{$S,k$}
 \If{\(S=\emptyset\)} \State \Return \(0\) if \(k=0\), and \(+\infty\) otherwise \EndIf
 \If{\(k=0\) or \(k>|S|\)} \State \Return \(+\infty\) \EndIf
 \If{\((S,k)\) has a stored value} \State \Return that value \EndIf
 \State \(v\gets+\infty\), \(B^\star\gets\varnothing\)
 \ForAll{\(B\subseteq S\) with \(a(S)\in B\), \(|B|\le |S|-k+1\)}
  \State \(u\gets\phi_{\mathrm{SSE}}(B)+F(S\setminus B,k-1)\)
  \If{\(u<v\)} \State \((v,B^\star)\gets(u,B)\) \EndIf
 \EndFor
 \State store \((v,B^\star)\) at \((S,k)\); \Return \(v\)
\EndFunction
\State evaluate \(F(V,K)\); recover a minimizing partition from stored choices
\end{algorithmic}
\end{algorithm}

The unrestricted recurrence remains exponential, its possible subproblems
are indexed by a subset and a remaining block count, giving
\begin{equation}
\mathcal S_K^{\mathrm{full}}=2^V\times\{0,\ldots,K\},
\qquad |\mathcal S_K^{\mathrm{full}}|=(K+1)2^n.
\label{eq:fixed-state-space-cardinality}
\end{equation}
This is an upper bound before feasibility and reachability are used. Ignoring
the count restriction, each nonempty set of size \(s\) permits
\(2^{s-1}\) blocks containing its least vertex, so
\begin{equation}
\sum_{\emptyset\ne S\subseteq V}2^{|S|-1}
=\frac12\sum_{s=1}^n\binom ns2^s
=\frac{3^n-1}{2}.
\label{eq:complete-transition-count}
\end{equation}
Thus reuse of subproblem solutions replaces direct enumeration of complete
partitions, but does not by itself make general SSE optimization tractable.

Block costs can also be reused: with
\(m_B=|B|\), \(\mathbf{s}_B=\sum_{i\in B}\mathbf{x}_i\), and
\(q_B=\sum_{i\in B}\|\mathbf{x}_i\|^2\), expansion about the sample mean gives
\begin{equation}
\phi_{\mathrm{SSE}}(B)=q_B-\|\mathbf{s}_B\|^2/m_B.
\label{eq:block-sse-moments}
\end{equation}
These moments add over disjoint sets, allowing a subset's moments to be
updated after including or removing an observation
(Appendix~\ref{app:sse-identities}).

\subsection{Scaffold preconditioning of subset dynamic programming}
\label{subsec:scaffold-subset-lattice}

A scaffold replaces arbitrary classes in the partition search by connected
classes. This can retain geometrically plausible clusters while removing
many alternatives, although preserving the unrestricted SSE optimum is a
separate coverage requirement. For remaining set \(S\), the admissible first
blocks are
\begin{equation}
\mathcal A_H(S,k)=\mathcal B(S,k)\cap\mathcal C(H)
=\{B\in\mathcal C(H):B\subseteq S,\ a(S)\in B,\ |B|\le |S|-k+1\}.
\label{eq:scaffold-transition-family}
\end{equation}
The loss and remaining-set update \(S\mapsto S\setminus B\) are unchanged, giving
\begin{equation}
F_H(S,k)=\min_{B\in\mathcal A_H(S,k)}
\{\phi_{\mathrm{SSE}}(B)+F_H(S\setminus B,k-1)\},
\label{eq:scaffold-subset-recurrence}
\end{equation}
with the boundary conditions of \(F\). Connectivity is required of each
selected block, not of the remaining set: even when \(H[S]\) is disconnected, \(S\) may still
admit partitions into several connected blocks. When \(H=K_n\),
\eqref{eq:scaffold-subset-recurrence} reduces to the unrestricted recurrence.

For the penalized estimator, let
\(\mathcal A_{H,\lambda}(S)=\{B\in\mathcal C(H):B\subseteq S,\ a(S)\in B\}\).
With \(F_{H,\lambda}(\emptyset)=0\), charge the penalty once for each
selected block:
\begin{equation}
F_{H,\lambda}(S)=\min_{B\in\mathcal A_{H,\lambda}(S)}
\{\phi_{\mathrm{SSE}}(B)+\lambda+F_{H,\lambda}(S\setminus B)\}.
\label{eq:penalized-subset-recurrence}
\end{equation}
The remaining count need not be recorded, reducing the ambient index set to
\begin{equation}
\mathcal S_\lambda^{\mathrm{full}}=2^V,
\qquad |\mathcal S_\lambda^{\mathrm{full}}|=2^n.
\label{eq:penalized-state-space-cardinality}
\end{equation}

\begin{algorithm}[!tbp]
\caption{Scaffold-constrained subset dynamic programming at fixed \(K\)}
\label{alg:scaffold-subset-dp}
\begin{algorithmic}[1]
\Require observations \(\mathbf{X}\), scaffold \(H\), count \(K\)
\State form \(\mathcal C(H)\), grouped by each block's least vertex
\Function{$F_H$}{$S,k$}
 \If{\(S=\emptyset\)} \State \Return \(0\) if \(k=0\), and \(+\infty\) otherwise \EndIf
 \If{\(k=0\) or \(k>|S|\)} \State \Return \(+\infty\) \EndIf
 \If{\((S,k)\) has a stored value} \State \Return that value \EndIf
 \If{\(k=1\)}
  \State \(v\gets\phi_{\mathrm{SSE}}(S)\) if \(S\in\mathcal C(H)\), otherwise \(+\infty\)
  \State store \(v\) and, if feasible, block \(S\); \Return \(v\)
 \EndIf
 \State \(v\gets+\infty\), \(B^\star\gets\varnothing\)
 \ForAll{\(B\in\mathcal A_H(S,k)\)}
  \State \(u\gets\phi_{\mathrm{SSE}}(B)+F_H(S\setminus B,k-1)\)
  \If{\(u<v\)} \State \((v,B^\star)\gets(u,B)\) \EndIf
 \EndFor
 \State store \((v,B^\star)\) at \((S,k)\); \Return \(v\)
\EndFunction
\State evaluate \(F_H(V,K)\)
\State if finite, recover a minimizing partition; otherwise report infeasibility
\end{algorithmic}
\end{algorithm}

Algorithm~\ref{alg:scaffold-subset-dp} displays the restriction through an
explicit collection of connected blocks. This collection need not be stored:
the same minimum can be evaluated by generating admissible blocks as needed.
The dictionary-free implementation and its block generator are detailed in
Appendix~\ref{subsec:lazy-scaffold-evaluation}.

\begin{theorem}[Exact optimization on a fixed scaffold]
\label{thm:exact-fixed-scaffold}
For any finite scaffold \(H\), \(S\subseteq V\), and integer \(k\ge0\),
let \(\mathcal P_{H,k}(S)\) be its connected
\(k\)-partitions of \(S\), including the empty partition for
\((S,k)=(\emptyset,0)\). Then
\eqref{eq:scaffold-subset-recurrence} and
\eqref{eq:penalized-subset-recurrence} satisfy
\begin{align*}
F_H(S,k)&=\min_{\ell\in\mathcal P_{H,k}(S)}
               \sum_{A\in\ell}\phi_{\mathrm{SSE}}(A),\\
F_{H,\lambda}(S)&=\min_{\ell\in\bigcup_{k=0}^{|S|}\mathcal P_{H,k}(S)}
               \sum_{A\in\ell}\{\phi_{\mathrm{SSE}}(A)+\lambda\}.
\end{align*}
In particular, \(F_H(V,K)=\mathrm{OPT}(H,K)\).
\end{theorem}

The proof applies the first-block decomposition and induction on \(|S|\)
(Appendix~\ref{app:proof-exact-scaffold-dp}). The next result explains how
a disconnected remaining set can be evaluated without treating the union of
its components as a new clustering problem. A connected block cannot join
different components of \(H[S]\);
only the allocation of the required number of clusters couples their fits.

\begin{proposition}[Component factorization of the remaining-set problem]
\label{prop:residual-component-factorization}
For \(S\ne\emptyset\), let \(S_1,\ldots,S_c\) be the components of
\(H[S]\) and define
\(\mathcal K_H(S,k)=\{\mathbf{k}\in\mathbb Z^c:
1\le k_j\le |S_j|,\ \sum_jk_j=k\}\). Restriction to these components gives
\begin{equation}
\mathcal P_{H,k}(S)\cong
\bigsqcup_{\mathbf{k}\in\mathcal K_H(S,k)}
\prod_{j=1}^c\mathcal P_{H,k_j}(S_j),
\label{eq:residual-component-bijection}
\end{equation}
and consequently
\begin{align}
F_H(S,k)&=\min_{\mathbf{k}\in\mathcal K_H(S,k)}
                    \sum_{j=1}^c F_H(S_j,k_j),
\label{eq:fixed-component-factorization}\\
F_{H,\lambda}(S)&=\sum_{j=1}^c F_{H,\lambda}(S_j),
\label{eq:penalized-component-factorization}\\
\mathcal P_{H,k}(S)\ne\emptyset
&\quad\Longleftrightarrow\quad c\le k\le |S|.
\label{eq:component-count-feasibility}
\end{align}
\end{proposition}

The feasibility inequalities illustrate that each component requires a block
and each block requires an observation; their sufficiency follows by cutting
a spanning forest (Appendix~\ref{app:proof-component-factorization}). They
also reject a candidate \(B\) when \(H[S\setminus B]\) has more than \(k-1\)
components. To recover a minimizing partition, follow the recorded block choices. Whenever the remaining graph has several connected components, use the recorded minimizing allocation of cluster counts to recover a partition within each component, then combine these partitions. Retaining every tied minimizing block choice and count allocation recovers all optimal unlabeled partitions.

\subsection{Computational support and the coverage-compression trade-off}

Removing edges can reduce both the partition family and the remaining-set
problems encountered in solving it. To make their relationship precise,
consider the constraied subset-DP in Algorithm~\ref{alg:scaffold-subset-dp},
before component factorization or additional feasibility pruning. Let
\(\mathcal R_{H,K}\) be the smallest set containing \((V,K)\) and satisfying
\begin{equation}
(S,k)\in\mathcal R_{H,K},\quad k\ge2,\quad B\in\mathcal A_H(S,k)
\quad\Longrightarrow\quad(S\setminus B,k-1)\in\mathcal R_{H,K}.
\label{eq:reachable-state-closure}
\end{equation}
Define \(C_s(H)=|\{B\in\mathcal C(H):|B|=s\}|\) and
\(C(H)=\sum_sC_s(H)\). The reached subproblems permit
\begin{equation}
T_K^{\mathrm{reach}}(H)=
\sum_{\substack{(S,k)\in\mathcal R_{H,K}\\k\ge2}}
|\mathcal A_H(S,k)|
+\sum_{(S,1)\in\mathcal R_{H,K}}\mathbf1\{S\in\mathcal C(H)\}
\label{eq:reached-transition-count}
\end{equation}
block-removal transitions, including connected terminal blocks.
Thus \(C(H)\), \(|\mathcal R_{H,K}|\), and
\(T_K^{\mathrm{reach}}(H)\) respectively count admissible blocks, remaining-set
problems, and supported transitions. Their joint monotonicity gives the
computational counterpart of optimum coverage.

\begin{proposition}[Monotonicity under scaffold inclusion]
\label{prop:scaffold-monotonicity}
For graphs on the same indexed observations with \(E_H\subseteq E_{H'}\),
\begin{align}
\mathcal C(H)&\subseteq\mathcal C(H'),&
\mathcal A_H(S,k)&\subseteq\mathcal A_{H'}(S,k),
\label{eq:monotone-block-transition}\\
\mathcal R_{H,K}&\subseteq\mathcal R_{H',K},&
\mathcal F(H)&\subseteq\mathcal F(H').
\label{eq:monotone-state-partition}
\end{align}
Hence
\begin{equation}
C_s(H)\le C_s(H'),\quad
T_K^{\mathrm{reach}}(H)\le T_K^{\mathrm{reach}}(H'),\quad
\mathrm{OPT}(H,K)\ge\mathrm{OPT}(H',K).
\label{eq:monotone-count-objective}
\end{equation}
If \(H\) contains an unrestricted optimum, so does \(H'\).
\end{proposition}

The proof follows connected-block inclusion through the recurrence
(Appendix~\ref{app:proof-scaffold-monotonicity}). In a nested family
\(H_0\subseteq\cdots\subseteq H_m\), let
\(j^\star(\mathbf X,K)=
\min\{j:\mathrm{OPT}(H_j,K)=\mathrm{OPT}_{\mathrm{unc}}(K)\}\), if it exists.
The threshold identifies the sparsest member of the nested family that preserves an unrestricted SSE optimum. Every preceding graph excludes all such optima, whereas this graph and every subsequent graph preserve at least one. Since adding edges can only enlarge the admissible block family, supported transitions and reachable states, choosing a later graph cannot reduce these support counts. The first covering graph therefore minimizes computational support within the covering portion of this family, though actual running time also depends on how that support is explored and stored.

For penalized optimization, define \(\mathcal R_{H,\lambda}\) by starting
at \(V\) and repeatedly removing any block in
\(\mathcal A_{H,\lambda}(S)\). Since these indices belong to the ambient
sets already counted in \eqref{eq:fixed-state-space-cardinality} and
\eqref{eq:penalized-state-space-cardinality},
\begin{equation}
|\mathcal R_{H,K}|\le(K+1)2^n,
\qquad |\mathcal R_{H,\lambda}|\le2^n.
\label{eq:reached-state-storage-bounds}
\end{equation}
Allocating values only when remaining-set problems are reached realizes any
reduction below these bounds; preallocating all subset entries does not.
Component factorization introduces component problems and minimizations
over count allocations, so its retained entries need not equal
\(\mathcal R_{H,K}\).

The dictionary-free evaluator stores solutions only for connected remaining sets
with \(1<k<|S|\), computes boundary values directly, and combines component
solutions without retaining a separate value for their union. Bounded
auxiliary caches reuse component decompositions and block SSEs. Discarding
a cached value may require recomputation, but removes no admissible
partition and does not change a completed minimum. Storage, repeated
evaluations, and partial block choices are therefore recorded separately.
Appendices~\ref{subsec:component-factorized-evaluation}
and~\ref{subsec:lazy-scaffold-evaluation} give the evaluation and caching
details used in the full-data comparisons.

A path ordered as \(1,\ldots,n\) illustrates the possible reduction. Every
connected first block of a remaining suffix is a prefix of that suffix, so
its removal leaves another suffix. There are at most
\((n+1)(K+1)\) reached pairs and \(n\) removals per pair, recovering
the \(O(nK)\) subproblem values and \(O(Kn^2)\) transitions of segmentation;
separately storing all candidate intervals requires \(O(n^2)\) entries.
Graph sparsity alone is insufficient: a path and a star both have \(n-1\)
edges, whereas \(C(\mathrm{path})=n(n+1)/2\) and
\(C(\mathrm{star})=2^{n-1}+n-1\). Coverage and computational support thus
depend on the arrangement of edges, not merely their number.

\subsection{The scaffold as proposal geometry for posterior partitions}
\label{subsec:scaffold-proposal-mcmc}

We also use geometric scaffolds to guide inference without restricting the statistical
model. For a mixture prior \(M\), auxiliary variables \(\eta\), and integrated
block likelihood
\(m(\mathbf{X}_A)=\int\prod_{i\in A}f(\mathbf{x}_i\mid\theta)\,H_0(d\theta)\),
the unrestricted collapsed posterior is
\begin{equation}
\pi_M(\ell,\eta\mid\mathbf{X})\propto
p_M(\ell,\eta)\prod_{A\in\ell}m(\mathbf{X}_A),
\qquad \ell\in\Pi(V).
\label{eq:unrestricted-mixture-posterior}
\end{equation}
We consider finite mixtures, Dirichlet-process (DP) mixtures, and mixtures
based on the normalized generalized gamma process (NGGP), formed by
normalizing a generalized gamma completely random measure by its total
mass \citep{favaro2013mcmc}. In these formulations, only the NGGP needs an
auxiliary variable, written \(\eta=u\) and held fixed during each partition
move; this argument is omitted for finite and DP mixtures.

A split-merge proposal selects two observations, indexed by
\(i<j\). For a split, they belong to the same block and initialize the
nonempty subclusters \(R_1=\{i\}\) and \(R_2=\{j\}\); for a merge, they
identify the two blocks to be combined. In a restricted-allocation split
\citep{jain2004split}, the baseline probability of allocating the next observation \(v\) to subcluster \(r\) is
\begin{equation}
q_{0,M}(v\to r\mid R_1,R_2,u)\propto
c_M(|R_r|;u)\frac{m(\mathbf{X}_{R_r\cup\{v\}})}{m(\mathbf{X}_{R_r})},
\label{eq:restricted-allocation-base}
\end{equation}
where \(c_M\) is the prior factor for adding an observation to an existing
block. Appendix~\ref{app:posterior-targets} details these factors.
Writing \(N_H(v)\) for the scaffold neighbours and \(d_H(v)=|N_H(v)|\), set
\begin{align}
s_H(v,r)&=\frac{|N_H(v)\cap R_r|}{\max\{1,d_H(v)\}},
\label{eq:scaffold-neighbour-score}\\
q_{H,\kappa}(v\to r\mid R_1,R_2,u)&=
\frac{q_{0,M}(v\to r\mid R_1,R_2,u)e^{\kappa s_H(v,r)}}
{\sum_{h=1}^2q_{0,M}(v\to h\mid R_1,R_2,u)e^{\kappa s_H(v,h)}}.
\label{eq:scaffold-restricted-gibbs}
\end{align}
For \(\kappa>0\), this favours allocation to the subcluster containing more
of \(v\)'s neighbours, without requiring either subcluster to be connected.
Graph distance also favours observation pairs that are far apart for splits and nearby for merges.
The generic comparator uses uniform pair-selection weights, index-ordered visits,
and \(\kappa=0\), removing all graph dependence.

Denote the move type, selected pair, and sequential allocations collectively
by the auxiliary proposal variables \(\omega\). Their joint probability
\(\widetilde q_H(\omega\mid\ell,u)\) includes the move-type and pair-selection
probabilities and, for a split, every allocation probability.
The order \(i<j\) fixes the subcluster labels; visits follow a fixed rule
determined by the parent block and the selected pair.
The corresponding reverse choices, denoted by \(\bar\omega\), merge the
proposed clusters using the same pair or reconstruct the original blocks
through a split. The Metropolis-Hastings acceptance probability accounts
for both sets of choices:
\begin{equation}
\alpha_H(\ell,\ell';\omega\mid u)
=1\wedge
\frac{\pi_M(\ell',u\mid\mathbf{X})
      \widetilde q_H(\bar\omega\mid\ell',u)}
     {\pi_M(\ell,u\mid\mathbf{X})
      \widetilde q_H(\omega\mid\ell,u)}.
\label{eq:scaffold-mh-acceptance}
\end{equation}
The complete proposal probabilities and visit rule are specified in
Appendix~\ref{app:posterior-targets}.

\begin{proposition}[Posterior invariance of the scaffold-guided move]
\label{prop:scaffold-mh-invariance}
Assume a proper target in \eqref{eq:unrestricted-mixture-posterior},
\(0<m(\mathbf{X}_A)<\infty\) for every nonempty block, and normalized
move-type, pair-selection, and allocation probabilities on their eligible sets.
With the forward and reverse auxiliary proposal variables defined above and \(u\) fixed, the
Metropolis-Hastings kernel is reversible with respect to
\(\pi_M(\ell\mid u,\mathbf{X})\), omitting \(u\) for finite and DP mixtures.
Composition with allocation and auxiliary-variable kernels preserving the
same joint target remains invariant.
\end{proposition}

Appendix~\ref{app:posterior-kernel-proof} proves detailed balance by pairing
forward and reverse choices and then summing over the auxiliary proposal
variables. Irreducibility and aperiodicity
require additional properties of the full sampler. 

\section{Coverage of an unrestricted SSE optimum}
\label{sec:global-exactness}

For fixed \(K\), let \(\mathrm{OPT}_{\mathrm{unc}}(K)=
\min_{|\ell|=K}\mathrm{SSE}(\ell)\) and let \(\mathcal M_K\) be its
minimizer set. A scaffold preserves this value when at least one
member of \(\mathcal M_K\) remains feasible; it need not retain every tied
optimum. We first relate this requirement to the connections within optimal
blocks, then examine when these connections can be retained under a sampling
model. The latter argument concerns the fitted cells of an empirical
optimum, without assuming that its blocks agree with generating labels.

\subsection{Exact coverage and deterministic connectivity}

Since \(\mathcal F(H)\) is a subset of the finite unrestricted partition family,
\begin{equation}
\mathrm{OPT}(H,K)=\mathrm{OPT}_{\mathrm{unc}}(K)
\quad\Longleftrightarrow\quad
\mathcal M_K\cap\mathcal F(H)\ne\emptyset.
\label{eq:realizability-criterion}
\end{equation}
If an unrestricted minimizer is feasible on the scaffold, it attains both
values. Conversely, equality implies that a scaffold minimizer is also an
unrestricted optimum. For a feasible fitted partition \(\widetilde\ell_H\), define
\(\Delta_{\mathrm{opt}}(H,K)=\mathrm{OPT}(H,K)-\mathrm{OPT}_{\mathrm{unc}}(K)\)
and \(\Delta_{\mathrm{search}}(\widetilde\ell_H;H,K)=
\mathrm{SSE}(\widetilde\ell_H)-\mathrm{OPT}(H,K)\). Then
\begin{equation}
\mathrm{SSE}(\widetilde\ell_H)-\mathrm{OPT}_{\mathrm{unc}}(K)
=\Delta_{\mathrm{opt}}(H,K)
+\Delta_{\mathrm{search}}(\widetilde\ell_H;H,K),
\label{eq:intro-gap-decomposition}
\end{equation}
where both terms are nonnegative. The first measures the increase in minimum
SSE caused by the scaffold restriction, while the second measures the failure
to attain that restricted minimum. Exact scaffold optimization eliminates
the second term, but the first vanishes only when coverage holds. A
restarted \(K\)-means candidate whose blocks are connected on the scaffold
provides an upper bound on \(\mathrm{OPT}(H,K)\); it establishes coverage
only if independently certified optimal. By contrast, the adjusted Rand
index (ARI) and normalized mutual information (NMI) measure agreement with
reference labels \citep{hubert1985comparing,vinh2010information}, so neither
is generally monotone in SSE.

For nested scaffolds \(\mathcal H=(H_0,\ldots,H_m)\) with \(H_m=K_n\),
adding edges cannot increase the restricted minimum or reduce computational
support (Proposition~\ref{prop:scaffold-monotonicity}). Once an unrestricted
optimum becomes feasible, it remains feasible on every subsequent graph.
The first covering graph therefore has the smallest support among the
covering members of this hierarchy. To identify this graph, define an
edge's birth index and a block's connectivity threshold by
\begin{align}
\iota_{\mathcal H}(e)&=\min\{j:e\in E(H_j)\},
\label{eq:edge-birth-index}\\
\kappa_{\mathcal H}(A)&=
\min_{T\in\mathsf{ST}(A)}\max_{e\in E(T)}\iota_{\mathcal H}(e),
\qquad |A|\ge2,
\label{eq:block-birth-threshold}
\end{align}
where \(\mathsf{ST}(A)\) contains the spanning trees on \(A\), and
\(\kappa_{\mathcal H}(A)=0\) for a singleton. Under this notation
\(\kappa_{\mathcal H}(\ell)=\max_{A\in\ell}\kappa_{\mathcal H}(A)\).

\begin{proposition}[Exact edge-birth threshold for a nested scaffold]
\label{prop:nested-edge-birth-threshold}
The induced graph \(H_j[A]\) is connected exactly when
\(j\ge\kappa_{\mathcal H}(A)\). Consequently, writing
\(\kappa_{\mathcal H,K}^\star(\mathbf X)=
\min_{\ell\in\mathcal M_K}\kappa_{\mathcal H}(\ell)\),
\begin{equation}
\mathrm{OPT}(H_j,K)=\mathrm{OPT}_{\mathrm{unc}}(K)
\quad\Longleftrightarrow\quad
j\ge\kappa_{\mathcal H,K}^\star(\mathbf X).
\label{eq:optimal-birth-threshold}
\end{equation}
Thus the first covering index is
\(j^\star(\mathbf X,K)=\kappa_{\mathcal H,K}^\star(\mathbf X)\).
\end{proposition}

A block is connected as soon as the scaffold contains a spanning tree on
its observations, and a partition is feasible once this holds for every
block. When there are several unrestricted SSE optima, coverage requires
only the first of these partitions to become feasible. Taking the minimum
of their thresholds therefore gives the first covering index.
Appendix~\ref{app:proof-bottleneck-lexicographic} proves the statement and
gives a recurrence for computing both this index and the unrestricted
value. For
\(H_q=Q_q(\mathbf X)\cup T_{\mathbf X}\), where \(Q_q\) is the union
\(q\)-nearest-neighbour graph and \(T_{\mathbf X}\) one common Euclidean
minimum spanning tree (MST), edge births are zero on that tree and otherwise
the smaller of the two directed neighbour ranks.

For a specified partition \(\ell=\{A_1,\ldots,A_K\}\), we compare the distances needed to connect observations within each block with the distances separating different blocks. For each block, construct a Euclidean minimum spanning tree, choosing one by a fixed rule if it is not unique. Define
\begin{align}
r_{\mathrm{in}}(\ell)
&=\max_k\max_{e\in\operatorname{MST}(A_k)}|e|,
\label{eq:within-bottleneck}\\
r_{\mathrm{out}}(\ell)
&=\min_{\substack{k\ne l\\\mathbf x\in A_k,\ \mathbf y\in A_l}}
\|\mathbf x-\mathbf y\|.
\label{eq:between-separation}
\end{align}
Here \(|e|\) denotes the length of an edge. Thus, \(r_{\mathrm{in}}(\ell)\) is the largest edge length needed to connect the blocks internally, whereas \(r_{\mathrm{out}}(\ell)\) is the shortest distance between observations in different blocks. A singleton requires no connecting edge and contributes zero to the first maximum; for a single-block partition, we set \(r_{\mathrm{out}}(\ell)=+\infty\).

\begin{proposition}[When the global MST covers an optimal partition]
\label{prop:mst-separated-coverage}
If
\begin{equation}
r_{\mathrm{in}}(\ell)<r_{\mathrm{out}}(\ell),
\label{eq:mst-separation-condition}
\end{equation}
every Euclidean MST on \(\mathbf X\) induces a connected subgraph on each
block of \(\ell\). If also \(\ell\in\mathcal M_K\), that MST and every
scaffold containing it preserve the unrestricted SSE optimum.
\end{proposition}

This formalizes the longstanding interpretation of MST clustering as
retaining short within-group connections and removing longer connections
between groups \citep{zahn1971graph}. Let \(G_r\) join observations at
distance at most \(r\). The first scale, \(r_{\mathrm{in}}\), is sufficient
to connect each block internally; the second, \(r_{\mathrm{out}}\), is the
first at which different blocks can be joined. For
\(r_{\mathrm{in}}\le r<r_{\mathrm{out}}\), the components of \(G_r\),
equivalently the single-linkage groups at threshold \(r\), are exactly
\(\ell\). The condition permits curved groups joined by chains of short
links, without requiring all within-group distances to be small. When the
partition is an SSE optimum, this condition guarantees coverage. When it
instead describes scientific groups, it expresses a grouping criterion
that need not agree with SSE minimization. The condition is sufficient but
not necessary for MST coverage; Appendix~\ref{app:proof-mst-separated-coverage}
gives the exchange proof.

\label{subsec:qnn-coverage}
Denser graphs can connect blocks without strict between-block separation.
Define the fill radius
\(h(A)=\sup_{\mathbf y\in\operatorname{conv}(A)}
\min_{\mathbf x\in A}\|\mathbf y-\mathbf x\|\) and
\(h(\ell)=\max_{A\in\ell}h(A)\). Let \(Q_q(\mathbf X)\) be the
undirected union \(q\)-nearest-neighbour graph, including the whole tied
shell: either endpoint may retain an edge. Set \(Q_0\) empty and
\(Q_q=K_n\) for \(q\ge n-1\). The count
\(N_i(r)=\#\{j\ne i:\|\mathbf x_i-\mathbf x_j\|\le r\}\) includes
observations from all blocks; write \(M_{\mathbf X}(r)=\max_iN_i(r)\).

\begin{proposition}[Radius and nearest-neighbour coverage]
\label{prop:radius-qnn-coverage}
For every partition \(\ell\) and \(r\ge0\),
\begin{equation}
\ell\in\mathcal F(G_r)\ \Longleftrightarrow\ r\ge r_{\mathrm{in}}(\ell),
\qquad r_{\mathrm{in}}(\ell)\le2h(\ell).
\label{eq:radius-bottleneck-fill}
\end{equation}
If \(M_{\mathbf X}(r)\le q\), then \(G_r\subseteq Q_q(\mathbf X)\).
Hence \(r_{\mathrm{in}}(\ell)\le r\) and \(M_{\mathbf X}(r)\le q\)
make \(\ell\) feasible on every \(H_q\supseteq Q_q(\mathbf X)\), giving
exact coverage when \(\ell\in\mathcal M_K\).
\end{proposition}

The proof in Appendix~\ref{app:proof-radius-qnn} considers the two components
obtained by removing a longest MST edge to establish both the radius
threshold and the fill bound. The sufficient rank
\(M_{\mathbf X}\{r_{\mathrm{in}}(\ell)\}\) need not be minimal, since
retaining every radius edge is stronger than retaining enough edges to
connect each block. When the optimal blocks are unknown, this result
characterizes their required connectivity but does not provide a computable
certificate of global optimality. For planar data in general position, the minimum spanning tree,
relative-neighbourhood graph (RNG), Gabriel graph, and Delaunay graph give
the familiar hierarchy
\(\mathrm{MST}\subseteq\mathrm{RNG}\subseteq\mathrm{Gabriel}
\subseteq\mathrm{Delaunay}\subseteq K_n\)
\citep{okabe2000spatial,berg2008computational}. Delaunay preserves a
two-cluster SSE optimum in general position; Appendix~\ref{app:delaunay-coverage} states and proves this special result.

\subsection{Model-dependent exact \texorpdfstring{\(q\)}{q}-NN coverage}
\label{subsec:gmm-qnn-coverage}

Under a sampling model, the observations determine both the optimal
partition and the scaffold. Connectivity of prespecified population groups
is therefore insufficient: the cells fitted to the same observations must
contain enough points to connect their assigned blocks. This also motivates
a computational consideration in model selection, since a statistically
appropriate model may have an optimum that is difficult to certify. We
examine whether a sparse scaffold preserves that optimum, while recognizing
that the cost of finding it also depends on the resulting partition search.

Let \(\mathbf x_1,\ldots,\mathbf x_n\) be independent draws from \(P\)
on \(\mathbb R^d\), and fix \(2\le K\le n\). For centre vectors
\(\mathbf c=(\mathbf c_1,\ldots,\mathbf c_K)\), define
\[
\mathcal Q_P(\mathbf c)=\mathbb E_P\min_k\|\mathbf x-\mathbf c_k\|^2,
\qquad
\mathcal Q_n(\mathbf c)=n^{-1}\sum_i\min_k\|\mathbf x_i-\mathbf c_k\|^2.
\]
Here \(K\) need not equal the number of mixture components. Population
geometry is useful because uniform convergence of the empirical risk keeps
its minimizing centres close to the population-optimal set
\citep{pollard1981strong,levrard2015nonasymptotic}. When the population
centres are separated and lie strictly inside the ball on which the density
is bounded below, nearby fitted centres also have cells with positive
probability in local neighbourhoods. An upper density bound then controls
how many observations compete for places in the neighbour lists needed to
retain connections within these cells.

Write \(B(\mathbf z,s)\) for the closed Euclidean ball of radius \(s\)
centred at \(\mathbf z\), and fix \(R>0\). Observations outside
\(B(\mathbf0,R)\) are not discarded: their probability
\(\tau_R=P\{\|\mathbf x\|>R\}\) appears in the bound.

\begin{assumption}[Density and optimal-centre regularity]
\label{ass:qnn-core-regularity}
The distribution has finite second moment and a density satisfying
\begin{equation}
0<\underline f_R\le f(\mathbf x)\quad(\|\mathbf x\|\le R),
\qquad f(\mathbf x)\le\overline f<\infty\quad(\mathbf x\in\mathbb R^d).
\label{eq:density-sandwich}
\end{equation}
Writing \(\mathcal C_R=B(\mathbf0,R)^K\) and
\(\mathcal Q_P^\star=\inf_{\mathbf c}\mathcal Q_P(\mathbf c)\), the set
\(\mathfrak C_{P,R}=\{\mathbf c\in\mathcal C_R:
\mathcal Q_P(\mathbf c)=\mathcal Q_P^\star\}\) is nonempty. Every member
has pairwise distinct centres strictly inside the ball; uniqueness is not
required.
\end{assumption}

For a distribution supported on a compact ball, the lower density bound
implies these conditions on the optimal centres
(Appendix~\ref{app:compact-centre-regularity}). Otherwise, they must hold
for the chosen ball. Define
\(d_\infty(\mathbf c,\mathbf c')=\max_k\|\mathbf c_k-\mathbf c'_k\|\).
Write \(d_{P,R}(\mathbf c)\) for the minimum of this distance to
\(\mathfrak C_{P,R}\). The optimal set includes all label
permutations. Continuity and compactness make the following risk gap positive
for \(\epsilon>0\), with value \(+\infty\) if its constraint set is empty:
\begin{align}
g_{P,R}(\epsilon)&=
\inf_{\substack{\mathbf c\in\mathcal C_R\\d_{P,R}(\mathbf c)\ge\epsilon}}
\{\mathcal Q_P(\mathbf c)-\mathcal Q_P^\star\},
\label{eq:population-risk-gap}\\
Z_{n,R}&=\sup_{\mathbf c\in\mathcal C_R}
|\mathcal Q_n(\mathbf c)-\mathcal Q_P(\mathbf c)|.
\label{eq:uniform-distortion-error}
\end{align}
For any empirical minimizer \(\widehat{\mathbf c}\) over \(\mathcal C_R\)
and population optimum \(\mathbf c^\star\),
\[
\mathcal Q_P(\widehat{\mathbf c})-\mathcal Q_P^\star
\le\mathcal Q_n(\widehat{\mathbf c})-\mathcal Q_n(\mathbf c^\star)+2Z_{n,R}
\le2Z_{n,R}.
\]
Thus, when \(2Z_{n,R}<g_{P,R}(\epsilon)\), every fitted optimum lies within
\(\epsilon\) of a population optimum after matching the centre labels.
This localization allows the geometry of the population-optimal cells to
be used for the fitted cells. The probability that the required risk bound
fails is
\(\zeta_{n,R}(\epsilon)=\mathbb P\{2Z_{n,R}\ge g_{P,R}(\epsilon)\}\).

The smallest separation between population-optimal centres and their
smallest distance from the boundary of the ball are, respectively,
\[
\Delta_*=\min_{\mathbf c\in\mathfrak C_{P,R}}\min_{k\ne l}
\|\mathbf c_k-\mathbf c_l\|>0,
\qquad
\rho_*=\min_{\mathbf c\in\mathfrak C_{P,R}}\min_k(R-\|\mathbf c_k\|)>0.
\]
Fix \(0<\epsilon_0<\min\{\Delta_*/4,\rho_*/2\}\) and put
\(\mathcal U_0=\{\mathbf c\in\mathcal C_R:d_{P,R}(\mathbf c)<\epsilon_0\}\).
For \(\mathbf c\in\mathcal U_0\), centre separation is at least
\(\Delta_*-2\epsilon_0\), and boundary distance at least
\(\rho_*-\epsilon_0\). Let \(D_k(\mathbf c)\) assign nearest-centre ties
to the smallest index, and let
\(C_k(\mathbf c)=B(\mathbf0,R)\cap\bigcap_{l\ne k}
\{\mathbf x:\|\mathbf x-\mathbf c_k\|\le\|\mathbf x-\mathbf c_l\|\}\)
be its closed, restricted cell. The difference from
\(D_k(\mathbf c)\cap B(\mathbf0,R)\) lies on bisector hyperplanes and has
zero probability. Closed cells can therefore be used for probability bounds
while the tie rule defines the actual assignments.

A ball centred at a fitted centre remains within its Voronoi cell if its
radius is no greater than half the distance to any other centre. Requiring
the radius also to be no greater than the distance to the boundary keeps
the ball within \(B(\mathbf0,R)\). The two bounds therefore give
\begin{equation}
B(\mathbf c_k,a_0)\subseteq C_k(\mathbf c),
\qquad
a_0=\min\{(\Delta_*-2\epsilon_0)/2,\rho_*-\epsilon_0\}>0.
\label{eq:regular-codebook-radius}
\end{equation}
Since each restricted cell is convex, the line segments joining any of its
points to this interior ball remain in the cell. This gives a lower bound
on the volume of cell neighbourhoods, including those centred on the
boundary. Denote with \(v_d\) the volume of the unit ball in \(\mathbb R^d\).

\begin{lemma}[Localization and neighbourhood probability]
\label{lem:codebook-localization-thickness}
Under Assumption~\ref{ass:qnn-core-regularity}, on
\(\{2Z_{n,R}<g_{P,R}(\epsilon_0)\}\), every empirical minimizer over
\(\mathcal C_R\) belongs to \(\mathcal U_0\). For every
\(\mathbf c\in\mathcal U_0\), \(k\le K\), \(\mathbf z\in C_k(\mathbf c)\)
and \(0<s\le2R+a_0\),
\begin{equation}
\operatorname{vol}\{C_k(\mathbf c)\cap B(\mathbf z,s)\}
\ge\theta_0v_ds^d,
\qquad \theta_0=\left(\frac{a_0}{2R+a_0}\right)^d.
\label{eq:uniform-cell-thickness}
\end{equation}
Consequently, this intersection has probability at least
\(\underline f_R\theta_0v_ds^d\).
\end{lemma}

To ensure that observations occur throughout the fitted cells, collect the regions
\(B(\mathbf0,R)\cap D_k(\mathbf c)\cap B(\mathbf z,r/4)\), for
\(\mathbf c\in\mathcal U_0\) and \(\mathbf z\in C_k(\mathbf c)\), in
\(\mathcal A_r\), where \(0<r\le4(2R+a_0)\). The preceding lemma bounds
their probabilities below by \(\eta_r\), defined below. A probability
\(\eta\)-net intersects every region of mass at least \(\eta\)
\citep{haussler1987epsilon}. Applying this requirement simultaneously over
\(\mathcal A_r\) accounts for the fact that the same observations both
determine the fitted cells and must provide their within-cell connections.

A family shatters a finite point set if its intersections with that set
produce every subset. Its Vapnik-Chervonenkis (VC) dimension is the largest
size of a shattered set. Each region here intersects one variable ball with
at most \(K-1\) halfspaces and the fixed ball. Its VC dimension is bounded
by \(V_{d,K}=\lceil2K(d+1)\log_2(2eK)\rceil\). For \(n\ge v\), define
\begin{align}
\eta_r&=\underline f_R\theta_0v_d(r/4)^d,
\label{eq:moving-cell-minimum-mass}\\
\Psi_n(\eta,v)&=\min\left\{1,
2\left(\frac{2en}{v}\right)^v e^{-n\eta/4}\right\}.
\label{eq:epsilon-net-failure}
\end{align}
The probability that the sample misses some required region is at most
\(\Psi_n(\eta_r,V_{d,K})\)
(Appendix~\ref{subsec:moving-cell-epsilon-net-proof}). Otherwise, when the
fitted centres are localized and all observations lie in the ball, each
assigned block has fill radius at most \(r/4\). Proposition~\ref{prop:radius-qnn-coverage} then supplies within-block connections at radius \(r\).

To retain these connections in neighbour lists, we must also account for
nearby observations belonging to other cells. The upper density bound gives
\(\overline\mu_n(r)=(n-1)\overline f v_dr^d\) as a bound on each
conditional expected neighbour count. Choose \(0<\delta_{\mathrm g}<1\)
as the allowed probability that some required connection is lost from the
neighbour lists. Bounding the probability of an excessive neighbour count
by \(\delta_{\mathrm g}/n\) for each observation gives the sufficient rank
in the following theorem.

\begin{theorem}[Exact coverage of an empirical SSE optimum]
\label{thm:qnn-mixture-coverage}
Under Assumption~\ref{ass:qnn-core-regularity}, use the notation of
Lemma~\ref{lem:codebook-localization-thickness}, and let
\(n\ge V_{d,K}\), \(0<r\le4(2R+a_0)\), and
\(0<\delta_{\mathrm g}<1\). Set
\begin{equation}
q_n(r,\delta_{\mathrm g})=
\min\left\{n-1,
\left\lceil\overline\mu_n(r)
+\sqrt{2\overline\mu_n(r)\log(n/\delta_{\mathrm g})}
+\frac23\log(n/\delta_{\mathrm g})\right\rceil\right\}.
\label{eq:qnn-exact-rank-threshold}
\end{equation}
For every integer \(q\ge q_n(r,\delta_{\mathrm g})\) and scaffold
\(H_q\supseteq Q_q(\mathbf X)\),
\begin{equation}
\mathbb P\{\mathrm{OPT}(H_q,K)>\mathrm{OPT}_{\mathrm{unc}}(K)\}
\le\zeta_{n,R}(\epsilon_0)+n\tau_R
+\Psi_n(\eta_r,V_{d,K})+\delta_{\mathrm g}.
\label{eq:qnn-exact-coverage-bound}
\end{equation}
\end{theorem}

The four terms respectively account for failure of the centre-localization
bound, observations outside the ball, cell neighbourhoods containing no
observations, and excessive neighbour counts. The result also applies when
the neighbour graph is augmented by an MST, without requiring mixture-label
recovery or a unique empirical optimum. Increasing \(r\) reduces the bound
on empty cell neighbourhoods, but also raises the sufficient neighbour
rank. Obtaining sparse coverage therefore requires a radius large enough
to connect the blocks while retaining those connections with a small
neighbour count.
Appendix~\ref{app:proof-empirical-qnn-coverage} proves that the resulting
connected blocks form an actual unrestricted optimum, not merely a nearby
partition. Appendix~\ref{app:scaffold-population-risk}
relates the two gaps in \eqref{eq:intro-gap-decomposition} to population risk.

When the distribution is supported entirely within \(B(\mathbf 0,R)\), no observations fall outside this ball, so the tail term \(\tau_R\) vanishes. Choosing \(nr_n^d\) proportional
to \(\log n\) fills the cell neighbourhoods while keeping their expected
neighbour counts of logarithmic order.

\begin{corollary}[Sparse exact scaffolds on compact support]
\label{cor:qnn-compact-exact}
Suppose \(P\) is supported on \(B(\mathbf0,R)\) and satisfies
\eqref{eq:density-sandwich}, with fixed \(R,d,K\). There is a deterministic
sequence \(q_n=O(\log n)\) such that every
\(H_{q_n}\supseteq Q_{q_n}(\mathbf X)\) satisfies
\begin{equation}
\mathbb P\{\mathrm{OPT}(H_{q_n},K)=\mathrm{OPT}_{\mathrm{unc}}(K)\}
\longrightarrow1.
\label{eq:compact-exactness-limit}
\end{equation}
For \(H_{q_n}=Q_{q_n}(\mathbf X)\cup\operatorname{MST}(\mathbf X)\),
\(|E(H_{q_n})|/\binom n2=O(\log n/n)\) almost surely. In particular,
these conclusions hold for any finite Gaussian mixture with positive weights
and positive-definite covariance matrices, conditioned on the ball. Neither
population-optimum uniqueness nor agreement between \(K\) and the generating
component count is required.
\end{corollary}

The proof in Appendix~\ref{app:proof-compact-qnn-corollary} permits
\(q_n=O_{d,K}([1+\overline f/(\underline f_R\theta_0)]\log n)\).
Low within-cell density or narrow cells increase the radius required to
avoid gaps; high ambient density raises the rank needed to retain connections
at that radius. A small population-risk gap instead demands a larger sample
before fitted centres reliably inherit the population geometry. For Gaussian
mixtures, density bounds are explicit in the weights, means, covariances and
truncation radius (Appendix~\ref{app:proof-gaussian-coverage-bounds}), but
the risk gap and \(\theta_0\) also depend on the optimal SSE centres.
Mixture overlap alone is therefore not a monotone coverage parameter, and
the theorem gives a sufficient rate rather than a sharp numerical rank.

A vanishing edge fraction does not imply polynomial complexity for subset
dynamic programming, since the numbers of connected blocks and remaining-set
subproblems must also be considered. The compact-support conclusion does
not automatically extend to untruncated mixtures either.
For a fixed ball, \(n\tau_R\) cannot vanish; allowing \(R=R_n\) reduces
this term but can worsen the density and geometric constants. Additional
control is needed before asserting the same logarithmic rank.

\section{Experiments}
\label{sec:experiments}

For clustering, restricting the feasible partition family can reduce computation
while preserving an SSE optimum, or exclude that optimum in favour of a
different geometric structure. Our experiments examine how scaffold choice affects both the computation of a clustering solution and its statistical interpretation. Synthetic examples allow us to distinguish failure to minimize SSE from disagreement between an SSE optimum and the groups of interest. Comparisons on full UCI datasets then examine the practical scope of exact scaffold-constrained fitting. We also consider a different use of the same geometry in posterior inference, where the scaffold guides split-merge proposals without excluding partitions from the posterior. Here the question is whether geometric guidance improves acceptance and posterior estimation.

\subsection{Objective-matched clustering comparisons}
\label{sec:exp-clustering}

We generated 20 independent samples of size \(n=15\) from each of six
two-dimensional designs, with \(K=3\) for all methods: unequal component radii, unequal weights, two outlier
pairs, rotated ellipses, and heterogeneous separated or overlapping ellipses
\citep{raykov2016simple}. The outlier design retains
five generating labels (three Gaussian components and two outlier groups),
so label agreement and three-block SSE optimization have different targets.

\begin{figure}[!tbp]
\centering
\includegraphics[width=\linewidth]{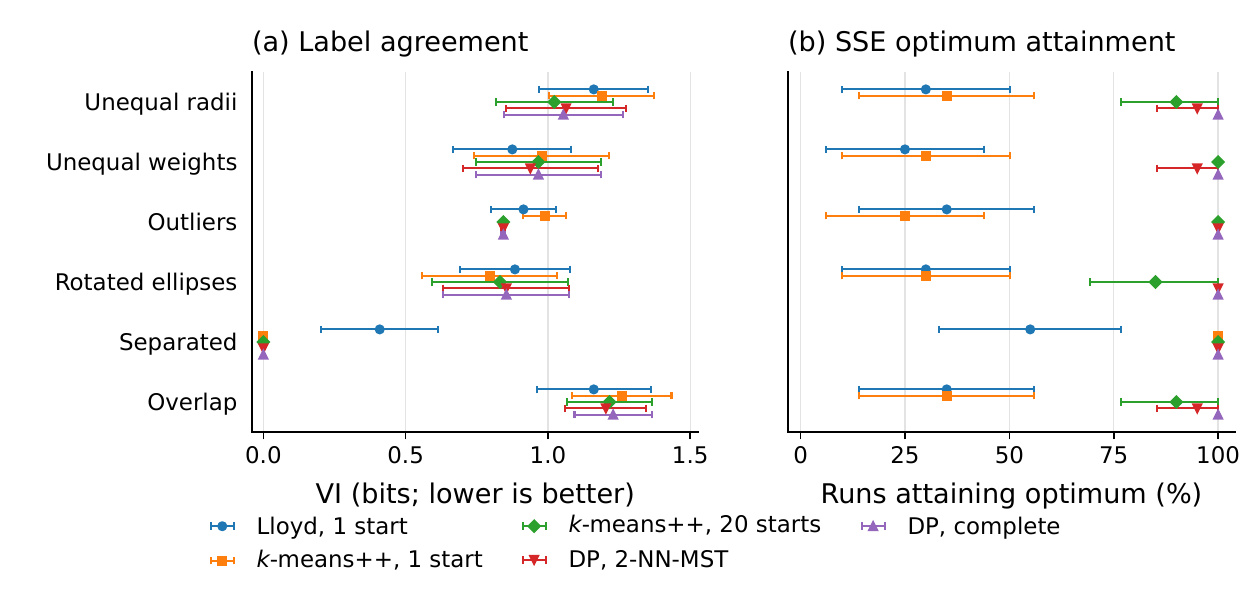}
\caption{Fixed-\(K\) comparison across the six Gaussian designs. Points are
means over 20 independently generated data sets. VI error bars are pointwise
95\% Monte Carlo intervals; oracle-hit error bars are normal-approximation
intervals for the observed proportions. The left panel compares partitions
with the generating labels by variation of information; the right panel gives
the frequency of attaining the complete-graph SSE optimum. The panels
evaluate different targets and need not rank methods identically.}
\label{fig:clustering-fixed-k}
\end{figure}

The complete-graph recurrence supplies the unrestricted optimum on each
sample. Across the 120 data sets, one-start random Lloyd and one-start
\(k\)-means++ attained it in 42 and 51 cases, with mean relative SSE gaps
of \(51.830\%\) and \(10.778\%\). Twenty-restart \(k\)-means++ attained
it in 113 cases, reducing the mean gap to negligible  \(0.065\%\). The exact
\(2\)-NN-MST recurrence attained it in 117 cases, with mean gap
\(0.032\%\); Gabriel and Delaunay agreed with the complete recurrence in
all 120 cases. Objective recovery need not improve label agreement:
mean variation of information (VI) to the generating labels was
\(0.813\) bits for restarted \(k\)-means++ and \(0.825\) bits for the
unrestricted optimum (Figure~\ref{fig:clustering-fixed-k}).

Table~\ref{tab:clustering-support} relates optimum coverage to computational
support. Counts are normalized by the corresponding complete-recurrence
count on each sample before averaging. For reached remaining sets,
\(100\%\) therefore denotes the complete
fixed-\(K\) recurrence. After its first block is removed, observation \(1\)
cannot occur in any remaining set, limiting the count to approximately half
of the \(2^n\) subsets. Edge density does not determine these reductions: Gabriel
has fewer edges than \(2\)-NN-MST but more connected blocks and supported
transitions, while reaching fewer remaining sets.
The central calculations are proportions
\begin{table}[!tbp]
\centering
\caption{Optimum coverage and computational support for the fixed-\(K\)
experiment. Coverage is the percentage of the 120 data sets for which the
scaffold recurrence equals the complete-graph SSE optimum. All remaining
entries are percentages of the corresponding complete-graph quantity.}
\label{tab:clustering-support}
\small
\begin{tabular}{lrrrrr}
\toprule
Scaffold & Coverage & Edges & Connected blocks & Reached sets
& Transitions\\
\midrule
\(2\)-NN-MST & 97.500 & 19.980 & 4.020 & 12.810 & 0.242\\
Gabriel    & 100.000 & 17.670 & 4.900 & 11.060 & 0.298\\
Delaunay    & 100.000 & 33.690 & 53.970 & 66.360 & 9.442\\
Complete    & 100.000 & 100.000 & 100.000 & 100.000 & 100.000\\
\bottomrule
\end{tabular}
\end{table}

Appendices~\ref{sec:supp-unknown-k}--\ref{sec:supp-structural-controls}
extend this comparison to penalized SSE and collapsed-MAP objectives,
component factorization, and structural controls. Factorization reproduces
objective values and partition counts in all 32 paired evaluations and
reduces block inspections, but not uniformly stored values. The structural
controls separate scaffold exclusion from local-search failure; the larger
separability sweep assesses candidate connectivity, not unrestricted
optimum coverage.

\subsection{Scaffold restriction for nonconvex groups}
\label{sec:exp-moons}

The two-moons construction consists of two interleaving curved groups whose
geometry is not well described by dispersion around a single centroid per
group. It allows us to examine a different role for the scaffold: restricting
the partition family to reflect local connectivity, rather than seeking to
preserve an unrestricted SSE optimum. The comparison illustrates why better
optimization of SSE need not recover the groups of interest, and how
scaffold choice can change that statistical interpretation.

For two moons, we fix \(n=60\), \(K=2\), and use five noise realizations at
each standard deviation \(0.020\), \(0.080\), and \(0.160\). Random-start
Lloyd, one-start and twenty-start \(k\)-means++, and exact scaffold fits
use the same coordinates; generating labels are used only for evaluation.
We compare the Euclidean MST \(T\) with the forest
\(H_{\mathrm{SL}}=T\setminus\{e_{\max}\}\), obtained by removing its
longest edge using distances alone.

This pruning rule follows the classical short-link interpretation of MST
clustering \citep{zahn1971graph}. For the generating partition
\(\ell_{\mathrm{gen}}\), suppose
\(r_{\mathrm{in}}(\ell_{\mathrm{gen}})<r_{\mathrm{out}}(\ell_{\mathrm{gen}})\),
with the scales defined in \eqref{eq:within-bottleneck} and
\eqref{eq:between-separation}. Proposition~\ref{prop:mst-separated-coverage}
then ensures that each moon induces a connected subtree of \(T\).
Each subtree is an MST on its observations, since replacing it by a shorter
tree would shorten \(T\). Its edges have length at most
\(r_{\mathrm{in}}\), whereas the unique joining edge has length at least
\(r_{\mathrm{out}}\). Removing the longest edge therefore recovers the
two moons. Short connections along a curved group can meet this condition
despite substantial dispersion about its centroid.

The intact tree admits 59 connected two-block partitions; the pruned
forest admits only its two-component partition. Consequently,
\(\mathrm{OPT}(T,2)\le\mathrm{OPT}(H_{\mathrm{SL}},2)\), but agreement
with generating labels need not follow this ordering. On the forest the
graph fixes the estimate before SSE is evaluated, reproducing single
linkage at \(K=2\); this is scaffold restriction, not an optimization
advantage over single linkage.

\begin{figure}[!tp]
\centering
\includegraphics[width=\linewidth]{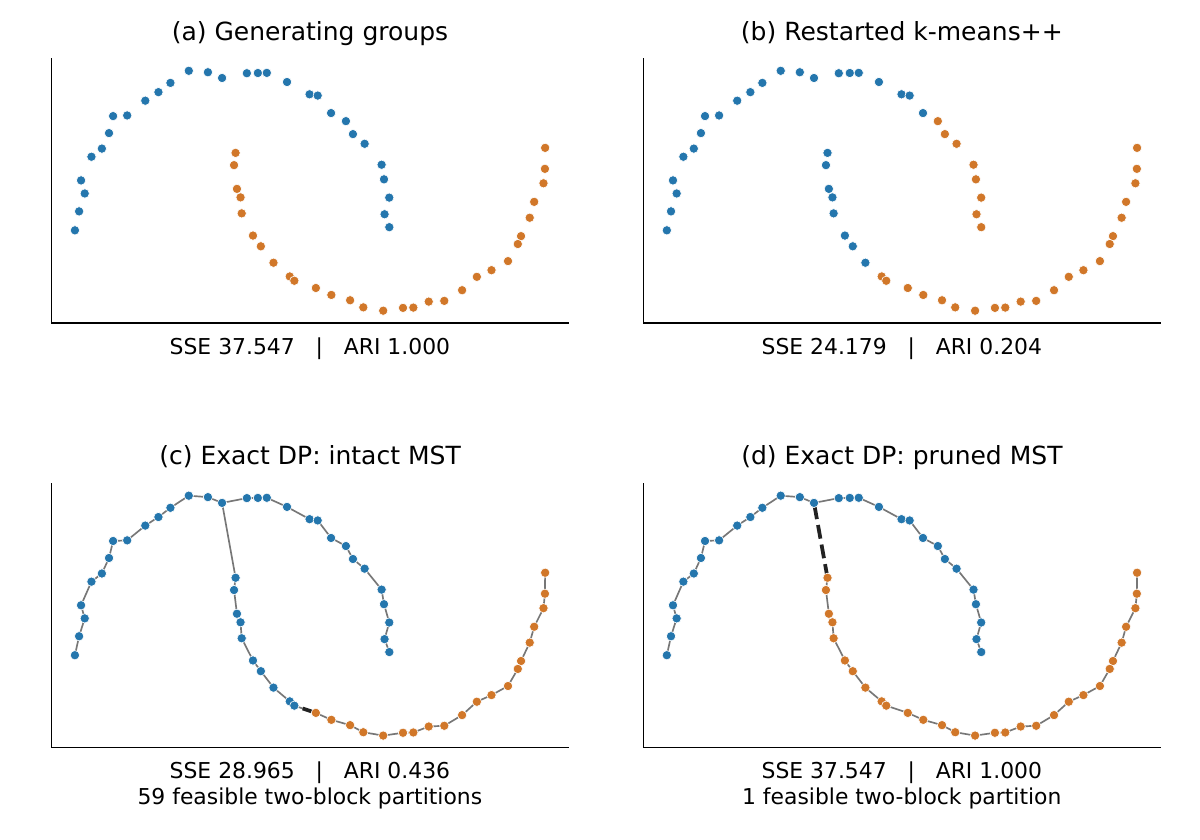}
\caption{Scaffold restriction on two moons, with \(n=60\), noise SD
\(0.020\), and the first prespecified seed. Colours show generating groups
in (a) and fitted blocks in (b), (c), and (d). Grey segments are scaffold edges.
In (c), the dashed MST edge separates the blocks selected by SSE but remains
part of the scaffold; any tree edge could have been cut. In (d), the dashed
longest edge is removed from the scaffold before optimization, forcing its
two components. Displayed scores are for this realization, not averages.}
\label{fig:moons-scaffold-choice}
\end{figure}

At noise \(0.020\), mean ARI is \(0.224\), \(0.224\), and \(0.211\)
for random Lloyd, one-start \(k\)-means++, and the restarted fit:
improved initialization does not recover the curved groups. Exact MST
optimization gives mean ARI \(0.445\) and SSE \(28.762\), compared with
SSE \(24.080\) for restarted \(k\)-means++. Although both moons are
MST-connected in every realization, they do not minimize SSE there.
Pruning recovers all five generating partitions, at mean SSE \(37.293\)
(Figure~\ref{fig:moons-scaffold-choice}). No unrestricted optimum is
computed in this comparison; exactness is conditional on the scaffold.

At noise \(0.160\), pruned-forest mean ARI falls to \(0.046\), compared
with \(0.240\) for the intact MST. Noise can increase
\(r_{\mathrm{in}}\) and decrease \(r_{\mathrm{out}}\), destroying their
separation (Appendix~\ref{sec:supp-moons}). Edge-length gaps
and threshold stability can inform pruning without labels, but removing
the longest edge always produces two components, even when these no longer
correspond to the generating groups. Exact optimization identifies the best
partition permitted by the scaffold; whether it captures the groups of
interest depends on both the SSE criterion and the graph. Here, making the
generating groups feasible on the intact MST does not ensure that SSE
selects them. Pruning can recover them by restricting the admissible
partitions, but its sensitivity to noise shows why the geometric basis
for that restriction also needs to be assessed.

\FloatBarrier
\subsection{Clustering comparisons on UCI data}
\label{sec:exp-uci}

We use five complete data sets from the University of California Irvine
(UCI) Machine Learning Repository: Iris (\(n=150\)), Wine (\(n=178\)),
Breast Cancer Wisconsin (\(n=569\)), Soybean Small (\(n=47\)), and Statlog
Vehicle (\(n=845\)) \citep{dua2019uci}. Every nonconstant feature is
centered and scaled to unit empirical variance, and \(K\) is the number
of reference classes. Random-start Lloyd, one-start and twenty-start
\(k\)-means++, and Ward agglomeration are compared with exact scaffold
recurrences. Each initialization fits all observations; the restarted
baseline retains the smallest SSE over 20 initializations.

The nested scaffolds are \(H_q=Q_q\cup T\), where \(Q_q\) is the union
\(q\)-nearest-neighbour graph, \(T\) is a common Euclidean MST, and
\(q\in\{0,2,5\}\), with \(H_0=T\). Since every data set has \(d>2\),
these graphs use the full feature space without a planar projection or
Delaunay comparison. The component-factorized recurrence generates blocks
on demand and caches only nontrivial connected subproblems. Completed
fits are exact on their scaffold; stopped searches supply no partition
(Appendix~\ref{sec:supp-uci}).

Table~\ref{tab:uci-method-comparison} reports SSE, its ratio to the one-cluster
loss \(\mathrm{TSS}=\sum_i\|\mathbf{x}_i-\bar{\mathbf{x}}\|^2\), and
reference-class agreement. ARI is chance-adjusted,
has null expectation zero, and may be negative, with range contained in
\([-1/2,1]\). Arithmetic-normalized NMI lies in \([0,1]\) without chance
adjustment. Larger values of either indicate stronger agreement, with one
denoting identical partitions up to relabelling. These indices assess a
different target from centroid loss and do not establish model adequacy.

\begingroup
\small
\setlength{\tabcolsep}{4pt}
\setlength{\LTcapwidth}{\textwidth}
\begin{longtable}{llrrrr}
\caption{Full-data UCI clustering comparison: every method uses all \(n\)
observations, with no subsampling. Starts count initializations, not data samples.
Smaller SSE and SSE/TSS indicate lower centroid loss; larger ARI and NMI
indicate stronger agreement with the reference classes.
\(H_q\) is union \(q\)-NN plus a common MST. Time limits mark
unfinished exact searches, not fitted partitions. All numerical DP entries
are exact on the stated scaffold. Computational settings and work counts
for these runs are reported in Appendix~\ref{sec:supp-uci}.}
\label{tab:uci-method-comparison}\\
\toprule
Data set & Method & SSE & SSE/TSS & ARI & NMI\\
\midrule
\endfirsthead
\caption[]{Full-data UCI clustering comparison (continued).}\\
\toprule
Data set & Method & SSE & SSE/TSS & ARI & NMI\\
\midrule
\endhead
\midrule
\multicolumn{6}{r}{\emph{Continued on next page.}}\\
\endfoot
\bottomrule
\endlastfoot
Iris ($n=150$) & Lloyd, random start & 140.902 & 0.235 & 0.645 & 0.661 \\*
$d=4,\ K=3$ & $k$-means++, 1 start & 139.820 & 0.233 & 0.620 & 0.659 \\*
 & $k$-means++, 20 starts & 139.820 & 0.233 & 0.620 & 0.659 \\*
 & Ward & 148.876 & 0.248 & 0.615 & 0.675 \\*
 & DP, MST & 142.030 & 0.237 & 0.594 & 0.647 \\*
 & DP, $H_2$ & 139.985 & 0.233 & 0.601 & 0.647 \\
\midrule
Wine ($n=178$) & Lloyd, random start & 1279.966 & 0.553 & 0.880 & 0.861 \\*
$d=13,\ K=3$ & $k$-means++, 1 start & 1277.928 & 0.552 & 0.897 & 0.876 \\*
 & $k$-means++, 20 starts & 1277.928 & 0.552 & 0.897 & 0.876 \\*
 & Ward & 1305.049 & 0.564 & 0.790 & 0.786 \\*
 & DP, MST & 1310.716 & 0.566 & 0.741 & 0.753 \\*
 & DP, $H_2$ & \multicolumn{4}{c}{Time limit} \\
\midrule
Breast cancer ($n=569$) & Lloyd, random start & 11595.683 & 0.679 & 0.677 & 0.562 \\*
$d=30,\ K=2$ & $k$-means++, 1 start & 11595.683 & 0.679 & 0.677 & 0.562 \\*
 & $k$-means++, 20 starts & 11595.461 & 0.679 & 0.671 & 0.555 \\*
 & Ward & 11866.537 & 0.695 & 0.575 & 0.457 \\*
 & DP, MST & 12983.354 & 0.761 & 0.396 & 0.405 \\*
 & DP, $H_2$ & \multicolumn{4}{c}{Time limit} \\
\midrule
Soybean ($n=47$) & Lloyd, random start & 387.251 & 0.392 & 0.748 & 0.847 \\*
$d=35,\ K=4$ & $k$-means++, 1 start & 387.251 & 0.392 & 0.748 & 0.847 \\*
 & $k$-means++, 20 starts & 367.144 & 0.372 & 1.000 & 1.000 \\*
 & Ward & 367.144 & 0.372 & 1.000 & 1.000 \\*
 & DP, MST & 367.144 & 0.372 & 1.000 & 1.000 \\*
 & DP, $H_2$ & 367.144 & 0.372 & 1.000 & 1.000 \\
\midrule
Vehicle ($n=845$) & Lloyd, random start & 6664.684 & 0.438 & 0.068 & 0.093 \\*
$d=18,\ K=4$ & $k$-means++, 1 start & 6674.433 & 0.439 & 0.055 & 0.083 \\*
 & $k$-means++, 20 starts & 5971.520 & 0.393 & 0.076 & 0.116 \\*
 & Ward & 6436.435 & 0.423 & 0.088 & 0.115 \\*
 & DP, MST & 6474.224 & 0.426 & 0.106 & 0.141 \\*
 & DP, $H_2$ & \multicolumn{4}{c}{Time limit} \\
\end{longtable}

\endgroup

On Soybean, exact MST and \(H_2\) fits reduce the single-start SSE
\(387.251\) by \(5.192\%\), to \(367.144\), with perfect label
agreement. On Vehicle, the exact MST SSE \(6474.224\) improves on random
Lloyd and one-start \(k\)-means++ by \(2.858\%\) and \(3.000\%\);
ARI \(0.106\) and NMI \(0.141\) also improve, although agreement remains
weak. On Iris, Wine, and Breast Cancer, exact MST fits instead have higher
SSE and lower ARI/NMI than both single-start baselines. A lower-loss fit
\(\widehat\pi_1\) gives
\(\mathrm{OPT}_{\mathrm{unc}}(K)\le L(\widehat\pi_1)<\mathrm{OPT}(T,K)\),
proving exclusion of every unrestricted optimum without computing one.

Enriching Iris's MST to \(H_2\) lowers exact SSE from \(142.030\) to
\(139.985\), below random Lloyd but above one-start \(k\)-means++
at \(139.821\). ARI improves slightly but remains below
both single-start fits. Restarts remove the scaffold fits' single-start
advantage on Soybean. On Vehicle, restarted \(k\)-means++ and Ward give
lower SSE (\(5971.520\) and \(6436.435\)) than the exact MST fit,
thereby also proving scaffold exclusion, but have lower ARI/NMI. Exactness
on a scaffold, unrestricted loss, and reference-label agreement thus
remain separate properties. These recorded fits do not estimate average
performance over initializations.

\subsection{Scaffold-informed posterior exploration}
\label{sec:exp-posterior}

The posterior experiment uses \(n=300\) bivariate observations from four
Gaussian components arranged as two separated pairs with overlap within
each pair. Within each finite, DP, or NGGP mixture model, collapsed Gibbs
sampling and four split-merge augmentations target the same unrestricted
posterior. We compare generic proposals, \(2\)-NN-MST and \(3\)-NN-MST
guidance from Section~\ref{subsec:scaffold-proposal-mcmc}, and a
vertex-permuted \(3\)-NN-MST control that preserves the weighted graph
but breaks its alignment with the data.

Four chains per method use matched dispersed initializations. Each augmented
sweep comprises an allocation Gibbs sweep and one split-merge attempt;
all NGGP methods use the same additional slice update for \(\log u\).
Comparisons use 500, 1,000, and 2,500 retained sweeps per chain.
A bounded cache reuses block likelihoods without enumerating connected
subsets or solving a subset DP. Appendix~\ref{sec:supp-posterior-diagnostics} specifies the protocol and
independent longer reference chains, which pass the stated diagnostics
but are not posterior oracles.

Accuracy concerns posterior recovery, not generating labels: we compare
posterior similarity matrices (PSMs) of pairwise co-clustering
probabilities, distributions of \(K\), and expected-VI decisions over a
common reference-derived candidate set \citep{wade2018bayesian}.

\begin{figure}[!tbp]
\centering
\includegraphics[width=\linewidth]{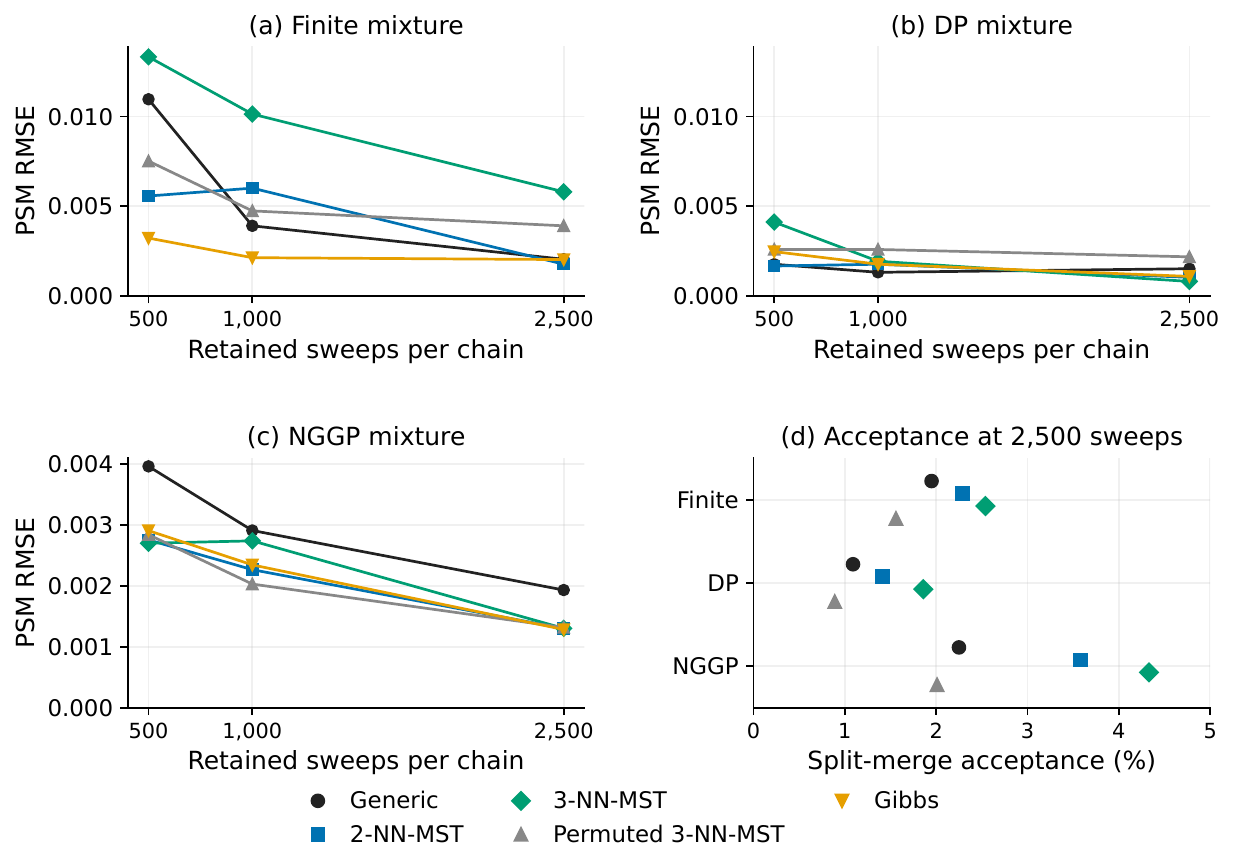}
\caption{Posterior accuracy and proposal acceptance on one shared
\(n=300\) data set. Panels (a) to (c) show off-diagonal PSM RMSE against
independent reference chains, using four evaluation chains at matched
retained-sweep counts; lower is better. The NGGP panel uses a finer vertical
scale. Panel (d) reports accepted
split-merge moves as a percentage of attempts over the 2,500-sweep prefixes,
excluding warm-up. Gibbs has no split-merge attempts.
The nested prefixes are not independent replicates; lines join observed
values and do not represent fitted convergence rates.}
\label{fig:posterior-benchmark}
\end{figure}

At 2,500 sweeps, \(3\)-NN guidance raises acceptance over generic
split-merge from \(1.950\%\) to \(2.540\%\) for the finite mixture,
\(1.090\%\) to \(1.860\%\) for DP, and \(2.250\%\) to \(4.330\%\)
for NGGP. The \(2\)-NN rates lie between these in each target; permutation
reduces acceptance below generic proposals throughout. Geometric alignment
therefore helps proposal acceptance on this sample, though rates remain modest.

Posterior accuracy shows no uniform corresponding gain. For NGGP, both
guided samplers have final PSM RMSE about \(1.300\times10^{-3}\), below
generic split-merge at \(1.931\times10^{-3}\), but similar to Gibbs and
the permuted control at \(1.283\times10^{-3}\) and
\(1.317\times10^{-3}\). The \(3\)-NN sampler has the smallest reference
error for the distribution of \(K\), whereas \(2\)-NN does not improve on generic
split-merge (Table~\ref{tab:posterior-scaling-diagnostics}).
Increasing from \(2\)-NN to \(3\)-NN reduces DP PSM RMSE from
\(1.038\times10^{-3}\) to \(8.104\times10^{-4}\), but increases finite-mixture
RMSE from \(1.771\times10^{-3}\) to \(5.798\times10^{-3}\).
At 2,500 sweeps every sampler selects the same three-cluster VI decision,
minimizing reference risk within its target-specific candidate set, not
necessarily over all partitions. The evidence supports improved acceptance
without a uniform posterior-accuracy advantage at matched iterations.

\FloatBarrier

\section{Discussion}
\label{sec:discussion}

Geometric scaffolds allow the statistical choice of admissible clusters
and the difficulty of exact optimization to be studied together. Our results
relate preservation of an unrestricted optimum to the candidate blocks and
subproblems considered by subset dynamic programming; the number of graph
edges alone does not determine this computational benefit. Restricting clusters to connected vertex sets reduces the support of the remaining-set recurrence while preserving unrestricted SSE whenever at least one minimizing partition remains feasible. Within a nested family, we identify the first scaffold that preserves this optimum and hence has the smallest computational support among the covering members. Component factorization and generation of connected blocks as needed allow this restriction to be used without first storing every candidate block. The full-data UCI comparisons demonstrate completed exact fits on some sparse scaffolds, extending the evaluation beyond samples small enough to compute an unrestricted oracle. The search can still be exponential, but admissible lower bounds could reduce it further while retaining a certificate of the restricted minimum.

The experiments also use exact optimization to distinguish failure to
minimize a criterion from a mismatch between that criterion and the groups
of interest. In the low-noise two-moons example, the generating partition is feasible on the minimum spanning tree (MST) but does not minimize its SSE. Removing the longest edge recovers the intended groups when the within-group connectivity scale is below the between-group distance; the same restriction deteriorates as noise closes this gap. Scaffold choice can therefore change which groups are selected, as well as how much computation is required. For planar two-cluster SSE in general position, Delaunay connectivity preserves an unrestricted optimum: the optimal blocks can be represented by complementary half-planes, each inducing a connected sample subgraph
(Appendix~\ref{app:delaunay-coverage}). For more than two clusters, the cells are intersections of half-planes, so the same argument
does not directly establish their connectivity. Corollary~\ref{cor:qnn-compact-exact} establishes a probabilistic link between
scaffold sparsity and exact optimum preservation for
fixed \(d,K\) and a distribution supported on a ball with density bounded
above and away from zero. A suitable \(q_n=O(\log n)\) preserves an
empirical SSE optimum with probability tending to one, using a vanishing
fraction of complete-graph edges, without requiring the optimal partition
to recover generating mixture labels. This is a sufficient asymptotic rate,
and fewer neighbours may suffice in a particular sample. Global density bounds
can be conservative when low-density regions between components do not
obstruct connectivity within optimal cells; the proof also does not exploit
MST augmentation. Accounting more closely for local cell geometry and the
connections supplied by the MST may therefore sharpen the guarantee.
Extending it to untruncated mixtures also requires control of the tail
probability without allowing the local probability bounds to deteriorate
too quickly.

In posterior inference, the scaffold instead guides proposed changes to
the partition without restricting the mixture posterior. Accounting for
the forward and reverse proposal probabilities in the Metropolis-Hastings
correction preserves this target. The comparisons for finite mixtures, Dirichlet-process mixtures, and normalized generalized gamma process mixtures show higher acceptance for geometric proposals, while their relative posterior accuracy at matched iteration counts depends on the mixture model and summary. Geometry can thus help construct moves that are more readily accepted, although this does not by itself ensure
more accurate posterior estimates.

\clearpage
\appendix
\section{Proofs for subset dynamic programming}
\label{app:subset-lattice-proofs}

\subsection{Proof of Theorem~\ref{thm:exact-fixed-scaffold}}
\label{app:proof-exact-scaffold-dp}

For \(S\subseteq V\), write
\[
\mathcal P_{H,k}(S)
:=
\{\ell\in\Pi_k(S):H[A]\text{ is connected for every }A\in\ell\}.
\]
We prove the fixed-count identity by strong induction on \(|S|\), using the
boundary conventions stated before
\eqref{eq:unrestricted-subset-recurrence}. If \(S=\emptyset\), if \(k=0\)
with \(S\ne\emptyset\), or if \(k>|S|\), the partition family and the
recurrence have respectively the same value, namely \(0\) only for
\((S,k)=(\emptyset,0)\) and \(+\infty\) otherwise. It remains to consider
\(S\ne\emptyset\) and \(1\le k\le |S|\).

Every \(\ell\in\mathcal P_{H,k}(S)\) has a unique block \(B_\ell\) containing
\(a(S)\). It belongs to \(\mathcal A_H(S,k)\): it is connected by
realizability, contains \(a(S)\), and has at most \(|S|-k+1\) vertices
because the other \(k-1\) blocks are nonempty. Define
\[
\Psi_S(\ell)
=
\bigl(B_\ell,\ell\setminus\{B_\ell\}\bigr).
\]
Its codomain is
\[
\bigsqcup_{B\in\mathcal A_H(S,k)}
\bigl(\{B\}\times\mathcal P_{H,k-1}(S\setminus B)\bigr).
\]
Conversely, if \(B\in\mathcal A_H(S,k)\) and
\(\ell'\in\mathcal P_{H,k-1}(S\setminus B)\), then
\(\{B\}\cup\ell'\) is an \(H\)-realizable \(k\)-partition of \(S\), and its
first block is \(B\). Thus
\[
\Psi_S^{-1}(B,\ell')=\{B\}\cup\ell',
\qquad
\Psi_S^{-1}\circ\Psi_S=\operatorname{id},
\qquad
\Psi_S\circ\Psi_S^{-1}=\operatorname{id}.
\]
The decomposition is therefore a genuine disjoint-union bijection, not only a
surjection. Block additivity and the induction hypothesis now give
\begin{align*}
\min_{\ell\in\mathcal P_{H,k}(S)}
\sum_{A\in\ell}\phi_{\mathrm{SSE}}(A)
&=
\min_{B\in\mathcal A_H(S,k)}
\min_{\ell'\in\mathcal P_{H,k-1}(S\setminus B)}
\left\{\phi_{\mathrm{SSE}}(B)
+\sum_{A\in\ell'}\phi_{\mathrm{SSE}}(A)\right\}\\
&=
\min_{B\in\mathcal A_H(S,k)}
\{\phi_{\mathrm{SSE}}(B)+F_H(S\setminus B,k-1)\}\\
&=F_H(S,k).
\end{align*}
Every remaining set \(S\setminus B\) in the inner minimum is a proper subset of \(S\), so the use
of the induction hypothesis is valid. Taking \(S=V\) proves the fixed-count
claim.

For the penalized recurrence, put
\(\mathcal P_H(S)=\bigsqcup_{k=0}^{|S|}\mathcal P_{H,k}(S)\) and
\[
J_\lambda(\ell)
:=
\sum_{A\in\ell}\{\phi_{\mathrm{SSE}}(A)+\lambda\}.
\]
For nonempty \(S\), taking the disjoint union of the preceding first-block
bijections over \(k\) gives
\[
\mathcal P_H(S)
\cong
\bigsqcup_{B\in\mathcal A_{H,\lambda}(S)}
\bigl(\{B\}\times\mathcal P_H(S\setminus B)\bigr).
\]
Since
\(J_\lambda(\{B\}\cup\ell')=
\phi_{\mathrm{SSE}}(B)+\lambda+J_\lambda(\ell')\), minimizing over this
finite disjoint union gives, with
\(G_\lambda(S):=\min_{\ell\in\mathcal P_H(S)}J_\lambda(\ell)\),
\[
G_\lambda(S)
=
\min_{B\in\mathcal A_{H,\lambda}(S)}
\{\phi_{\mathrm{SSE}}(B)+\lambda+G_\lambda(S\setminus B)\}.
\]
At \(S=\emptyset\), the only admitted object is the empty partition and
\(G_\lambda(\emptyset)=0=F_{H,\lambda}(\emptyset)\). Strong induction on
\(|S|\) now permits substitution of
\(G_\lambda(S\setminus B)=F_{H,\lambda}(S\setminus B)\) in the display for
every candidate \(B\); comparison with
\eqref{eq:penalized-subset-recurrence} gives
\(G_\lambda(S)=F_{H,\lambda}(S)\). This proves every stated boundary and
recurrence case.

\subsection{Proof of Proposition~\ref{prop:residual-component-factorization}}
\label{app:proof-component-factorization}

Let \(S_1,\ldots,S_c\) be the nonempty components of \(H[S]\), in any fixed
order. If a connected block \(A\subseteq S\) met two components, a path in
\(H[A]\subseteq H[S]\) would join them, a contradiction. Hence every block
of \(\ell\in\mathcal P_{H,k}(S)\) lies in exactly one \(S_j\). Define
\[
R(\ell)
=
(\ell|_{S_1},\ldots,\ell|_{S_c}),
\qquad
\ell|_{S_j}:=\{A\in\ell:A\subseteq S_j\}.
\]
If \(k_j=\#(\ell|_{S_j})\), then \(1\le k_j\le|S_j|\) and
\(\sum_jk_j=k\), so \(R\) has the codomain on the right of
\eqref{eq:residual-component-bijection}. In the other direction define
\[
U(\ell_1,\ldots,\ell_c)=\bigcup_{j=1}^c\ell_j.
\]
The vertex sets \(S_j\) are disjoint, so this union is an \(H\)-realizable
partition of \(S\) with \(\sum_jk_j=k\). Directly,
\(U\{R(\ell)\}=\ell\) and
\(R\{U(\ell_1,\ldots,\ell_c)\}=(\ell_1,\ldots,\ell_c)\), proving the
bijection.

For a fixed allocation \(\mathbf{k}\), block additivity gives
\[
\min_{(\ell_1,\ldots,\ell_c)\in
   \prod_j\mathcal P_{H,k_j}(S_j)}
\sum_{j=1}^c\sum_{A\in\ell_j}\phi_{\mathrm{SSE}}(A)
=
\sum_{j=1}^cF_H(S_j,k_j).
\]
The equality holds because the product is finite and the choices of the
\(\ell_j\)'s are independent: choosing a minimizer in every component attains
the sum of the componentwise lower bounds. Minimizing this display over
\(\mathbf{k}\in\mathcal K_H(S,k)\) proves
\eqref{eq:fixed-component-factorization}. The same argument without a count
allocation gives
\[
\min_{(\ell_1,\ldots,\ell_c)\in\prod_j\mathcal P_H(S_j)}
\sum_{j=1}^cJ_\lambda(\ell_j)
=
\sum_{j=1}^cF_{H,\lambda}(S_j),
\]
which is \eqref{eq:penalized-component-factorization}.

For feasibility, every component of \(H[S]\) needs at least one block and
every block is nonempty, so \(c\le k\le |S|\) is necessary. Conversely,
choose a spanning forest of \(H[S]\), with \(c\) trees and \(|S|-c\) edges.
The inequalities give \(0\le k-c\le |S|-c\), so delete any \(k-c\) forest
edges. Each deletion increases the number of components by one, leaving
exactly \(k\) nonempty vertex sets connected by forest edges. Their induced
subgraphs in \(H\) are therefore connected, proving sufficiency.

\subsection{Proof of Proposition~\ref{prop:scaffold-monotonicity}}
\label{app:proof-scaffold-monotonicity}

Let \(E_H\subseteq E_{H'}\). Every path in \(H[A]\) is also a path in
\(H'[A]\), so
\(\mathcal C(H)\subseteq\mathcal C(H')\). Intersecting these dictionaries
with the same first-block condition \(a(S)\in B\) and size bound gives
\(\mathcal A_H(S,k)\subseteq\mathcal A_{H'}(S,k)\) for every \((S,k)\).

To prove the reached-state inclusion, stratify each closure by the length of a
path from \((V,K)\). The length-zero state is common. If an \(H\)-reachable
state at length \(t+1\) is obtained from an \(H\)-reachable \((S,k)\) at
length \(t\) by removing \(B\in\mathcal A_H(S,k)\), the induction hypothesis
and the transition inclusion show that the same removal is available in
\(H'\). Hence \(\mathcal R_{H,K}\subseteq\mathcal R_{H',K}\). Likewise,
if every block of \(\ell\) is connected in \(H\), every block is connected in
\(H'\), proving \(\mathcal F(H)\subseteq\mathcal F(H')\).

Taking cardinalities yields the \(C_s\) inequalities. For
\(T_K^{\mathrm{reach}}\), every old reached state is also new and has a nested
transition family, while states belonging only to
\(\mathcal R_{H',K}\) contribute nonnegative terms. The terminal indicator
is also monotone because \(S\in\mathcal C(H)\) implies
\(S\in\mathcal C(H')\). Thus
\(T_K^{\mathrm{reach}}(H)\le T_K^{\mathrm{reach}}(H')\).

Finally, minimization over nested feasible families reverses the inequality:
\(\mathrm{OPT}(H,K)\ge\mathrm{OPT}(H',K)\), with the convention
\(\min\emptyset=+\infty\). If \(H\) covers an unrestricted optimum, then
\[
\mathrm{OPT}_{\mathrm{unc}}(K)
\le \mathrm{OPT}(H',K)
\le \mathrm{OPT}(H,K)
=\mathrm{OPT}_{\mathrm{unc}}(K),
\]
so equality holds throughout and coverage persists.

\section{Deterministic coverage proofs}
\label{app:deterministic-coverage-proofs}

\subsection{Nested thresholds and the lexicographic recurrence}
\label{app:proof-bottleneck-lexicographic}

For Proposition~\ref{prop:nested-edge-birth-threshold}, nestedness gives
\(e\in E(H_j)\) exactly when \(\iota_{\mathcal H}(e)\le j\). A finite
graph is connected exactly when it contains a spanning tree, so for
\(|A|\ge2\),
\[
H_j[A]\text{ connected}
\ \Longleftrightarrow\
\exists T\in\mathsf{ST}(A):\max_{e\in E(T)}\iota_{\mathcal H}(e)\le j
\ \Longleftrightarrow\ \kappa_{\mathcal H}(A)\le j.
\]
A singleton is connected at every index. Applying this criterion to every
block gives \(\ell\in\mathcal F(H_j)\) exactly when
\(\kappa_{\mathcal H}(\ell)\le j\). By
\eqref{eq:realizability-criterion}, coverage is equivalent to existence of
\(\ell\in\mathcal M_K\) satisfying this inequality. The nonempty finite
set \(\mathcal M_K\) attains its minimum threshold, proving
\eqref{eq:optimal-birth-threshold}.

For \(H_q=T_{\mathbf X}\cup Q_q(\mathbf X)\), under the tied-shell
convention, define
\(\rho_i(j)=1+\#\{s\notin\{i,j\}:\|\mathbf x_i-\mathbf x_s\|
<\|\mathbf x_i-\mathbf x_j\|\}\). Then
\begin{equation}
\iota_{q\mathrm{MST}}(\{\mathbf x_i,\mathbf x_j\})=
\begin{cases}
0,&\{\mathbf x_i,\mathbf x_j\}\in E(T_{\mathbf X}),\\
\min\{\rho_i(j),\rho_j(i)\},&\text{otherwise}.
\end{cases}
\label{eq:qnn-mst-edge-birth}
\end{equation}
The strict comparison includes a whole tied shell at its first admissible
rank. With exactly \(q\) neighbours retained under a deterministic tie
order, use that order's ordinal ranks instead. The standard
minimum-bottleneck property of an MST identifies
\(\kappa_{\mathcal H}(A)\) as the largest edge weight in an MST on \(A\)
weighted by these birth indices; Kruskal's algorithm obtains such a tree
\citep{kruskal1956shortest,cormen2022introduction}.

To minimize the threshold over tied SSE optima, compare partitions by
\begin{equation}
\left(\sum_{A\in\ell}\phi_{\mathrm{SSE}}(A),
\max_{A\in\ell}\kappa_{\mathcal H}(A)\right)
\label{eq:lexicographic-first-cover}
\end{equation}
in lexicographic order. Let \(L(S,k)=(L_1,L_2)\) be the minimum pair over
\(k\)-partitions of \(S\), with \(L(\emptyset,0)=(0,0)\) and
\((+\infty,+\infty)\) for infeasible states. The remaining-set recurrence is
\[
L(S,k)=\mathop{\operatorname{lexmin}}_{B\in\mathcal B(S,k)}
\left(\phi_{\mathrm{SSE}}(B)+L_1(S\setminus B,k-1),
\max\{\kappa_{\mathcal H}(B),L_2(S\setminus B,k-1)\}\right).
\]
For fixed \(B\), combining with a pair \((s,t)\) from the problem on
\(S\setminus B\) preserves a strict inequality in \(s\), and when \(s\)
ties is nondecreasing in \(t\).
It therefore suffices to retain the minimum pair for this remaining-set problem. Induction on
\(|S|\), using the first-block partition bijection in
Appendix~\ref{app:proof-exact-scaffold-dp}, proves the recurrence: each
partition contributes its unique first block and a partition of the smaller
remaining set, while each such pair gives a valid partition. At \((V,K)\), the
first coordinate is \(\mathrm{OPT}_{\mathrm{unc}}(K)\); minimizing the
second over its ties gives \(\kappa_{\mathcal H,K}^\star(\mathbf X)\).

\subsection{MST, radius, and nearest-neighbour connectivity}
\label{app:proof-mst-separated-coverage}

For Proposition~\ref{prop:mst-separated-coverage}, let \(T\) be any
Euclidean MST and suppose \(T[A_k]\) is disconnected. A within-block MST
contains an edge \(e=\{u,v\}\) joining two of these components, with
\(|e|\le r_{\mathrm{in}}(\ell)\). The \(T\)-path between its endpoints
must leave \(A_k\), so it contains an edge \(f\) between blocks with
\(|f|\ge r_{\mathrm{out}}(\ell)>|e|\). Adding \(e\) and removing \(f\)
gives a shorter spanning tree, a contradiction. Every block is therefore
connected in \(T\). If \(\ell\in\mathcal M_K\),
\eqref{eq:realizability-criterion} and graph inclusion prove coverage by
\(T\) and its supergraphs.

\label{app:proof-radius-qnn}
For Proposition~\ref{prop:radius-qnn-coverage}, take a longest edge
\(e=\{\mathbf u,\mathbf v\}\) in an MST of a nonsingleton block \(A\),
with length \(b_A\). Removing it yields sets \(A_1,A_2\), containing
\(\mathbf u,\mathbf v\), respectively. Every edge between these sets has
length at least \(b_A\), since otherwise exchanging it for \(e\) would
shorten the tree. Thus \(G_r[A]\) is disconnected for \(r<b_A\), whereas
for \(r\ge b_A\) it contains the entire within-block MST.

For the fill bound, let \(\mathbf z=(\mathbf u+\mathbf v)/2\). If
\(\mathbf x\in A_1\), then
\[
\|\mathbf x-\mathbf z\|
\ge\|\mathbf x-\mathbf v\|-\|\mathbf v-\mathbf z\|
\ge b_A/2.
\]
For \(\mathbf x\in A_2\), use \(\mathbf u\) instead. The endpoints
attain equality, and \(\mathbf z\in\operatorname{conv}(A)\), whence
\(b_A/2=\min_{\mathbf x\in A}\|\mathbf z-\mathbf x\|\le h(A)\).
For a singleton both quantities are zero. Maximizing over blocks proves
\eqref{eq:radius-bottleneck-fill}, including coincident observations.
Finally, a radius-\(r\) neighbour satisfies
\[
\rho_i(j)=1+\#\{s\notin\{i,j\}:\|\mathbf x_i-\mathbf x_s\|
<\|\mathbf x_i-\mathbf x_j\|\}\le N_i(r)\le M_{\mathbf X}(r).
\]
If \(M_{\mathbf X}(r)\le q\), every radius edge is therefore retained in
\(Q_q\). Combining the radius criterion with this inclusion proves the
remaining claims.

\section{Probabilistic nearest-neighbour coverage proofs}
\label{app:probabilistic-coverage-proofs}

\subsection{Population regularity and fitted-cell geometry}
\label{app:proof-risk-gap}

For \(\mathbf c,\mathbf c'\in\mathcal C_R\), the minimum inequality and
difference of squares give
\[
\left|\min_k\|\mathbf x-\mathbf c_k\|^2-
\min_k\|\mathbf x-\mathbf c'_k\|^2\right|
\le\max_k\left|\|\mathbf x-\mathbf c_k\|^2-
\|\mathbf x-\mathbf c'_k\|^2\right|
\le2(\|\mathbf x\|+R)d_\infty(\mathbf c,\mathbf c').
\]
The coefficient is integrable, so \(\mathcal Q_P\) is continuous.
Consequently \(\mathfrak C_{P,R}\) is compact and \(d_{P,R}\) is
continuous. The set \(\{\mathbf c\in\mathcal C_R:
d_{P,R}(\mathbf c)\ge\epsilon\}\) is compact and contains no optimum.
If nonempty, the continuous excess risk attains a positive minimum there;
otherwise the gap is \(+\infty\). This proves the risk-gap claim without
uniqueness.

\label{app:compact-centre-regularity}
Under compact ball support and positive lower density, projecting centres
onto the ball cannot increase risk; compactness therefore gives a minimum
\(Q_j^\star\) for every finite centre count \(j\). For an optimal
\(j\)-centre configuration, choose an interior point \(\mathbf z\) not
among its centres, put \(a=\min_k\|\mathbf z-\mathbf c_k\|>0\), and
choose \(0<t<\min\{a/4,R-\|\mathbf z\|\}\). Adding \(\mathbf z\)
as a centre reduces risk on \(B(\mathbf z,t)\) by at least
\[
P\{B(\mathbf z,t)\}\{(a-t)^2-t^2\}
\ge\underline f_Rv_dt^d(a^2-2at)>0.
\]
Hence \(Q_{j+1}^\star<Q_j^\star\). Coincident centres or a zero-mass
tie-broken cell would permit deleting a centre without changing risk, so
neither can occur in a \(K\)-centre optimum. For a positive-mass cell
\(D_k\), its conditional mean \(\mathbf m_k\) satisfies
\[
\int_{D_k}\|\mathbf x-\mathbf c_k\|^2dP
=\int_{D_k}\|\mathbf x-\mathbf m_k\|^2dP
+P(D_k)\|\mathbf c_k-\mathbf m_k\|^2.
\]
Recentring and nearest-centre reassignment cannot increase risk, forcing
\(\mathbf c_k=\mathbf m_k\) at an optimum. Absolute continuity gives
\(\|\mathbf x\|<R\) almost surely, and therefore
\(\|\mathbf c_k\|\le\mathbb E(\|\mathbf x\|\mid D_k)<R\).
All centre conditions follow. Compactness of the optimal set makes the
positive continuous separation and boundary margins attain positive minima
\(\Delta_*\) and \(\rho_*\).

\label{app:proof-localization-thickness}
To prove Lemma~\ref{lem:codebook-localization-thickness}, continuity on
\(\mathcal C_R\) first gives an empirical minimizer. The main-text risk
inequality shows that \(2Z_{n,R}<g_{P,R}(\epsilon_0)\) places every such
minimizer in \(\mathcal U_0\). For \(\mathbf c\in\mathcal U_0\), choose
a population optimum within \(\epsilon_0\). The triangle inequality gives
centre separation at least \(\Delta_*-2\epsilon_0\) and boundary distance
at least \(\rho_*-\epsilon_0\). Thus for
\(\mathbf x\in B(\mathbf c_k,a_0)\), one has \(\|\mathbf x\|\le R\)
and, for \(l\ne k\),
\[
\|\mathbf x-\mathbf c_l\|
\ge\|\mathbf c_l-\mathbf c_k\|-\|\mathbf x-\mathbf c_k\|
\ge\Delta_*-2\epsilon_0-a_0
\ge a_0\ge\|\mathbf x-\mathbf c_k\|.
\]
This proves
\begin{equation}
B(\mathbf c_k,a_0)\subseteq C_k(\mathbf c).
\label{eq:uniform-cell-inball}
\end{equation}
For \(\mathbf z\in C_k(\mathbf c)\) and \(0<s\le2R+a_0\), set
\(\lambda=s/(2R+a_0)\). Convexity gives
\[
B\{(1-\lambda)\mathbf z+\lambda\mathbf c_k,\lambda a_0\}
\subseteq C_k(\mathbf c)\cap B(\mathbf z,s),
\]
because this ball is the convex image of \(B(\mathbf c_k,a_0)\), and
each of its points is at distance at most
\(\lambda(\|\mathbf c_k-\mathbf z\|+a_0)\le\lambda(2R+a_0)=s\)
from \(\mathbf z\). Its volume is
\(v_d(\lambda a_0)^d=\theta_0v_ds^d\). Integration of
\(f\ge\underline f_R\) proves the lemma's probability bound.

\subsection{Uniform sampling within fitted cells}
\label{subsec:moving-cell-epsilon-net-proof}

For the fixed radius \(r\), the family
\begin{equation}
\mathcal A_r=\{B(\mathbf0,R)\cap D_k(\mathbf c)\cap B(\mathbf z,r/4):
\mathbf c\in\mathcal U_0,\ \mathbf z\in C_k(\mathbf c),\ k\le K\}
\label{eq:moving-voronoi-range-class}
\end{equation}
is specified before sampling. Fitted centres subsequently select members
of this family, so a uniform bound is required. For a family \(\mathcal A\),
let \(\Pi_{\mathcal A}(m)=\sup_{|S|=m}|\{A\cap S:A\in\mathcal A\}|\)
count its largest number of membership patterns on \(m\) points. Balls
and affine halfspaces in \(\mathbb R^d\) have VC dimension \(d+1\).
Sauer's inequality, \(\Pi_{\mathcal A}(m)\le(em/v)^v\) for
\(m\ge v\ge\operatorname{VCdim}(\mathcal A)\), gives
\begin{equation}
\Pi_{\mathcal A_r}(m)\le
\left(\frac{em}{d+1}\right)^{K(d+1)},\qquad m\ge d+1.
\label{eq:moving-range-growth-bound}
\end{equation}
Indeed, intersecting one variable ball and \(K-1\) halfspaces multiplies
their numbers of membership patterns; intersection with a fixed ball or
restriction of the parameters can only reduce that number. Open versus
closed halfspaces have the same finite-set patterns under perturbation.

Put \(p=K(d+1)\) and \(L=\log_2(2eK)\). Shattering \(m\ge d+1\)
points would imply \(2^{m/p}\le eK(m/p)\). The function \(2^y/y\)
increases for \(y>1/\log2\), and at \(y=2L\) its value is
\((2eK)^2/(2L)>eK\). Thus shattering is impossible for \(m\ge2pL\);
the smaller case \(m<d+1\) satisfies the same bound. Therefore
\begin{equation}
\operatorname{VCdim}(\mathcal A_r)\le
V_{d,K}=\lceil2K(d+1)\log_2(2eK)\rceil.
\label{eq:moving-range-vc-dimension}
\end{equation}

The standard double-sampling \(\eta\)-net inequality
\citep{haussler1987epsilon} states that, for \(n\ge v\), VC dimension at
most \(v\ge1\), and \(0<\eta\le1\) with \(n\eta\ge8\log2\),
\begin{equation}
\mathbb P^*\{\exists A\in\mathcal A:P(A)\ge\eta,\ S_n\cap A=\varnothing\}
\le\min\{1,2\Pi_{\mathcal A}(2n)e^{-n\eta/4}\}
\le\Psi_n(\eta,v).
\label{eq:epsilon-net-derived-bound}
\end{equation}
Here \(S_n\) is an independent sample from \(P\), and \(\mathbb P^*\)
is outer probability. For \(n\eta<8\log2\), the final expression equals
one, since \(n\ge v\ge1\); it remains a valid bound. In our family,
membership is Borel measurable jointly in the observation and parameters,
as is its integrated probability. The missing-region event is a projection
of a Borel set, hence analytic and measurable under the completed sampling
law. Ordinary probability may therefore replace outer probability.

Cell boundaries are finite unions of null hyperplanes, so Lemma~
\ref{lem:codebook-localization-thickness} gives, for every member of
\(\mathcal A_r\),
\[
P\{B(\mathbf0,R)\cap D_k(\mathbf c)\cap B(\mathbf z,r/4)\}
=P\{C_k(\mathbf c)\cap B(\mathbf z,r/4)\}
\ge\underline f_R\theta_0v_d(r/4)^d=\eta_r.
\]
Thus \(0<\eta_r\le1\), and \eqref{eq:epsilon-net-derived-bound} yields
\begin{equation}
\mathbb P\{\exists A\in\mathcal A_r:S_n\cap A=\varnothing\}
\le\Psi_n(\eta_r,V_{d,K}).
\label{eq:uniform-moving-cell-fill}
\end{equation}

\subsection{Uniform neighbour counts}
\label{subsec:radius-neighbour-count-proof}
For the neighbour counts, write
\(p_+(r)=\sup_{\mathbf x}P\{B(\mathbf x,r)\}\),
\(\mu_n(r)=(n-1)p_+(r)\), and define
\begin{equation}
b(\lambda,u)=\lambda+\sqrt{2\lambda u}+2u/3.
\label{eq:direct-binomial-threshold}
\end{equation}
\begin{lemma}[Uniform radius-neighbour count]
\label{lem:uniform-radius-neighbour-count}
For every \(u\ge0\),
\begin{equation}
\mathbb P\{\max_iN_i(r)>b(\mu_n(r),u)\}\le ne^{-u}.
\label{eq:uniform-radius-neighbour-count}
\end{equation}
Under \eqref{eq:density-sandwich},
\(\mu_n(r)\le\overline\mu_n(r)=(n-1)\overline f v_dr^d\).
\end{lemma}

Conditionally on \(\mathbf x_i=\mathbf x\), the other observations are
independent, giving \(N_i(r)\sim\operatorname{Bin}(n-1,p_{\mathbf x}(r))\)
with \(p_{\mathbf x}(r)=P\{B(\mathbf x,r)\}\) and mean
\(\lambda_{\mathbf x}\le\mu_n(r)\). Binomial Bernstein concentration
\citep[Chapter~2]{boucheron2013concentration} bounds its upper tail by
\(\exp\{-t^2/[2(\lambda_{\mathbf x}+t/3)]\}\) at mean plus \(t>0\).
For \(u>0\), substituting
\(t=\sqrt{2\lambda_{\mathbf x}u}+2u/3\) gives
\[
t^2-2u(\lambda_{\mathbf x}+t/3)
=(2u/3)\sqrt{2\lambda_{\mathbf x}u}\ge0,
\]
so the conditional probability of exceeding \(b(\lambda_{\mathbf x},u)\)
is at most \(e^{-u}\). A zero mean makes the count zero almost surely.
Monotonicity of \(b\) in its first argument, integration over \(\mathbf x_i\),
and a union bound over \(i\) prove the lemma; the counts need not be
independent. For \(u=0\) the bound is trivial. Finally,
\(p_{\mathbf x}(r)=\int_{B(\mathbf x,r)}f(\mathbf y)d\mathbf y
\le\overline f v_dr^d\), proving the stated mean bound.

\subsection{Exact empirical coverage}
\label{app:proof-empirical-qnn-coverage}

In Theorem~\ref{thm:qnn-mixture-coverage}, let \(\mathcal E_R\) mean that
all observations lie in the ball, \(\mathcal E_{\mathrm L}\) that
\(2Z_{n,R}<g_{P,R}(\epsilon_0)\), \(\mathcal E_{\mathrm F}\) that every
region in \(\mathcal A_r\) contains an observation, and
\(\mathcal E_{\mathrm C}\) that
\(M_{\mathbf X}(r)\le b\{\overline\mu_n(r),\log(n/\delta_{\mathrm g})\}\).
The preceding results give
\[
\mathbb P(\mathcal E_R^c)\le n\tau_R,\quad
\mathbb P(\mathcal E_{\mathrm L}^c)=\zeta_{n,R}(\epsilon_0),\quad
\mathbb P(\mathcal E_{\mathrm F}^c)\le\Psi_n(\eta_r,V_{d,K}),\quad
\mathbb P(\mathcal E_{\mathrm C}^c)\le\delta_{\mathrm g}.
\]
We show that their intersection implies coverage.

An empirical global minimum exists because
\[
n\min_{\mathbf c}\mathcal Q_n(\mathbf c)
=\min_{\mathbf z\in\{1,\ldots,K\}^n}
\min_{\mathbf c_1,\ldots,\mathbf c_K}
\sum_i\|\mathbf x_i-\mathbf c_{z_i}\|^2.
\]
There are finitely many assignments; each nonempty group is minimized at
its sample mean, and unused centres may be placed at the origin. On
\(\mathcal E_R\) these choices are in \(\mathcal C_R\), so choose a global
empirical minimizer \(\widehat{\mathbf c}\) there. On
\(\mathcal E_{\mathrm L}\), it belongs to \(\mathcal U_0\).

Assign \(A_k=\{i:\mathbf x_i\in D_k(\widehat{\mathbf c})\}\). For
every \(\mathbf z\in C_k(\widehat{\mathbf c})\), the corresponding member
of \(\mathcal A_r\) contains an observation from \(A_k\) on
\(\mathcal E_{\mathrm F}\). Thus \(A_k\ne\emptyset\) and
\[
\sup_{\mathbf z\in C_k(\widehat{\mathbf c})}
\operatorname{dist}(\mathbf z,\mathbf X_{A_k})\le r/4.
\]
Since the closed cell is convex and contains \(\mathbf X_{A_k}\), it also
contains its convex hull. Hence \(h(\mathbf X_{A_k})\le r/4\), and the
exactly \(K\)-block partition \(\widehat\ell\) satisfies
\(r_{\mathrm{in}}(\widehat\ell)\le2h(\widehat\ell)\le r/2\le r\).
The uniform event \(\mathcal E_{\mathrm F}\), rather than an argument for
a fixed cell, justifies selecting \(\widehat{\mathbf c}\) from these same data.

Recentring the assigned nonempty blocks and comparing with arbitrary
\(K\)-partitions gives
\[
\mathrm{SSE}(\widehat\ell)
\le n\mathcal Q_n(\widehat{\mathbf c})
=n\min_{\mathbf c}\mathcal Q_n(\mathbf c)
\le\mathrm{OPT}_{\mathrm{unc}}(K)
\le\mathrm{SSE}(\widehat\ell).
\]
The middle inequality holds because any partition's means supply a centre
configuration with no larger nearest-centre cost. Equality throughout proves
\(\widehat\ell\in\mathcal M_K\).

If the uncapped threshold in \eqref{eq:qnn-exact-rank-threshold} is below
\(n-1\), \(\mathcal E_{\mathrm C}\) implies
\(M_{\mathbf X}(r)\le q_n(r,\delta_{\mathrm g})\); otherwise
\(q_n=n-1\) and the neighbour graph is complete. In either case,
Proposition~\ref{prop:radius-qnn-coverage} makes \(\widehat\ell\) feasible
on every \(H_q\supseteq Q_q\) with \(q\ge q_n\). The union bound on the
four complementary events proves \eqref{eq:qnn-exact-coverage-bound}.

\subsection{Compact support and the Gaussian-mixture specialization}
\label{app:proof-compact-qnn-corollary}
\label{app:proof-gaussian-coverage-bounds}

Under the compact-support hypotheses of Corollary~\ref{cor:qnn-compact-exact},
the preceding regularity argument verifies Assumption~
\ref{ass:qnn-core-regularity}. Put
\(\kappa=\underline f_R\theta_0v_d/4^d\), choose
\(A>4(V_{d,K}+2)/\kappa\), and set
\(r_n^d=A\log n/n\), \(\delta_{\mathrm g}=n^{-2}\), and
\(q_n=q_n(r_n,n^{-2})\) for all sufficiently large \(n\). The finitely
many initial ranks may be set to \(n-1\).

For completeness, the uniform risk discrepancy has an explicit bound.
Writing \(f_{\mathbf c}(\mathbf x)=\min_k\|\mathbf x-\mathbf c_k\|^2\),
on ball support \(0\le f_{\mathbf c}\le4R^2\) and
\(|f_{\mathbf c}-f_{\mathbf c'}|\le4R d_\infty(\mathbf c,\mathbf c')\).
For \(0<t\le4R^2\), a maximal \(t/(16R)\)-separated subset of the ball
covers it at that distance. Disjoint balls of half that radius, contained
in the correspondingly enlarged ball, bound its size by
\((1+32R^2/t)^d\). Taking \(K\)-fold products covers
\(\mathcal C_R\) in \(d_\infty\) by at most
\(N_t=(1+32R^2/t)^{dK}\) configurations. Approximating any configuration
changes its empirical and population averages by at most \(t/4\) each.
Therefore \(Z_{n,R}\ge t\) implies a discrepancy of at least \(t/2\)
at a covering configuration. Hoeffding's inequality for bounded variables
\citep[Chapter~2]{boucheron2013concentration} and a union bound give
\begin{equation}
\mathbb P\{Z_{n,R}\ge t\}
\le2\left(1+\frac{32R^2}{t}\right)^{dK}
\exp\left\{-\frac{nt^2}{32R^4}\right\}.
\label{eq:compact-uniform-risk-concentration}
\end{equation}
The exclusion set for \(g_0=g_{P,R}(\epsilon_0)\) is nonempty: coinciding
all centres at a boundary point places the configuration at distance at
least \(\rho_*>2\epsilon_0\) from the optimal set. Since risks lie in
\([0,4R^2]\), \(0<g_0\le4R^2\). Substituting \(t=g_0/2\) yields
\begin{equation}
\zeta_{n,R}(\epsilon_0)
\le2\left(1+\frac{64R^2}{g_0}\right)^{dK}
\exp\left\{-\frac{ng_0^2}{128R^4}\right\}\longrightarrow0.
\label{eq:compact-codebook-localization-bound}
\end{equation}

For large \(n\), the radius is in the allowed range, and
\[
\eta_{r_n}=\kappa A\frac{\log n}{n},\qquad
\Psi_n(\eta_{r_n},V_{d,K})
\le2\left(\frac{2en}{V_{d,K}}\right)^{V_{d,K}}n^{-\kappa A/4}
\longrightarrow0.
\]
Our choice makes the last expression \(O(n^{-2-\varepsilon})\) for some
\(\varepsilon>0\). Also \(\overline\mu_n(r_n)\le A\overline f v_d\log n\)
and \(\log(n/\delta_{\mathrm g})=3\log n\), so all terms in the rank
threshold are logarithmic:
\begin{equation}
q_n(r_n,n^{-2})=O(\log n).
\label{eq:compact-q-rate}
\end{equation}
Together with \(\tau_R=0\) and \(\delta_{\mathrm g}\to0\), these bounds
prove \eqref{eq:compact-exactness-limit}. More explicitly, choosing
\(A=C_{d,K}/\kappa\), with \(C_{d,K}>4(V_{d,K}+2)\), gives
\(\overline\mu_n(r_n)\le L_0\log n\), where
\(L_0=C_{d,K}4^d\overline f/(\underline f_R\theta_0)\).
The remaining terms are at most \(\sqrt{6L_0}\log n\) and \(2\log n\).
Using \(\sqrt{x}\le(1+x)/2\) proves the stated
\(O_{d,K}([1+\overline f/(\underline f_R\theta_0)]\log n)\) dependence.

Absolute continuity makes each equality between two distances from one
sample point a null event; a finite union excludes all such ties almost
surely. The union neighbour graph then has at most \(nq_n\) edges, and
an MST adds at most \(n-1\). Hence
\begin{equation}
|E(H_{q_n})|\le nq_n+n-1,\qquad
\frac{|E(H_{q_n})|}{\binom n2}
\le\frac{2q_n}{n-1}+\frac2n=O(\log n/n)
\quad\text{almost surely}.
\label{eq:qnn-edge-compression}
\end{equation}

For the Gaussian specialization, let
\(\widetilde f(\mathbf x)=\sum_{a=1}^J\pi_a
\phi(\mathbf x;\boldsymbol\mu_a,\boldsymbol\Sigma_a)\), with positive
weights summing to one and positive-definite covariance matrices. Its
density conditioned on the ball is
\(f_R(\mathbf x)=\mathbf1\{\|\mathbf x\|\le R\}\widetilde f(\mathbf x)/p_R\),
where \(p_R=P_{\widetilde f}\{\|\mathbf x\|\le R\}>0\). Valid constants are
\begin{align}
\overline f_R&=\frac{(2\pi)^{-d/2}}{p_R}
\sum_{a=1}^J\pi_a|\boldsymbol\Sigma_a|^{-1/2},
\label{eq:truncated-gmm-upper-density}\\
\underline f_R&=\frac{(2\pi)^{-d/2}}{p_R}
\sum_{a=1}^J\pi_a|\boldsymbol\Sigma_a|^{-1/2}
\exp\left\{-\frac{(R+\|\boldsymbol\mu_a\|)^2}
{2\lambda_{\min}(\boldsymbol\Sigma_a)}\right\}.
\label{eq:truncated-gmm-lower-density}
\end{align}
The Gaussian exponential is at most one, proving the upper bound globally.
For \(\|\mathbf x\|\le R\),
\[
(\mathbf x-\boldsymbol\mu_a)^\mathsf T\boldsymbol\Sigma_a^{-1}
(\mathbf x-\boldsymbol\mu_a)
\le\frac{\|\mathbf x-\boldsymbol\mu_a\|^2}{\lambda_{\min}(\boldsymbol\Sigma_a)}
\le\frac{(R+\|\boldsymbol\mu_a\|)^2}{\lambda_{\min}(\boldsymbol\Sigma_a)},
\]
which proves the positive lower bound term by term. Taking
\(\overline f=\overline f_R\), compact regularity supplies all remaining
hypotheses, regardless of whether \(K=J\).

For the untruncated mixture, setting \(p_R=1\) gives the corresponding
density bounds. Conditional on component \(a\), write
\(\mathbf x=\boldsymbol\mu_a+\boldsymbol\Sigma_a^{1/2}\mathbf z\), where
\(\mathbf z\sim N(\mathbf0,\mathbf I_d)\). If \(R>\|\boldsymbol\mu_a\|\),
\(\|\mathbf x\|>R\) implies
\(\|\mathbf z\|^2>(R-\|\boldsymbol\mu_a\|)^2/
\lambda_{\max}(\boldsymbol\Sigma_a)\), using the triangle inequality and
the largest-eigenvalue bound. For smaller \(R\) a probability bound of
one suffices. With \((t)_+=\max(t,0)\), summing gives
\[
\tau_R\le\sum_{a=1}^J\pi_a\Pr\left\{\chi_d^2>
\frac{(R-\|\boldsymbol\mu_a\|)_+^2}{\lambda_{\max}(\boldsymbol\Sigma_a)}\right\}.
\]
Here \(\chi_d^2\) has \(d\) degrees of freedom. For this unbounded model,
the population-optimal centre conditions in the chosen ball remain explicit
assumptions; the compact-support argument does not verify them.

\section{Posterior invariance}

\subsection{Proof of Proposition~\ref{prop:scaffold-mh-invariance}}
\label{app:posterior-kernel-proof}

For distinct endpoint partitions \(\ell,\ell'\), let
\(\Omega_u(\ell,\ell')\) be the finite set of oriented auxiliary traces that
map \(\ell\) to \(\ell'\) at fixed \(u\). The deterministic convention
\(i<j\) labels the subclusters, and the deterministic scan fixes the order of
the allocation bits. A split trace reverses to the merge with the same
ordered seed pair. A merge trace reverses to the forced split assigning each
visited observation to its original source block. Applying reversal twice
restores the seed pair and every allocation decision, so
\[
\omega\mapsto\bar\omega:
\Omega_u(\ell,\ell')\longrightarrow\Omega_u(\ell',\ell)
\]
is an involutive bijection. Strict positivity of the block marginals,
\(\epsilon\), and the restricted-allocation weights supplies reverse support
whenever both endpoint moves are eligible. Eligibility is symmetric at the
two endpoints: in the finite model a forward split starts with \(K<L\) and
reverses by an always-available merge, whereas a forward merge starts with
\(K\ge2\) and reverses by a split from \(K-1<L\); for the DP and NGGP the
same count argument uses the upper bound \(K\le n\). We work on the
target-positive state space, so the current target value and every sampled
forward proposal factor are positive. A proposed path without positive
reverse support is assigned acceptance probability zero.

For \(\ell'\ne\ell\), define the endpoint kernel
\[
P_H(\ell,\ell'\mid u)
=
\sum_{\omega\in\Omega_u(\ell,\ell')}
\widetilde q_H(\omega\mid\ell,u)
\alpha_H(\ell,\ell';\omega\mid u),
\]
and set the diagonal probability to one minus the off-diagonal row sum. For
one paired path, writing \(\pi_u(\ell)=\pi_M(\ell\mid u,\mathbf{X})\),
\begin{align}
&\pi_u(\ell)\widetilde q_H(\omega\mid\ell,u)
\alpha_H(\ell,\ell';\omega\mid u)
\notag\\
&\quad=
\min\{\pi_u(\ell)\widetilde q_H(\omega\mid\ell,u),
    \pi_u(\ell')\widetilde q_H(\bar\omega\mid\ell',u)\}
\notag\\
&\quad=
\pi_u(\ell')\widetilde q_H(\bar\omega\mid\ell',u)
\alpha_H(\ell',\ell;\bar\omega\mid u).
\label{eq:scaffold-mh-detailed-balance}
\end{align}
The conditional normalizing constant cancels, so the joint target values in
\eqref{eq:scaffold-mh-acceptance} give the same ratio. Summing
\eqref{eq:scaffold-mh-detailed-balance} over
\(\Omega_u(\ell,\ell')\) and changing variables through the reversal
bijection proves off-diagonal detailed balance. The diagonal identity is
automatic, proving reversibility and invariance.

If each allocation and auxiliary-variable kernel preserves the joint target,
so does their composition with the split-merge kernel:
\(\pi P_1P_2P_3=((\pi P_1)P_2)P_3=\pi\).
Appendix~\ref{supp:auxiliary-update} verifies the NGGP update.
Irreducibility and aperiodicity are additional requirements on the full
composition, not consequences of the guided subkernel alone.

\section{Additional experimental results}
\label{sec:supp-scaled-controls}

\subsection{Objective-matched selection of the number of clusters}
\label{sec:supp-unknown-k}

Unknown-\(K\) comparisons use two block-additive criteria. For penalized SSE,
\[
L_\lambda(\ell)
=
\operatorname{SSE}(\ell)+\lambda|\ell|,
\qquad
\lambda=\sum_{j=1}^d\widehat{\operatorname{Var}}(\mathbf{X}_{\cdot j}),
\]
the comparator is DP-means and the complete recurrence exactly minimizes
\(L_\lambda\). For model-based MAP clustering, the score is the collapsed
normal-inverse-Wishart (NIW) likelihood with a Chinese restaurant process
(CRP) partition prior. We set
\(\alpha=1\), \(\boldsymbol{\mu}_0=\bar{\mathbf{x}}\), \(\kappa_0=0.1\),
\(\nu_0=d+2\), and choose the inverse-Wishart scale so that its prior mean is
the regularized empirical covariance. Three-restart MAP-DP coordinate ascent
is compared with exact recurrences for this same collapsed posterior score.

Table~\ref{tab:unknown-k-results} reports gaps only to the oracle for the
corresponding objective. DP-means cannot therefore be assessed by a
collapsed-posterior gap, nor MAP-DP by a penalized-SSE gap. Delaunay matched
the complete recurrence in all 120 data sets for both objectives and is
omitted because its rows would duplicate the complete results.

\begin{table}[!tbp]
\centering
\caption{Unknown-\(K\) objective-matched comparisons over 120 data sets. For
penalized SSE, the reported gap is the mean relative objective gap. For
collapsed MAP, it is the mean loss in log posterior, in nats. VI is measured
against the generating labels and is not part of either optimized objective.}
\label{tab:unknown-k-results}
\small
\setlength{\tabcolsep}{4pt}
\begin{tabular}{@{}L{1.3in}L{1.55in}rrrr@{}}
\toprule
Objective & Method & Oracle hits & Matched gap & Mean \(\widehat K\)
& Mean VI\\
\midrule
Penalized SSE
 & DP-means  & \(37/120\) & \(9.101\%\) & 3.600 & 0.883\\
 & \(2\)-NN-MST DP & \(117/120\) & \(0.025\%\) & 3.817 & 0.750\\
 & Gabriel DP  & \(120/120\) & \(0\%\) & 3.833 & 0.751\\
 & Complete DP  & \(120/120\) & \(0\%\) & 3.833 & 0.751\\
\addlinespace
Collapsed NIW-CRP MAP
 & MAP-DP coordinate ascent & \(92/120\) & \(0.186\) & 2.683 & 1.046\\
 & \(2\)-NN-MST DP  & \(107/120\) & \(0.047\) & 2.608 & 1.011\\
 & Gabriel DP  & \(119/120\) & \(0.001\) & 2.592 & 1.053\\
 & Complete DP  & \(120/120\) & \(0\) & 2.583 & 1.057\\
\bottomrule
\end{tabular}
\end{table}

The penalized results show that the remaining-set recurrence is not confined to
fixed \(K\): it gives a direct exact comparator for the DP-means objective.
The collapsed-MAP results make the separate statistical point. Moving toward
the exact posterior mode improves the optimized score, but does not
necessarily improve VI to a generating allocation. This is expected when the
fitted likelihood, partition prior, and data-generating labels define
different targets.

\subsection{Component-factorization ablation}
\label{sec:supp-frontier-pilot}

The implementation ablation uses three Gaussian components with covariance
\(\mathbf I_2\), means on an equilateral triangle of side \(\Delta\), and
weights \(((1-w)/2,(1-w)/2,w)\). One common stream of allocation uniforms
and Gaussian residuals generates \(n=12\) samples in four cells,
\((\Delta,w)\in\{1,3\}\times\{0.1,1/3\}\), without enforcing component
counts; the fitted count is \(K=3\).
The basic subset-DP and dictionary-based component-factorized evaluators solve
\(q\)-NN-MST for \(q\in\{0,1,2,3,5\}\), RNG, Gabriel, and Delaunay on each
sample. All 32 paired evaluations agree in objective value and feasible
partition count; independently recomputed SSE also agrees with the
backtracked solutions. This is a paired implementation check, not a
replicated estimate of model-dependent coverage.

\begin{table}[!tbp]
\centering
\caption{Component-factorized evaluation relative to basic subset-DP evaluation on
the same \(n=12\) sample and scaffold. The first two numerical columns are mean
paired ratios, expressed as percentages; the last two are mean additional
bookkeeping counts. Each row averages four geometry cells.
The \(q=1\) row duplicates MST and is omitted. These are operation and
storage counts, not total-memory or runtime speedups.}
\label{tab:frontier-factorization}
\small
\begin{tabular}{lrrrr}
\toprule
Scaffold & \shortstack{Stored state\\pairs (\%)}
& \shortstack{Inspected block\\candidates (\%)}
& \shortstack{Component-cache\\masks}
& \shortstack{Convolution\\relaxations}\\
\midrule
MST       & 55.800 & 23.800 & 111.000 & 70.000\\
RNG       & 75.400 & 26.900 & 237.000 & 198.500\\
Gabriel     & 83.900 & 28.800 & 536.000 & 432.500\\
\(2\)-NN-MST  & 70.100 & 20.300 & 476.250 & 361.500\\
\(3\)-NN-MST  & 95.800 & 44.500 & 1179.000 & 1072.000\\
Delaunay    & 102.200 & 67.800 & 1658.750 & 1342.500\\
\(5\)-NN-MST  & 100.900 & 89.000 & 1957.750 & 909.500\\
\bottomrule
\end{tabular}
\end{table}

Factorization reduces block inspections in these cells, but eight of the
32 evaluations retain more state-count pairs, including all four
\(5\)-NN-MST cases. Component-cache entries and convolution relaxations are
additional resources, so the ratios are not total-memory or runtime
speedups. This ablation retains union-state solutions; the selective cache
used for full-data UCI fits does not
(Appendix~\ref{subsec:lazy-scaffold-evaluation}).

\subsection{Structural controls and candidate connectivity}
\label{sec:supp-structural-controls}

\begin{figure}[!tbp]
\centering
\includegraphics[width=0.76\linewidth]{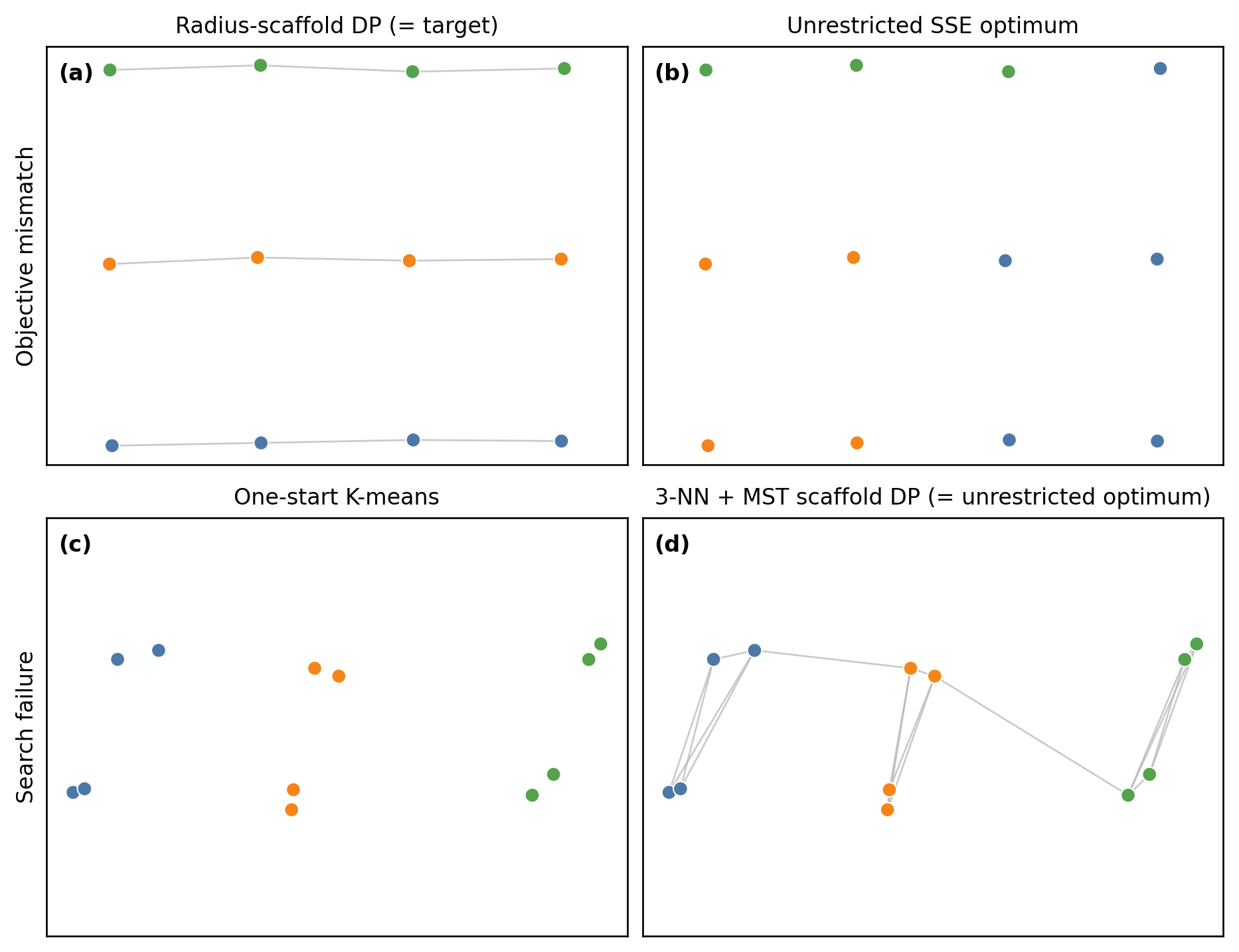}
\caption{The objective-mismatch chain control (top) and the local-search
offset-column control (bottom). Grey segments are scaffold edges.}
\label{fig:synthetic-structural-controls-supp}
\end{figure}

Two \(n=12\), \(K=3\) controls separate objective mismatch from search
failure (Figure~\ref{fig:synthetic-structural-controls-supp}). Across 20
samples, a radius-1.35 graph recovers the stipulated chain labels exactly
but excludes every unrestricted SSE optimum. The \(2\)-NN-MST and Delaunay
recurrences recover those optima in all samples, with mean target ARI
\(0.212\). In the offset-column control, \(3\)-NN-MST and Delaunay recover
both the optimum and target labels in every sample, whereas one-start Lloyd
has mean ARI \(0.763\). The chain is thus a negative control for the
centroid criterion, not an optimization failure.

\begin{figure}[!tbp]
\centering
\includegraphics[width=0.88\linewidth]{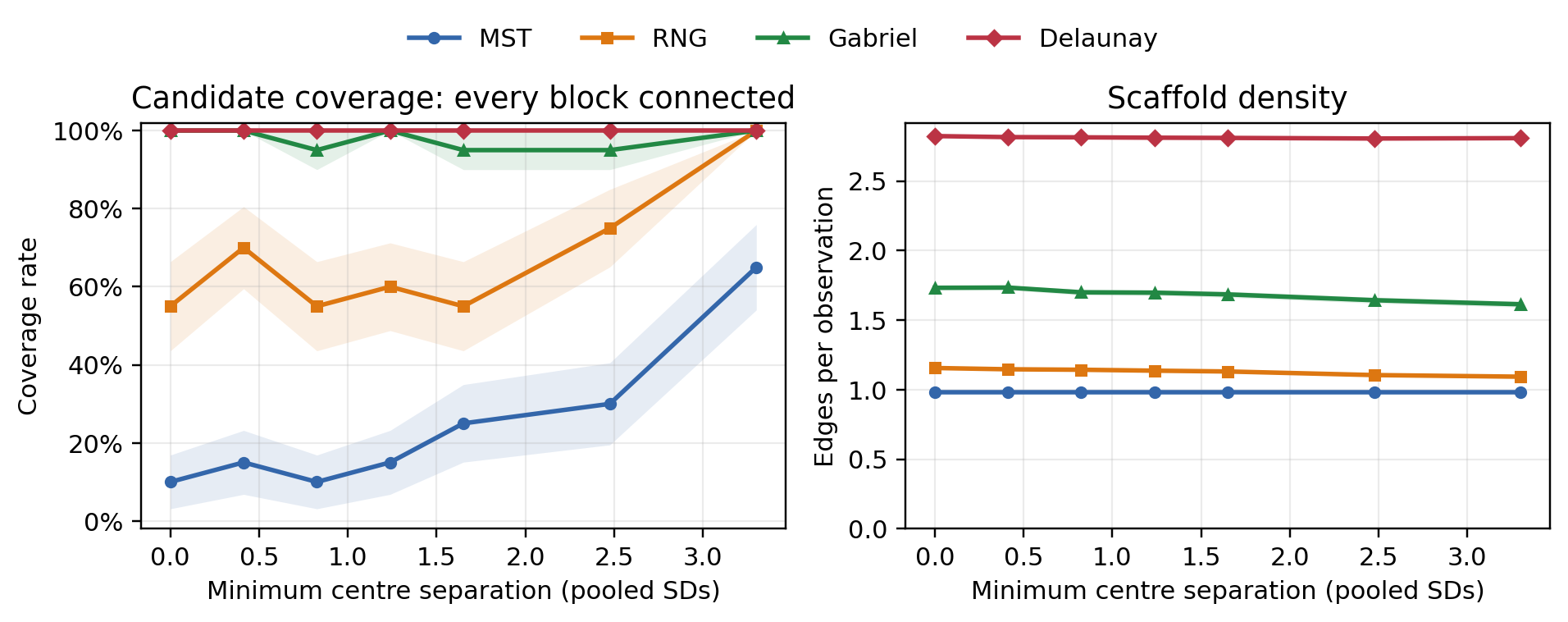}
\caption{Candidate connectivity at \(n=60\) over 20 data sets per separation
level (left), with edges per observation shown only as graph density (right).}
\label{fig:gmm-separability-supp}
\end{figure}

At \(n=60\), the separability sweep checks only a fitted candidate
\(\ell_c\), through
\begin{equation}
\operatorname{Cov}^{\mathrm{cand}}_H(\ell_c)
:=\mathbf{1}\{\ell_c\in\mathcal F(H)\}.
\label{eq:candidate-coverage}
\end{equation}
This does not certify coverage of an unknown unrestricted optimum.
Figure~\ref{fig:gmm-separability-supp} varies between-component separation
while holding centered within-component residuals fixed. Candidate coverage
usually increases, but not monotonically: moving the centres also changes the
fitted Voronoi allocation and the particular boundary observations that must
be joined. This is consistent with the separate roles of within-cell
sampling, local cell volume, neighbour counts, and centre stability in
Section~\ref{subsec:gmm-qnn-coverage}.

\subsection{Two-moons construction and noise sensitivity}
\label{sec:supp-moons}

Each sample contains 30 equally spaced locations on each of the semicircles
\((\cos t,\sin t)\) and \((1-\cos t,0.5-\sin t)\), with
\(t\in[0,\pi]\). Independent bivariate Gaussian perturbations have covariance
\(\sigma^2\mathbf I_2\), for
\(\sigma\in\{0.020,0.080,0.160\}\); observations are shuffled and are
not standardized. Five seeds are paired across noise levels and methods.
The experiment varies noise around fixed arc locations, rather than drawing
independent locations along the arcs. The displayed realization is fixed
as the first seed at the lowest noise level, without selection by fit quality.

The pruned forest is constructed by removing the longest edge from the same
MST used for the intact-tree fit. The connected-partition recurrence returns
exactly one feasible two-block partition for every forest, and its result
agrees with an independent connected-component calculation in all 15 cases.
On the intact MST, the recurrence counts exactly 59 feasible partitions.
The \(2\)-NN-MST fits retain the same MST and add union-neighbour edges;
their SSE cannot exceed the intact-tree minimum. These comparisons share
the same observations and SSE block costs.

\begin{table}[!tp]
\centering
\small
\caption{Two-moons results at \(n=60\), \(K=2\). Entries are means over
five paired realizations at each noise level; every displayed fit completed.
The pruned-MST fit is identical to single linkage and is shown only once.
SSE is the optimized loss for the DP fits; ARI and VI assess generating-label
agreement. Variability is illustrated in Figure~\ref{fig:moons-noise-supp}.}
\label{tab:moons-scaffold-results}
\begin{tabular}{llrrr}
\toprule
Noise SD & Method & SSE & ARI & VI (bits)\\
\midrule
0.020 & Lloyd & 24.095 & 0.224 & 1.642 \\
 & $k$-means++, 1 start & 24.093 & 0.224 & 1.642 \\
 & $k$-means++, 20 starts & 24.080 & 0.211 & 1.663 \\
 & DP, intact MST & 28.762 & 0.445 & 0.989 \\
 & DP, pruned MST & 37.293 & 1.000 & 0.000 \\
 & DP, $2$-NN-MST & 28.762 & 0.445 & 0.989 \\
\midrule
0.080 & Lloyd & 24.513 & 0.231 & 1.625 \\
 & $k$-means++, 1 start & 24.491 & 0.231 & 1.629 \\
 & $k$-means++, 20 starts & 24.452 & 0.238 & 1.618 \\
 & DP, intact MST & 29.327 & 0.455 & 0.978 \\
 & DP, pruned MST & 37.884 & 0.821 & 0.254 \\
 & DP, $2$-NN-MST & 29.308 & 0.446 & 0.988 \\
\midrule
0.160 & Lloyd & 25.906 & 0.274 & 1.545 \\
 & $k$-means++, 1 start & 25.878 & 0.274 & 1.546 \\
 & $k$-means++, 20 starts & 25.855 & 0.254 & 1.585 \\
 & DP, intact MST & 26.137 & 0.240 & 1.591 \\
 & DP, pruned MST & 50.935 & 0.046 & 1.245 \\
 & DP, $2$-NN-MST & 26.137 & 0.240 & 1.591 \\
\bottomrule
\end{tabular}

\end{table}

\begin{figure}[!tp]
\centering
\includegraphics[width=\linewidth]{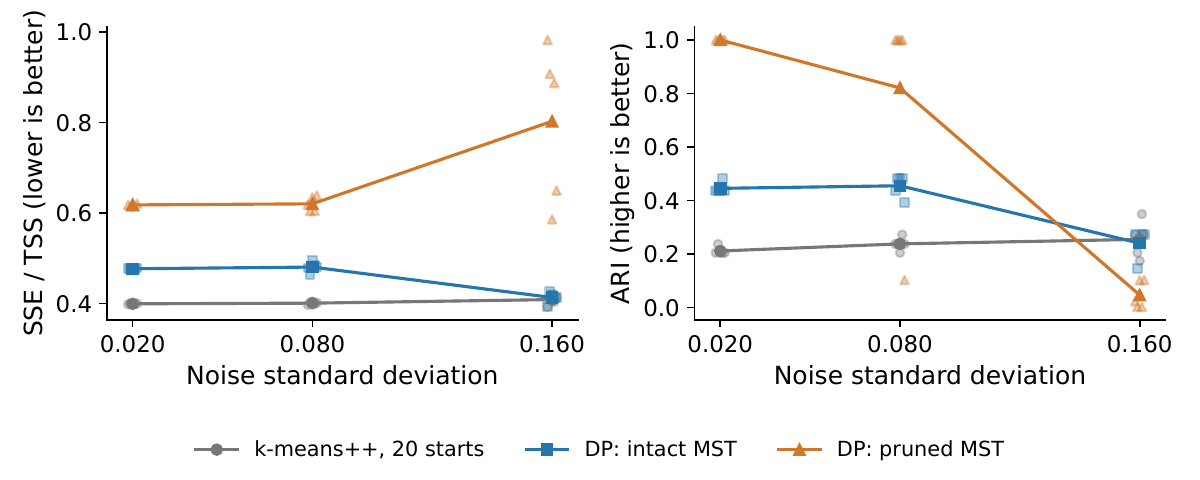}
\caption{Noise sensitivity of the two-moons comparison. Small translucent
points are individual realizations, slightly displaced horizontally for
visibility; larger joined points are means over five realizations.
Lower SSE/TSS and higher ARI measure different properties. Pruning the MST
improves label recovery at low noise but becomes unreliable as noise increases.}
\label{fig:moons-noise-supp}
\end{figure}

\FloatBarrier
All MST and \(2\)-NN-MST searches finish within the five-second solver
allowance, with limits of \(200{,}000\) retained plus active states and
\(2\times10^6\) generator nodes. In the wider scaffold sweep, all 45
searches on union \(5\)-NN, union \(5\)-NN plus MST, and union \(10\)-NN
plus MST reach the time limit; no fitted scores are assigned to them.
Unaugmented union \(2\)-NN has exactly two components and recovers the
moons in all five lowest-noise samples, but has more than two components
and hence no feasible two-block partition in the ten other samples.
These outcomes distinguish graph feasibility from incomplete optimization.

\FloatBarrier
\subsection{UCI computation and supporting comparisons}
\label{sec:supp-uci}

Full-data optimization uses the dictionary-free evaluator described in
Appendix~\ref{subsec:lazy-scaffold-evaluation}.
Blocks are streamed rather than
collected in a global dictionary, disconnected remaining sets are factorized,
and memo entries are retained only for connected \((S,k)\) with
\(1<k<|S|\). Unless specified below, the limits are \(500{,}000\) retained plus active memo
entries, \(10^7\) generator search nodes, and a polled worker-memory ceiling
of \(2048\) MiB. Each worker has a 300-second allowance, with an internal
solver deadline of 270 seconds. Connectivity and scalar SSE caches retain
at most \(4096\) and \(50{,}000\) entries, respectively. Distance ties are
resolved by observation index, excluding self-neighbours, and the same
Euclidean MST augments each union-neighbour graph.

The Iris \(H_2\) run allows \(2{,}000{,}000\) memo entries,
\(2\times10^8\) generator nodes and 1620 solver seconds.
Wine and Breast Cancer on \(H_2\), and Vehicle on the MST, allow
\(2{,}000{,}000\) memo entries and \(5\times10^8\) generator nodes;
their solver deadlines are 108 minutes for the two \(H_2\) runs and
270 minutes for Vehicle. No reported run reaches its memory or state ceiling.

Table~\ref{tab:uci-full-work} gives the computational work and stopping
outcomes for the runs presented in Table~\ref{tab:uci-method-comparison}.
Each data-set/scaffold pair contributes one run; scores are reported only
for completed searches.

\begin{table}[!tp]
\centering
\small
\caption{Computational work for the final full-data scaffold-DP runs
reported in Table~\ref{tab:uci-method-comparison}. Memo counts retained
completed state-count pairs; Calls includes repeated and boundary evaluations.
Blocks counts emitted candidate blocks across subproblems, not globally distinct
sets, and Nodes counts generator search nodes, including unsuccessful branches.
For stopped searches these are work observed before termination, not final
support sizes. Seconds measure the solver, excluding graph construction and
worker startup.}
\label{tab:uci-full-work}
\small
\setlength{\tabcolsep}{5pt}
\begin{tabular}{lrlrrrrr}
\toprule
Data set & $q$ & Outcome & Memo & Calls & Blocks & Nodes & Seconds\\
\midrule
Iris & 0 & Exact & 78 & 23,234 & 11,103 & 22,394 & 0.215 \\
 & 2 & Exact & 5,325 & 9,624,602 & 5,056,399 & 14,132,817 & 504.050 \\
\addlinespace[2pt]
Wine & 0 & Exact & 75 & 36,243 & 15,650 & 31,529 & 0.389 \\
 & 2 & Time limit & 87,906 & 142,109,069 & 50,016,325 & 138,449,989 & 6480.018 \\
\addlinespace[2pt]
Breast cancer & 0 & Exact & 1 & 569 & 568 & 1,137 & 0.052 \\
 & 2 & Time limit & 0 & 11,187,579 & 11,187,578 & 28,531,351 & 6480.000 \\
\addlinespace[2pt]
Soybean & 0 & Exact & 92 & 37,457 & 10,395 & 21,118 & 0.143 \\
 & 2 & Exact & 252 & 888,691 & 238,581 & 537,426 & 8.560 \\
\addlinespace[2pt]
Vehicle & 0 & Exact & 2,326 & 339,818,361 & 89,657,557 & 179,423,226 & 12373.865 \\
 & 2 & Time limit & 58 & 1,223,910 & 372,203 & 775,100 & 270.000 \\
\bottomrule
\end{tabular}

\end{table}

The work counts distinguish retained storage from search effort. For
example, Vehicle's completed MST fit retains \(2326\) state-count pairs
but generates \(89{,}657{,}557\) candidate blocks. All stopped searches
reach time or search-node limits, rather than storage ceilings.
The recurrence counts feasible partitions as well as minimizing SSE;
it does not prune using an incumbent objective value.

Independent MST-cut enumeration reproduces the objectives and agreement
scores for Iris, Wine, Breast Cancer and Soybean, over \(11{,}026\),
\(15{,}576\), \(568\) and \(15{,}180\) partitions, respectively.
Vehicle's remaining-set recurrence counts \(\binom{844}{3}=99{,}846{,}044\)
partitions, as required for a tree. Independent cut enumeration is omitted
because this exceeds its \(500{,}000\)-partition ceiling, leaving no separate
enumeration check of Vehicle's minimum SSE. On Soybean and Iris, \(H_2\)
admits \(312{,}896\) and \(5{,}051{,}075\) partitions, respectively. These
are constrained partition counts, not certificates of unrestricted
optimum coverage.

\begin{table}[!tp]
\centering
\caption{Fitting times in seconds for the full-data comparators and
independent MST-cut checks used in the reported evaluation.
KM++ (1) and KM++ (20) use one and twenty initializations, respectively.
Times exclude data preprocessing and worker startup; the cut times also exclude
graph construction. These single-run, hardware-dependent measurements are
computational diagnostics, not replicated speed comparisons.}
\label{tab:uci-full-times}
\small
\setlength{\tabcolsep}{5pt}
\begin{tabular}{lrrrrr}
\toprule
Data set & Lloyd & KM++ (1) & KM++ (20) & Ward & MST cuts\\
\midrule
Iris & 0.010 & 0.000 & 0.003 & 0.021 & 0.517 \\
Wine & 0.010 & 0.000 & 0.003 & 0.022 & 0.848 \\
Breast cancer & 0.010 & 0.001 & 0.004 & 0.023 & 0.091 \\
Soybean & 0.009 & 0.000 & 0.003 & 0.020 & 0.426 \\
Vehicle & 0.036 & 0.001 & 0.008 & 0.035 & Not run \\
\bottomrule
\end{tabular}

\end{table}

The runs used Python 3.13.9, NumPy 2.3.5, SciPy 1.16.3 and
scikit-learn 1.7.2 on macOS arm64, with numerical-library thread limits set
to one per worker. Comparator fitting times are given in
Table~\ref{tab:uci-full-times}.

A separate exact-coverage study uses 20 stratified \(n=15\) samples per
data set, retaining all features and standardizing before subsampling.
Its union-neighbour graphs are connected by shortest inter-component
bridges, not by a common MST, so nesting is not assumed. The complete
recurrence supplies the unrestricted oracle for
Table~\ref{tab:uci-exact-subsamples-supp}; these results are not pooled with
the full-data comparisons.

\begin{table}[!tp]
\centering
\caption{Optimum-coverage rates (percent) on 20 full-dimensional, stratified
\(n=15\) samples.}
\label{tab:uci-exact-subsamples-supp}
\small
\begin{tabular}{lcrrr}
\toprule
Dataset & \(K\) & MST & Union \(2\)-NN & Union \(5\)-NN\\
\midrule
Iris     & 3 & 85 & 95 & 100\\
Wine     & 3 & 55 & 90 & 100\\
Breast Cancer & 2 & 65 & 95 & 100\\
Soybean    & 4 & 95 & 100 & 100\\
Vehicle    & 4 & 75 & 90 & 100\\
\bottomrule
\end{tabular}
\end{table}

Across the 100 exact samples, MST, union \(2\)-NN, and union \(5\)-NN cover
75, 94, and 100 optima. With the block dictionary's empty base entry and
reached masks both normalized by \(2^{15}\), their respective connected-block
ranges are 1.110 to 2.020\%, 4.700 to 10.290\%, and 59.560 to 85.440\%, while
reached-mask ranges are
1.230 to 2.900\%, 6.390 to 8.340\%, and 35.720 to 45.290\%.

\FloatBarrier

\subsection{Posterior protocol and validation}
\label{sec:supp-posterior-diagnostics}

The \(n=300\) design fixes component counts at
\((84,66,84,66)\), draws independently within each component, and shuffles
the observations. The means are
\((-3,-0.65)^\top\), \((-3,0.65)^\top\),
\((3,-0.65)^\top\), and \((3,0.65)^\top\), with common covariance
\(\operatorname{diag}(0.55^2,0.72^2)\).
The finite model has four available components and symmetric Dirichlet
parameter \(\gamma=1\); the number occupied may be smaller.
The DP concentration is \(\alpha=1\), and the NGGP parameters are
\((a,\sigma,\tau)=(1.5,0.25,1)\). All targets use the same collapsed Gaussian
block likelihood, with empirical prior hyperparameters
\[
\boldsymbol{\mu}_0=\bar{\mathbf{x}},\qquad
\kappa_0=0.2,\qquad
\nu_0=d+3,\qquad
\boldsymbol{\Lambda}_0=\widehat{\boldsymbol{\Sigma}}+10^{-3}\mathbf{I}_d .
\]
The scaffold proposal uses \(\kappa=2\) in
\eqref{eq:scaffold-restricted-gibbs}, \(\epsilon=0.05\) in
\eqref{eq:scaffold-split-anchor-weight} to
\eqref{eq:scaffold-merge-anchor-weight}, and distance offset
\(\epsilon_D=10^{-12}\). Union nearest-neighbour edges are augmented by the
Euclidean MST computed from all pairwise distances. The \(2\)-NN-MST and
\(3\)-NN-MST graphs have 421 and 584 edges, respectively.

The four evaluation chains start at one cluster, all singletons, one-start
\(k\)-means++, and a random four-group allocation. For the finite target,
a cyclic allocation to its four available components replaces the
inadmissible singleton initialization. Initial partitions are matched across
samplers within each target. Each NGGP chain starts at \(u=1\) and uses
the log-scale slice update derived in Appendix~\ref{supp:auxiliary-update},
with fixed initial bracket width one. The update is identical across all
NGGP samplers. After 1,000 warm-up sweeps, each evaluation chain retains
2,500 states and each reference chain retains 20,000, without thinning;
move counters are reset at the warm-up boundary. The finite and DP results
use saved archives with 100 warm-up sweeps. These archives were collected
under limits of 60 seconds per evaluation chain and 120 seconds per reference
chain, but comparisons use only common iteration prefixes. Their counters
first appear at sweep 101, so we subtract that boundary and analyse the
next 500, 1,000, or 2,500 states. Every archive contains the required prefix;
no chains are omitted. Warm-up and evaluation prefixes are matched across
samplers within each target, not across targets. Equal retained sweeps
match allocation updates and split-merge opportunities, not arithmetic work.

Each target has eight separate reference chains, four generic and four
\(3\)-NN-guided, with independent random seeds and the same four initialization
types. Diagnostics use \(K\), log posterior, VI to generating labels, and
eight pairwise co-clustering indicators selected from the reference PSM
closest to probability \(1/2\); NGGP additionally includes \(u\) and \(\log u\).
Globally constant statistics are excluded. We compute rank-normalized split
\(\widehat R\) and bulk effective sample size (ESS) for each statistic.
Reference validation requires \(\max\widehat R\le1.01\), minimum bulk ESS
at least 100 for finite and DP mixtures and 1,000 for NGGP, and at least
100 retained draws per chain, together with PSM RMSE at most \(0.03\) and
total variation for \(K\) at most \(0.05\) between the generic and guided
reference groups. NGGP also requires tail ESS at least 1,000 for both
\(u\) and \(\log u\), taking the smaller ESS of indicators below their
pooled 5th and 95th percentiles. These thresholds were specified before
their respective runs. The finite and DP references have maximum \(\widehat R\)
\(1.002\) and \(1.003\), minimum ESS 4,102 and 11,881, and cross-group
PSM RMSE \(2.200\times10^{-3}\) and \(9.000\times10^{-4}\), respectively.
The NGGP reference has maximum \(\widehat R=1.000\), minimum bulk ESS
49,418, and auxiliary bulk and tail ESS 138,735 and 112,221, respectively.
Its cross-group PSM RMSE is \(6.310\times10^{-4}\) and total variation for \(K\) is
\(5.610\times10^{-3}\). All reference checks pass; these diagnostics support the
comparison but do not certify convergence.

For four evaluation chains of length \(T\), let
\(\widehat P_{ij}=(4T)^{-1}\sum_{c=1}^4\sum_{t=1}^T
\mathbf{1}\{i\sim_{\ell_{ct}}j\}\), and let \(\widehat p_k\) be the analogous
frequency of \(K=k\). With reference quantities denoted by superscript
\({\rm ref}\), the accuracy measures are
\[
{\rm RMSE}
=\left\{\binom n2^{-1}\sum_{i<j}
(\widehat P_{ij}-\widehat P^{\rm ref}_{ij})^2\right\}^{1/2},
\qquad
{\rm TV}_K=\tfrac12\sum_k|\widehat p_k-\widehat p^{\rm ref}_k|.
\]
The reference PSM averages the eight chain-specific PSMs with equal weight,
using every retained draw. NGGP reference \(K\) frequencies also use every
draw; finite and DP reference \(K\) frequencies use 500 evenly spaced draws
per chain. All reference VI risks use 500 draws per chain. VI estimation uses
300 evenly spaced draws per evaluation chain. Candidate sets combine
the 64 most frequent reference partitions with average-linkage cuts of
\(1-\widehat{\mathbf P}^{\rm ref}\), then remove duplicates; this gives
68, 69, and 71 candidates for finite, DP, and NGGP mixtures, respectively.
Writing these sets as \(\mathcal A\), the reported excess risk is
\(\widehat r_{\rm ref}(\widehat\ell_T)
-\min_{\ell\in\mathcal A}\widehat r_{\rm ref}(\ell)\), where
\(\widehat\ell_T=\arg\min_{\ell\in\mathcal A}\widehat r_T(\ell)\)
and each \(\widehat r\) is the corresponding Monte Carlo average of VI
in bits. No Jensen-bound approximation is used for this finite-set
minimization. Zero excess risk therefore means agreement with the reference
candidate-set decision, not exact minimization over all partitions.
For NGGP, all methods select the same three-cluster decision at every
checkpoint; finite and DP methods agree with it at 1,000 and 2,500 sweeps.
At 500 sweeps the two guided DP methods instead select two clusters, with
additional reference VI risk \(1.900\times10^{-3}\) bits. The NGGP reference PSM has
estimated root-mean-square Monte Carlo standard error \(3.990\times10^{-4}\), obtained
from variation across its eight independent chain means. This measures
reference precision, not a bias bound or uncertainty across data sets, and
does not resolve very small differences between evaluation methods.

\begin{table}[tbp]
\centering
\caption{The \(n=300\) comparison at 2,500 retained sweeps per chain.
Acceptance pools 10,000 attempted split-merge moves across four chains.
PSM RMSE and total variation for \(K\) compare with the independent
Monte Carlo reference; smaller is better. All three references pass their
diagnostic checks. The NGGP rows use the log-scale slice update throughout.
Gibbs has no split-merge acceptance rate.}
\label{tab:posterior-scaling-diagnostics}
\small
\setlength{\tabcolsep}{4pt}
\begin{tabular}{llrrrrr}
\toprule
Target & Sampler & $\max\widehat R$ & Min. ESS & Accept. (\%) & PSM RMSE & TV for $K$\\
\midrule
Finite & Gibbs & 1.008 & 808 & n/a & 2.029e-03 & 1.495e-02 \\
 & Generic & 1.002 & 1218 & 1.950 & 2.045e-03 & 3.450e-03 \\
 & Guided $q=2$ & 1.007 & 948 & 2.290 & 1.771e-03 & 1.600e-03 \\
 & Guided $q=3$ & 1.004 & 1274 & 2.540 & 5.798e-03 & 9.700e-03 \\
 & Permuted $q=3$ & 1.006 & 665 & 1.560 & 3.906e-03 & 1.410e-02 \\
\midrule
DP & Gibbs & 1.010 & 1759 & n/a & 1.092e-03 & 1.215e-02 \\
 & Generic & 1.006 & 2671 & 1.090 & 1.519e-03 & 1.310e-02 \\
 & Guided $q=2$ & 1.002 & 2178 & 1.410 & 1.038e-03 & 2.205e-02 \\
 & Guided $q=3$ & 1.006 & 1980 & 1.860 & 8.104e-04 & 1.560e-02 \\
 & Permuted $q=3$ & 1.002 & 2114 & 0.890 & 2.181e-03 & 1.375e-02 \\
\midrule
NGGP & Gibbs & 1.001 & 3098 & n/a & 1.283e-03 & 1.032e-02 \\
 & Generic & 1.000 & 3407 & 2.250 & 1.931e-03 & 1.246e-02 \\
 & Guided $q=2$ & 1.001 & 3119 & 3.580 & 1.300e-03 & 1.328e-02 \\
 & Guided $q=3$ & 1.001 & 3562 & 4.330 & 1.303e-03 & 8.400e-03 \\
 & Permuted $q=3$ & 1.001 & 3399 & 2.010 & 1.317e-03 & 1.013e-02 \\
\bottomrule
\end{tabular}

\end{table}

All 15 target-sampler combinations pass the evaluation thresholds
\(\max\widehat R\le1.05\) and minimum bulk ESS at least 100 at each of the
three checkpoints; Table~\ref{tab:posterior-scaling-diagnostics} reports the
final checkpoint. No uncertainty
intervals across data sets are reported for this single simulated case.

All samplers use the same bounded cache of sufficient statistics and
collapsed likelihoods for at most 8,192 visited blocks. Eviction changes
recomputation, not the likelihood or target; no subset-DP solution is cached.
Regression checks compare cached and uncached Gibbs updates, proposal
probabilities, and Metropolis ratios for every admissible partition of five
observations under all three targets, including forced cache eviction.
Exhaustive four-observation transition checks for the cached generic and
guided kernels give detailed-balance, stationarity, and row-sum errors below
\(10^{-12}\). The NGGP check conditions on \(u=1.7\), so it validates the
partition subkernel rather than the continuous auxiliary update. Separate
slice-update checks compare its conditional distribution with numerical
quadrature and, at \(\tau=0\), with the identity
\((a/\sigma)U^\sigma\mid\ell\sim\operatorname{Gamma}(K,1)\).
These are implementation checks, not additional sampling-performance experiments.

\FloatBarrier

\section{Implementation details}
\label{supp:implementation}

\subsection{Component-factorized evaluation of remaining-set problems}
\label{subsec:component-factorized-evaluation}

The dictionary-based factorization ablation in
Appendix~\ref{sec:supp-frontier-pilot} uses the first-block rule on connected
remaining sets and component convolution on disconnected remaining sets. With components
\(S_1,\ldots,S_c\), initialize \(D_0(0)=0\) and \(D_0(t)=+\infty\) otherwise,
then compute
\[
D_j(t)=\min_{1\le r\le\min\{|S_j|,t\}}
\{D_{j-1}(t-r)+F_H(S_j,r)\},\qquad F_H(S,k)=D_c(k).
\]
Unlike the selective-storage implementation in
Appendix~\ref{subsec:lazy-scaffold-evaluation}, this ablation stores both
component and union solutions. Algorithm~\ref{alg:lazy-scaffold-evaluation} gives the
dictionary-free evaluator used for full-data UCI fits.

The component decomposition depends only on \(S\), so it is cached by the
mask encoding that set rather than by \((S,k)\). Subproblem solutions remain
indexed by \((S,k)\): the optimum on a component depends on neither the
larger remaining set from which it was obtained nor the other components. Recursion is acyclic because every component at a
genuine split and every set left after a block removal is a proper subset of
the current remaining set. In the penalized version, recursively evaluate every
component with the penalized hybrid routine, sum the returned
\(F_{H,\lambda}(S_j)\), memoize that sum at \(S\), and record its decomposition
into \(S_1,\ldots,S_c\). With \(n_j=|S_j|\), a direct
fixed-count convolution truncated at \(k\) uses at most
\[
O\left(k\sum_{j=1}^c\min\{n_j,k\}\right)
\subseteq O(ck^2)
\]
min-plus relaxations beyond the memoized component solves; the penalized
combination uses \(O(c)\) additions. Finding the components costs
\(O\{|S|+|E(H[S])|\}\) once per cached mask. The complete-graph worst case is
unchanged, because no nontrivial component split then occurs.

\begin{samepage}
The split rule must \emph{replace}, rather than supplement, the loop over first-block choices
at a disconnected state. Either rule alone represents every unlabeled
partition once, but adding their outputs would double count under counting,
enumeration, or tied-optimum arithmetic. For
\(N_H(S,k):=|\mathcal P_{H,k}(S)|\), the correct component rule is
\begin{equation}
\label{eq:component-count-convolution}
N_H(S,k)
=
\sum_{\mathbf{k}\in\mathcal K_H(S,k)}
\prod_{j=1}^c N_H(S_j,k_j),
\end{equation}
with no multinomial factor for interleavings.
The stored component decomposition allows each selected component partition
to be reconstructed separately and their blocks combined; sorting these blocks by their least vertices restores
the order imposed by the first-block rule when a canonical ordered representation is
needed. To enumerate every tied optimum, the first-block and convolution steps
must retain all minimizing block and allocation predecessors, respectively.
\end{samepage}

Finally, component calls can create masks that are not members of the
block-removal closure \(\mathcal R_{H,K}\) in
\eqref{eq:reachable-state-closure}. Hence \(\mathcal R_{H,K}\) and
\(T_K^{\mathrm{reach}}(H)\) continue to describe the basic subset-DP
evaluator, but need not equal the stored states and operations of the hybrid
evaluator. Nor does the nesting in
\eqref{eq:monotone-state-partition} automatically transfer to its recursion
DAG: on two isolated vertices with \(K=2\), for example, a component split may
evaluate both singleton masks, whereas after the connecting edge is added the
remaining-set recurrence reaches only the singleton that excludes
the least-indexed observation. A benchmark of the hybrid evaluator should therefore report its
component splits and convolution relaxations separately rather than reuse the
transition count for the basic remaining-set recurrence.

\subsection{Dictionary-free evaluation, block generation, and caching}
\label{subsec:lazy-scaffold-evaluation}

The small-sample evaluator in Algorithm~\ref{alg:scaffold-subset-dp} builds
its dictionary by adjacent-vertex additions and duplicate removal; a rooted
spanning tree provides such an ordering for every connected set. With
\(a=\min S\), it scans the smaller of the indexed list
\(\mathcal C_a(H)=\{B\in\mathcal C(H):\min B=a\}\) and all subsets of \(S\)
containing \(a\), checking containment, connectivity, and cardinality.
The index reduces lookup work but still requires the dictionary in advance.
For full-data evaluation, we instead generate feasible blocks only when
needed by the current remaining-set problem.

Write \(c_H(R)\) for the number of components of \(H[R]\), with
\(c_H(\emptyset)=0\). For connected \(H[S]\) and \(1<k<|S|\),
the generator emits exactly the family
\begin{equation}
\label{eq:lazy-extendable-blocks}
\mathcal G_H(S,k)=\{B\subseteq S:\min S\in B,\ H[B]\text{ connected},\
c_H(S\setminus B)\le k-1\le |S\setminus B|\}.
\end{equation}
After choosing \(B\), the remaining set \(R=S\setminus B\) must form
\(k-1\) further blocks. A connected block cannot span distinct components
of \(H[R]\), so \(c_H(R)\le k-1\). Since the blocks are disjoint and
nonempty, \(k-1\le |R|\). These are the two inequalities in
\eqref{eq:lazy-extendable-blocks}. They also suffice: cutting
\(k-1-c_H(R)\) edges from spanning trees of the components of \(H[R]\)
produces \(k-1\) connected blocks
(Appendix~\ref{app:proof-component-factorization}). Thus
\(\mathcal G_H(S,k)\) contains exactly the blocks satisfying the first-block rule with a feasible
continuation.

Algorithm~\ref{alg:lazy-scaffold-evaluation} evaluates
\(F_H(S,k)\), the minimum SSE for partitioning the remaining observation
indices \(S\) into exactly \(k\) nonempty, scaffold-connected blocks, with
\(H\) and the block cost fixed. Subproblems are evaluated and blocks
generated only as needed, without changing the exact objective.
A finite return also carries a minimizing partition; an infeasible call
returns \(+\infty\). Arithmetic uses only the objective component, and the
initial call \((V,K)\) yields \(\mathrm{OPT}(H,K)\) if the search completes.

\begin{algorithm}[!tbp]
\small
\caption{Dictionary-free fixed-count evaluator with selective memoization}
\label{alg:lazy-scaffold-evaluation}
\begin{algorithmic}[1]
\Require scaffold \(H\), scalar block cost \(\phi_{\mathrm{SSE}}\), target \(K\)
\State initialize empty solution table \(M\) and bounded component and cost caches
\Function{$F_H$}{$S,k$}
 \State check the computational budget; if exhausted, stop as incomplete
 \If{\((S,k)\in M\)} \State \Return \(M[S,k]\) \EndIf
 \If{\(S=\emptyset\)} \State \Return \(0\) if \(k=0\), else \(+\infty\) \EndIf
 \If{\(k\le0\) or \(k>|S|\)} \State \Return \(+\infty\) \EndIf
 \State identify the connected components \(S_1,\ldots,S_c\) of \(H[S]\)
 \If{\(c>k\)} \State \Return \(+\infty\) \EndIf
 \If{\(k=|S|\)} \State \Return \(0\) and the singleton partition \EndIf
 \If{\(c>1\)}
  \State convolve \(F_H(S_j,r)\) over component counts
  using \eqref{eq:fixed-component-factorization}
  \State \Return the minimum and a minimizing partition, without storing \((S,k)\)
 \EndIf
 \If{\(k=1\)} \State \Return \(\phi_{\mathrm{SSE}}(S)\) and the single block \(S\) \EndIf
 \State reserve one active table entry; set \(v\gets+\infty\)
 \ForAll{\(B\) streamed from \(\mathcal G_H(S,k)\)}
  \State \(u\gets\phi_{\mathrm{SSE}}(B)+F_H(S\setminus B,k-1)\)
  \State if \(u<v\), set \(v\gets u\) and retain the corresponding minimizing partition
 \EndFor
 \State release the active slot and store the completed solution in \(M[S,k]\)
 \State \Return the stored solution
\EndFunction
\State \Return \(F_H(V,K)\)
\end{algorithmic}
\end{algorithm}

\begin{samepage}
To generate \(\mathcal G_H(S,k)\) without listing connected sets in advance,
we decide which observations belong to the next block. Fix \(a=\min S\)
and \(r=k-1\). At any stage, \(B\) contains the vertices already included
and \(E\) those already excluded; initially \(B=\{a\}\) and
\(E=\emptyset\). The set \(B\) remains connected, and a possible final
block \(B'\) must satisfy \(B\subseteq B'\subseteq S\setminus E\).\par
\end{samepage}

The feasibility conditions in \eqref{eq:lazy-extendable-blocks} allow some
decisions to be made without branching. Put \(R=S\setminus B\). If
\(|R|<r\), no extension of \(B\) can leave enough observations for the
remaining clusters, so this search branch is discarded. A vertex
unreachable from \(a\) in \(H[S\setminus E]\) cannot belong to a connected
extension avoiding \(E\), so it is also excluded. After these exclusions,
let \(p\) be the number of components of \(H[R]\) that intersect \(E\).
Each contains a vertex that must remain outside \(B'\), and removing more
vertices from \(R\) cannot join distinct components. Every continuation
therefore leaves at least \(p\) connected components, ruling it out if \(p>r\).

If \(p=r\), these components already require all \(r\) remaining blocks.
Any component of \(H[R]\) not intersecting \(E\) must then be included in
\(B\): leaving even one of its vertices outside would create another
component outside \(B\). This inclusion preserves connectivity because, when
\(H[S]\) and \(H[B]\) are connected, every component of \(H[R]\) has an
edge to \(B\). We update \(R\) and repeat the size check after these forced
inclusions. Each of these rules consequently preserves every feasible final
block consistent with the decisions already made.

\begin{samepage}
When \(|R|=r\), the current \(B\) is passed to the recurrence: its complement
can be partitioned into \(r\) singletons, and any further inclusion would
leave too few observations. Otherwise, the undecided neighbours of \(B\)
form the frontier \((N_H(B)\cap R)\setminus E\). If this set is empty,
there is no larger connected extension avoiding \(E\), since any path to an
added vertex would need a first step through the frontier. We then pass
\(B\) to the recurrence exactly when \(c_H(R)\le r\).\par
\end{samepage}

For a nonempty frontier, choose its least-indexed vertex \(v\) and explore the two
possibilities: include \(v\) in \(B\), or add it to \(E\). These choices
divide the possible final blocks into disjoint families. Every branch fixes
another vertex, so the search terminates; every feasible block follows one
unique sequence of choices and survives the feasibility checks above.
Together with the terminal tests, this proves that the generator returns
every member of \(\mathcal G_H(S,k)\) exactly once. A depth-first traversal
passes blocks directly to the recurrence, without keeping a list of previous
outputs for duplicate removal.

\begin{samepage}
If \(H[S]\) is a tree, root it at the least-indexed observation \(a=\min S\). For a vertex \(v\)
to belong to a connected block containing \(a\), every vertex on the unique
path from \(a\) to \(v\) must also belong to that block. Starting from
\(B=\{a\}\), consider an undecided vertex \(v\) whose parent is already
in \(B\). Including \(v\) preserves connectivity and allows each of its
children to be considered next, without yet deciding whether they belong to
\(B\). Excluding \(v\) forces its entire subtree to remain outside
\(B\): every path from \(a\) to a descendant of \(v\) passes through
\(v\). These excluded observations are not discarded from the clustering
problem; they remain available for the other \(r=k-1\) blocks.\par
\end{samepage}

The excluded subtrees are disjoint. Each attaches to \(B\) only through
the parent of its root, so each is a separate component of
\(H[S\setminus B]\). Each such component may subsequently be split into several
clusters; it need not be a single cluster. Thus excluding \(t\) subtrees
requires at least \(t\) remaining blocks, giving \(t\le r\). Once
\(t=r\), every undecided vertex must be included in \(B\), since any
further exclusion would create an additional component of \(H[S\setminus B]\). When all
decisions have been made, \(c_H(S\setminus B)=t\), and the excluded
vertices must also satisfy \(r\le |S\setminus B|\). These are exactly
the component and size conditions in \eqref{eq:lazy-extendable-blocks}.
Every connected block satisfying the first-block rule determines
a unique set of excluded subtree roots whose parents are in the block, so
this generator is also exhaustive and duplicate-free. Only block generation
changes; the objective and component recurrence remain the same. Both
generators have depth at most \(|S|\) and retain \(O(|S|)\) pending search
records, each containing vertex masks. This bounds the number of records,
not their storage in constant-size words independently of \(|S|\).

The same evaluation also counts feasible partitions. Counts are added over
first-block choices. For a disconnected remaining set with components
\(S_1,\ldots,S_c\), they are combined by
\eqref{eq:component-count-convolution}. Suppose the first \(j-1\) components
have used \(t\) blocks and component \(S_j\) receives \(r_j\). Its capacity requires
\(1\le r_j\le |S_j|\), while the remaining components require
\(c-j\le k-t-r_j\le\sum_{h>j}|S_h|\). Combining these inequalities gives
\(\max\{1,k-t-\sum_{h>j}|S_h|\}\le r_j\le
\min\{|S_j|,k-t-(c-j)\}\). Restricting the convolution to this range
omits no feasible allocation.

Generating blocks on demand avoids storing a global dictionary, but the
recursion can still request the solution of the same remaining-set problem repeatedly. To reuse
these calculations, the solution table stores completed subproblem solutions indexed by
\((S,k)\), restricted to connected \(H[S]\) with \(1<k<|S|\).
Each entry contains the objective, a minimizing tuple of blocks, and the
feasible-partition count. For a disconnected remaining set, the minimizing
component partitions are combined and returned from a temporary convolution;
one-cluster and singleton-partition cases are obtained directly. These
solutions need no separate persistent table entry, and reconstruction does
not require references to subproblems whose solutions were never stored.

Two smaller caches reuse auxiliary calculations: component decompositions
are indexed by \(S\), independently of \(k\), and SSE values by the block
mask \(B\). Each cache has a fixed entry limit and discards the least
recently used entry when full. A later request for a discarded value
recomputes the same quantity. To avoid subtracting large raw second moments,
SSE is evaluated after translating the observations by one block member and
centering at the translated block mean. The routine holds a private copy of
the observations so that external edits cannot invalidate these cached costs.

These storage choices do not change the recurrence values. For an induction
on \(|S|\), the empty and other boundary cases agree with the original
recurrence. A first-block step leaves a proper subset \(S\setminus B\),
and a genuine component split uses proper subsets \(S_j\); their values
therefore agree by the induction hypothesis. Since the generator and
convolution consider exactly the feasible alternatives, both the minimum
and the sums and products defining the partition count agree as well.
Recomputing a value instead of retaining it can increase work, but cannot
change either result.

Neither on-demand generation nor selective storage removes the possibility
of exponential computation and memory use. Before evaluating a subproblem
that will be memoized, the implementation reserves an entry for its eventual
solution; the state limit counts both completed entries and these active
reservations. Recurrence calls, emitted blocks and generator search nodes
are recorded separately because repeated calculations contribute work without
adding distinct stored states. The auxiliary-cache limits likewise control
entry counts rather than total process memory. If a state, search-node or
time limit is reached, the search is reported as incomplete, with no fitted
partition or exactness claim. Appendix~\ref{sec:supp-uci} gives the limits
and external worker-memory checks used in the UCI evaluation.

\section{Mixture priors and proposal specification}
\label{app:posterior-proofs}

\subsection{Collapsed targets and proposal specification}
\label{app:posterior-targets}

An exchangeable partition probability function (EPPF) assigns to each
particular partition of \(n\) observations the probability
\(p(n_1,\ldots,n_K)\), where \(n_j\) is the size of its \(j\)th block.
Exchangeability makes this probability depend only on the block sizes,
and \(p\) is symmetric in its arguments. It is not the aggregate
probability of all partitions with those sizes.
For a partition \(\ell\) with \(K\) blocks and block sizes \(n_A\), the three
partition priors used in the experiments have, up to factors independent of
\(\ell\), the following forms:
\begin{align}
p_{\mathrm F}(\ell)
&\propto
(L)_K\prod_{A\in\ell}\frac{\Gamma(n_A+\gamma)}{\Gamma(\gamma)},
&&1\le K\le L,
\label{eq:finite-partition-prior-app}\\
p_{\mathrm{DP}}(\ell)
&\propto
\alpha^K\prod_{A\in\ell}\Gamma(n_A),
\label{eq:dp-partition-prior-app}\\
p_{\mathrm{NGGP}}(\ell,u)
&\propto
u^{n-1}
\exp\left[-\frac a\sigma\{(u+\tau)^\sigma-\tau^\sigma\}\right]
a^K
\prod_{A\in\ell}
\frac{\Gamma(n_A-\sigma)}
   {\Gamma(1-\sigma)(u+\tau)^{n_A-\sigma}},
&&u>0.
\label{eq:nggp-partition-u-prior-app}
\end{align}
Here \((L)_K=L!/(L-K)!\), \(L\ge1\), \(\gamma>0\), \(\alpha>0\),
\(a>0\), \(0<\sigma<1\), and \(\tau\ge0\). The finite expression follows
by integrating symmetric Dirichlet weights and summing over the \((L)_K\)
injective assignments of occupied unlabeled blocks to component labels. The
DP expression is proportional to the Chinese-restaurant EPPF, whose
normalizing factor is \(\Gamma(\alpha)/\Gamma(\alpha+n)\), and
\eqref{eq:nggp-partition-u-prior-app} is the auxiliary-variable joint density
used by the implementation \citep{neal2000markov,favaro2013mcmc}.

During a restricted split, both subclusters have already been seeded; assigning
one observation to an existing subcluster of current size \(s\) changes the
corresponding block factor by
\begin{align*}
\frac{\Gamma(s+1+\gamma)}{\Gamma(s+\gamma)}&=s+\gamma,
&
\frac{\Gamma(s+1)}{\Gamma(s)}&=s,
&
\frac{\Gamma(s+1-\sigma)}{\Gamma(s-\sigma)}
\frac1{u+\tau}&=\frac{s-\sigma}{u+\tau}.
\end{align*}
Thus the prior factors in \eqref{eq:restricted-allocation-base} are
\begin{equation}
c_M(s;u)=
\begin{cases}
s+\gamma,&M=\text{finite mixture},\\
s,&M=\text{DP mixture},\\
(s-\sigma)/(u+\tau),&M=\text{NGGP}.
\end{cases}
\label{eq:restricted-allocation-cohesion}
\end{equation}
The parameters are the finite symmetric Dirichlet mass \(\gamma>0\), DP
concentration \(\alpha>0\), and NGGP mass \(a>0\), discount \(0<\sigma<1\),
and tilt \(\tau\geq0\). New-component factors do not enter this allocation
because both subclusters are already seeded.

For completeness, the collapsed Gaussian likelihood uses
\[
\boldsymbol{\Sigma}\sim\operatorname{IW}_d(\nu_0,\boldsymbol{\Lambda}_0),
\qquad
\boldsymbol{\mu}\mid\boldsymbol{\Sigma}\sim N_d(\boldsymbol{\mu}_0,\boldsymbol{\Sigma}/\kappa_0),
\qquad
\mathbf{x}_i\mid\boldsymbol{\mu},\boldsymbol{\Sigma}\sim N_d(\boldsymbol{\mu},\boldsymbol{\Sigma}),
\]
with \(\kappa_0>0\), \(\nu_0>d-1\), and \(\boldsymbol{\Lambda}_0\succ0\). For a nonempty block
\(A\), put
\begin{align*}
n_A&=|A|,&
\bar{\mathbf{x}}_A&=n_A^{-1}\sum_{i\in A}\mathbf{x}_i,&
\mathbf{S}_A&=\sum_{i\in A}(\mathbf{x}_i-\bar{\mathbf{x}}_A)(\mathbf{x}_i-\bar{\mathbf{x}}_A)^\mathsf T,\\
\kappa_A&=\kappa_0+n_A,&
\nu_A&=\nu_0+n_A,&
\boldsymbol{\Lambda}_A&=\boldsymbol{\Lambda}_0+\mathbf{S}_A+
\frac{\kappa_0n_A}{\kappa_A}
(\bar{\mathbf{x}}_A-\boldsymbol{\mu}_0)(\bar{\mathbf{x}}_A-\boldsymbol{\mu}_0)^\mathsf T.
\end{align*}
The standard normal-inverse-Wishart marginal likelihood is
\citep{murphy2012machine}
\begin{equation}
m(\mathbf{X}_A)
=
\pi^{-n_Ad/2}
\left(\frac{\kappa_0}{\kappa_A}\right)^{d/2}
\frac{\Gamma_d(\nu_A/2)}{\Gamma_d(\nu_0/2)}
\frac{|\boldsymbol{\Lambda}_0|^{\nu_0/2}}{|\boldsymbol{\Lambda}_A|^{\nu_A/2}},
\label{eq:niw-block-marginal-app}
\end{equation}
where \(\Gamma_d\) is the multivariate gamma function
\citep{murphy2012machine}. Positive-definiteness
of \(\boldsymbol{\Lambda}_0\) makes every determinant and block marginal positive and
finite. This formula makes the ratios in
\eqref{eq:restricted-allocation-base} reproducible; \(\boldsymbol{\Lambda}_0\) is the
inverse-Wishart scale matrix, not in general its prior mean.

Assume that the proposal scaffold \(H\) is connected, as is the
\(q\)-NN-MST graph in Section~\ref{sec:exp-posterior}. For \(n\ge2\), write \(D_H(i,j)\) for
Euclidean-edge-weighted shortest-path distance and \(\bar D_H\) for its median
over distinct unordered vertex pairs. To cover coincident observations or a
zero median, fix \(\epsilon_D>0\) and put
\(D_H^\star=\max(\bar D_H,\epsilon_D)\). Conditional on a split, eligible
unordered seed pairs \(\{i,j\}\) satisfy \(i\sim_\ell j\); conditional on a merge
they satisfy \(i\not\sim_\ell j\). Every selected pair is then oriented by
the convention \(i<j\): \(i\) initializes subcluster \(R_1\), \(j\)
initializes subcluster \(R_2\), and this orientation is retained by the reverse move.
The scaffold proposal
normalizes the weights
\begin{align}
w_{\mathrm S}(i,j\mid\ell)
&\propto
\epsilon+\max\{D_H(i,j)/D_H^\star,10^{-8}\},
&& i\sim_\ell j,
\label{eq:scaffold-split-anchor-weight}\\
w_{\mathrm M}(i,j\mid\ell)
&\propto
\epsilon+\exp\{-D_H(i,j)/D_H^\star\},
&& i\not\sim_\ell j,
\label{eq:scaffold-merge-anchor-weight}
\end{align}
with \(\epsilon>0\). A split is disabled at \(|\ell|=L\) for an
\(L\)-component finite mixture. Otherwise the move-type probabilities are
\(p_{\mathrm S}=p_{\mathrm M}=1/2\) when both operations are available and
one for the sole available operation. If neither operation is available
(in particular, when \(n=1\), or when \(L=1\) and \(|\ell|=1\)), the proposal
is the identity kernel and no seed pair is selected and no distance scale is evaluated. After
choosing the seed pair \((i,j)\), a
scaffold split visits
\(A_\ell(i)\setminus\{i,j\}\), the remaining observations in their common
parent block (with \(A_\ell(i)\) denoting the block containing \(i\)), in increasing
\(\min\{D_H(v,i),D_H(v,j)\}\), with index tie-breaking, and applies
\eqref{eq:scaffold-restricted-gibbs}. The distance ordering and index
tie-breaking fix the same observation sequence for the forward split and the
forced reverse allocation, allowing both proposal probabilities to be evaluated.

Writing \(o\in\{\mathrm S,\mathrm M\}\) for move type and
\(\mathbf b_{1:m}\) for a split's sequential allocation decisions, let
\(\omega=(o,i,j,\mathbf b_{1:m})\), with no allocation vector for a merge.
If \(p_o\) denotes move-type probability and \(w_o\) the normalized seed-pair
probability, the history mass used by the acceptance ratio is
\begin{equation}
\label{eq:auxiliary-path-proposal-mass}
\widetilde q_H(\omega\mid\ell,u)
=
p_o(\ell)w_o(i,j\mid\ell)
\begin{cases}
\prod_{t=1}^m q_{H,\kappa}
\{v_t\to b_t\mid R_{1,t-1},R_{2,t-1},u\},&o=\mathrm S,\\
1,&o=\mathrm M.
\end{cases}
\end{equation}

\subsection{NGGP auxiliary-variable update}
\label{supp:auxiliary-update}

For the NGGP auxiliary update, put \(\eta=\log u\). The conditional density
with respect to \(d\eta\) is proportional to
\(p_{\mathrm{NGGP}}(\ell,e^\eta)e^\eta\), since the collapsed likelihood
does not involve \(u\). Taking logarithms in
\eqref{eq:nggp-partition-u-prior-app} gives, up to an additive constant,
\begin{equation}
\label{eq:nggp-log-u-conditional}
h(\eta)=n\eta-(n-K\sigma)\log(e^\eta+\tau)
-\frac a\sigma(e^\eta+\tau)^\sigma.
\end{equation}
The coefficient \(n\), rather than \(n-1\), includes the transformation
Jacobian. With \(u=e^\eta\), differentiation yields
\[
h''(\eta)=-(n-K\sigma)\frac{\tau u}{(u+\tau)^2}
-a u(u+\tau)^{\sigma-2}(\tau+\sigma u)<0.
\]
Thus \(h\) is strictly concave and tends to \(-\infty\) in both tails
under the stated parameter assumptions. The implementation draws
\(E\sim\operatorname{Exp}(1)\), sets \(b=h(\eta)-E\) and \(y=e^b\), and applies
stepping-out and shrinkage slice sampling \citep{neal2003slice}.
A randomly offset bracket of width one is expanded until both endpoints
lie outside \(\{s:h(s)\ge b\}\), a bounded interval. Uniform proposals
from the bracket are accepted inside this interval; a rejected proposal
replaces the endpoint on its side of the current \(\eta\). No part of the
slice is removed, so the accepted point is uniform on the slice. Alternating
this update with the height draw preserves the joint density proportional
to \(\mathbf{1}\{0<y<e^{h(\eta)}\}\), whose \(\eta\)-marginal is the
required conditional. A numerical evaluation limit raises an error rather
than returning an approximate update.
If an allocation kernel, this auxiliary kernel, and the split-merge kernel each preserve the
joint target, then their composition does too, since
\(\pi P_1P_2P_3=((\pi P_1)P_2)P_3=\pi\). Irreducibility and aperiodicity are
additional properties of the full composition, not consequences of the
scaffold-guided subkernel alone.

\section{Segmentation and geometric background}
\label{app:preliminaries}

\subsection{One-dimensional segmentation}
\label{app:bellman-1d}

For sorted scalar observations \(x_{[1]}\le\cdots\le x_{[n]}\), an
unrestricted \(K\)-cluster SSE optimum can be chosen with contiguous blocks
\citep{wang2011ckmeans,nielsen2014optimal}. Writing \(e(i,j)\) for the SSE
of interval \([i{:}j]\), the standard segmentation recurrence is
\[
D(j,k)=\min_{k\le i\le j}\{D(i-1,k-1)+e(i,j)\},
\qquad D(0,0)=0,
\]
with other infeasible boundary states set to \(+\infty\)
\citep{bellman1961approximation,bellman1969curve}.
It evaluates \(O(nK)\) prefix states with at most \(n\) possible final
intervals per state, despite there being \(2^{n-1}\) unrestricted
segmentations. In Section~\ref{sec:subset-lattice}, the naturally ordered
path yields the suffix version of exactly this recurrence. The path has
\(n(n+1)/2\) connected blocks; materializing these blocks is distinct from
storing the \(O(nK)\) DP values.

\subsection{SSE identities}
\label{app:sse-identities}

For a nonempty indexed block \(A\), let
\[
m_A=|A|,\qquad \mathbf{s}_A=\sum_{i\in A}\mathbf{x}_i,\qquad
q_A=\sum_{i\in A}\|\mathbf{x}_i\|^2,\qquad \boldsymbol{\mu}_A=\mathbf{s}_A/m_A.
\]
Expanding the square and using \(\sum_{i\in A}(\mathbf{x}_i-\boldsymbol{\mu}_A)=\mathbf{0}\) gives
\[
\sum_{i\in A}\|\mathbf{x}_i-\mathbf{c}\|^2
=\mathrm{SSE}(A)+m_A\|\mathbf{c}-\boldsymbol{\mu}_A\|^2,
\qquad
\mathrm{SSE}(A)=q_A-\|\mathbf{s}_A\|^2/m_A.
\]
For disjoint nonempty blocks \(A,B\), the corresponding ANOVA identity is
\begin{equation}
\mathrm{SSE}(A\cup B)
=\mathrm{SSE}(A)+\mathrm{SSE}(B)
+\frac{m_Am_B}{m_A+m_B}\|\boldsymbol{\mu}_A-\boldsymbol{\mu}_B\|^2.
\label{eq:anova-merge-app}
\end{equation}
The last term is the standard Ward merge cost
\citep{ward1963hierarchical,murtagh2014ward}. The first identity proves that
block means minimize assigned SSE; the second gives the cached subset weight
used by the recurrence. These uses require no hierarchical merge path.

\subsection{Voronoi cells and Delaunay adjacency}
\label{app:voronoi-delaunay}

For distinct sites \(\mathbf{c}_1,\ldots,\mathbf{c}_K\), the closed Voronoi cell is
\[
\mathrm{Vor}(\mathbf{c}_k;\mathbf{c})
=\{\mathbf{x}:\|\mathbf{x}-\mathbf{c}_k\|\le\|\mathbf{x}-\mathbf{c}_l\|\ \text{for every }l\}.
\]
A fixed priority rule resolves boundary ties. Nearest-centre assignment
minimizes SSE at fixed centres, while the centering identity above minimizes
it at fixed assignments \citep{lloyd1982least}. If reassignment leaves
all \(K\) cells nonempty, both nonincreasing steps preserve the value of a
global optimum. Nonemptiness is verified explicitly in the proofs of
Theorems~\ref{thm:delaunay-k2} and \ref{thm:qnn-mixture-coverage}.

For planar observations in general position, two sites are joined by a
Delaunay edge exactly when their sample Voronoi cells share an edge,
equivalently when a circle through the two sites has empty interior
\citep{okabe2000spatial,berg2008computational}. These standard duality and
empty-region properties give
\[
\mathrm{MST}\subseteq\mathrm{RNG}\subseteq\mathrm{Gabriel}
\subseteq\mathrm{Delaunay}\subseteq K_n.
\]
The Delaunay graph is a straight-line triangulation with \(O(n)\) edges.
This sparsity does not bound its connected-subset count; the half-plane
argument needed for optimum coverage is proved separately in
Appendix~\ref{app:proof-halfplane-connectivity}.

\section{Planar Delaunay coverage}
\label{supp:delaunay}

\subsection{Planar proximity graphs and Delaunay coverage}
\label{app:delaunay-coverage}

\begin{theorem}[Delaunay coverage for \(K=2\)]
\label{thm:delaunay-k2}
For \(n\ge2\) distinct observations \(\mathbf{X}\subset\mathbb R^2\) in general
position (no three collinear and no four cocircular), at least one global
two-cluster SSE minimizer is realizable in the Delaunay graph. Hence
\(\mathrm{OPT}(\mathrm{Delaunay},2) = \mathrm{OPT}_{\mathrm{unc}}(2)\).
\end{theorem}

The argument first establishes half-plane connectivity for a triangulation,
then applies it to the two nonempty cells of an optimal nearest-centroid
assignment.

\subsubsection{Half-plane connectivity in a triangulation}
\label{app:proof-halfplane-connectivity}

Let \(T\) be a straight-line triangulation of a finite planar point set with
no three collinear vertices, and let \(A\) be an open or closed half-plane.
We prove that \(T[\mathbf{X}\cap A]\) is connected whenever \(\mathbf{X}\cap A\ne\emptyset\).
If \(|\mathbf{X}|\le2\), the triangulation is the single vertex or the edge joining the
two vertices, and the assertion is immediate. Hence assume \(|\mathbf{X}|\ge3\); the
no-three-collinear condition makes \(\operatorname{conv}(\mathbf{X})\) two-dimensional.
Take \(\mathbf{u},\mathbf{v}\in \mathbf{X}\cap A\). If \(\mathbf{u}=\mathbf{v}\) there is nothing to prove. Write the
half-plane as \(A=\{\mathbf{x}:g(\mathbf{x})>0\}\) or \(A=\{\mathbf{x}:g(\mathbf{x})\ge0\}\) for a nonconstant
affine function \(g\).

If \(A\cap\operatorname{conv}(\mathbf{X})\) has empty planar interior, \(A\) must be a
closed supporting half-plane and its intersection with \(\mathbf{X}\) lies on the
supporting line. General position leaves at most two such vertices; when
there are two they are consecutive convex-hull vertices and hence joined by a
triangulation edge. Assume henceforth that the intersection has nonempty
interior. Convexity supplies a polygonal path \(\boldsymbol{\gamma}\) from \(\mathbf{u}\) to \(\mathbf{v}\)
whose interior lies in the interior of
\(A\cap\operatorname{conv}(\mathbf{X})\). Because the triangulation has finitely many
vertices and edges, its interior waypoints may be perturbed within this open
set so that \(\boldsymbol{\gamma}\) avoids every triangulation vertex other than its
endpoints and crosses every encountered edge transversely. This follows by
avoiding the finite union of lines and point conditions that would create a
nongeneric intersection.

The path passes through a finite sequence of adjacent triangles. At a
crossing of a shared edge \([\mathbf{p},\mathbf{q}]\), write
\(\mathbf{y}=(1-t)\mathbf{p}+t\mathbf{q}\), \(0<t<1\). Since \(\mathbf{y}\in A\),
\[
g(\mathbf{y})=(1-t)g(\mathbf{p})+tg(\mathbf{q})
\]
implies that at least one endpoint of the crossed edge belongs to \(A\); in
the open case at least one endpoint has strictly positive \(g\). Choose such
an endpoint \(\mathbf{w}\). Endpoints chosen at two consecutive crossings, together
with \(\mathbf{u}\) at the first triangle and \(\mathbf{v}\) at the last, are vertices of a
common triangle. Any two distinct vertices of a triangle are joined by an
edge of \(T\). The chosen vertices therefore form a walk entirely in
\(T[\mathbf{X}\cap A]\); deleting repetitions gives a path from \(\mathbf{u}\) to \(\mathbf{v}\).
Since \(\mathbf{u},\mathbf{v}\) were arbitrary, the induced graph is connected.

\subsubsection{Proof of Theorem~\ref{thm:delaunay-k2}}
\label{app:proof-delaunay-k2}

For \(n=2\), the only two-block partition consists of two singletons and is
Delaunay-realizable. Assume \(n\ge3\), choose
\(\ell^\star=\{A_1,A_2\}\in\mathcal M_2\), and let \(\boldsymbol{\mu}_1,\boldsymbol{\mu}_2\) be its
centroids. Assign each observation to a nearest one of these centroids using
the fixed rule
\[
A'_1=\{i:\|\mathbf{x}_i-\boldsymbol{\mu}_1\|^2\le\|\mathbf{x}_i-\boldsymbol{\mu}_2\|^2\},
\qquad
A'_2=V\setminus A'_1.
\]
Both sets are nonempty. To see this, let \(\bar{\mathbf{x}}\) be the mean of all
observations. General position implies that some \(\mathbf{x}_j\ne\bar{\mathbf{x}}\). Applying
the ANOVA identity \eqref{eq:anova-merge-app} with
\(A=V\setminus\{j\}\) and \(B=\{j\}\), then solving for
\(\mathrm{SSE}(A)\), gives the valid two-block competitor
\[
\mathrm{SSE}(V\setminus\{j\})+
\mathrm{SSE}(\{j\})
=
\mathrm{OPT}_{\mathrm{unc}}(1)
-\frac{n}{n-1}\|\mathbf{x}_j-\bar{\mathbf{x}}\|^2
<\mathrm{OPT}_{\mathrm{unc}}(1).
\]
Thus \(\mathrm{OPT}_{\mathrm{unc}}(2)<
\mathrm{OPT}_{\mathrm{unc}}(1)\). If, say, \(A'_2=\emptyset\), then
\[
\mathrm{OPT}_{\mathrm{unc}}(1)
\le\sum_{i=1}^n\|\mathbf{x}_i-\boldsymbol{\mu}_1\|^2
\le\mathrm{OPT}_{\mathrm{unc}}(2),
\]
a contradiction; the other case is symmetric.
In particular, \(\boldsymbol{\mu}_1\ne\boldsymbol{\mu}_2\), because equality of the two centroids under
the stated weak/strict tie rule would put every observation in \(A'_1\) and
leave \(A'_2\) empty.

Let \(\boldsymbol{\mu}'_a\) be the mean of \(A'_a\). Nearest-centroid assignment followed
by exact recentering gives
\begin{align*}
\mathrm{SSE}\{A'_1,A'_2\}
&=\sum_{a=1}^2\sum_{i\in A'_a}\|\mathbf{x}_i-\boldsymbol{\mu}'_a\|^2\\
&\le\sum_{a=1}^2\sum_{i\in A'_a}\|\mathbf{x}_i-\boldsymbol{\mu}_a\|^2\\
&\le\sum_{a=1}^2\sum_{i\in A_a}\|\mathbf{x}_i-\boldsymbol{\mu}_a\|^2
=\mathrm{OPT}_{\mathrm{unc}}(2).
\end{align*}
The reverse inequality holds by definition, so
\(\ell'=\{A'_1,A'_2\}\in\mathcal M_2\).

Finally,
\[
\|\mathbf{x}-\boldsymbol{\mu}_1\|^2-\|\mathbf{x}-\boldsymbol{\mu}_2\|^2
=2\mathbf{x}^\mathsf T(\boldsymbol{\mu}_2-\boldsymbol{\mu}_1)+\|\boldsymbol{\mu}_1\|^2-\|\boldsymbol{\mu}_2\|^2
\]
is affine in \(\mathbf{x}\). Hence \(A'_1\) and \(A'_2\) are the intersections of
the data with complementary closed and open half-planes. A Delaunay
triangulation is a straight-line triangulation on \(\mathbf{X}\); applying
the geometric background in Appendix~\ref{app:voronoi-delaunay} and
Appendix~\ref{app:proof-halfplane-connectivity} to both cells gives
\(\ell'\in\mathcal F(\mathrm{Delaunay})\cap\mathcal M_2\).
Equation~\eqref{eq:realizability-criterion} proves the equality of optimum
values.

\section{Scaffold gaps and population risk}
\label{app:scaffold-population-risk}

Suppose \(P\) is supported on \(B(\mathbf0,R)\), and collect the means of
a scaffold-feasible \(K\)-partition \(\widetilde\ell_H\) in
\(\widetilde{\mathbf c}_H\). Convexity places these means in
\(\mathcal C_R\). For any \(\mathbf c^\star\in\mathfrak C_{P,R}\),
nearest-centre assignment gives
\[
\mathcal Q_n(\widetilde{\mathbf c}_H)
\le n^{-1}\mathrm{SSE}(\widetilde\ell_H),\qquad
\mathcal Q_n(\mathbf c^\star)\ge n^{-1}\mathrm{OPT}_{\mathrm{unc}}(K).
\]
For the second inequality, nearest-centre assignment uses at most \(K\)
nonempty groups; recentering and, if necessary, splitting them to obtain
exactly \(K\le n\) groups cannot increase SSE. Hence
\begin{align}
0\le\mathcal Q_P(\widetilde{\mathbf c}_H)-\mathcal Q_P^\star
&\le2Z_{n,R}+\mathcal Q_n(\widetilde{\mathbf c}_H)
-\mathcal Q_n(\mathbf c^\star)\notag\\
&\le2Z_{n,R}+
\frac{\Delta_{\mathrm{opt}}(H,K)+
\Delta_{\mathrm{search}}(\widetilde\ell_H;H,K)}{n}.
\label{eq:scaffold-population-risk}
\end{align}
Exact optimization removes the search term, and coverage removes the
empirical approximation term due to the scaffold. Neither assertion
identifies SSE-optimal groups with latent mixture components.

\newpage
\bibliographystyle{plainnat}
\bibliography{reference}
\end{document}